%% file: main_coap.tex
\RequirePackage{fix-cm}
\documentclass[smallextended]{svjour3}       % onecolumn (second format)
\smartqed  % flush right qed marks, e.g. at end of proof
\usepackage{amsmath}
\usepackage{amsfonts,amssymb}
\usepackage{color}
\usepackage{graphicx}
\usepackage{epstopdf}
\usepackage{algorithm}
\usepackage{algorithmic} % algorithmic
\usepackage[caption=false]{subfig}
\renewcommand{\b}[1]{\ensuremath{\mathbf{#1}}}

\newcommand{\bh}[1]{\ensuremath{\widehat{\mathbf{#1}}}}

\newcommand{\bs}[1]{\ensuremath{{\mathbf{#1}}}}
\newcommand{\bsk}[1]{\ensuremath{{\mathbf{#1}_k}}}
\newcommand{\bskt}[1]{\ensuremath{{\mathbf{#1}^T_k}}}
\newcommand{\bsko}[1]{\ensuremath{{\mathbf{#1}_{k+1}}}}				% J.B., 06/11/19
\newcommand{\bskot}[1]{\ensuremath{{\mathbf{#1}^T_{k+1}}}}				% J.B., 06/11/19

\newcommand{\bko}[1]{\ensuremath{\mathbf{#1}_{k+1}}}
\newcommand{\bkmo}[1]{\ensuremath{\mathbf{#1}_{k-1}}}						% J.B., 06/06/19
\newcommand{\bkoex}[2]{\ensuremath{\mathbf{#1}^{\scriptstyle \text{#2} }_{k+1}}} 	% #2			% J.B., 01/18/19
\newcommand{\bk}[1]{\ensuremath{\mathbf{#1}_k}}
\newcommand{\bkex}[2]{\ensuremath{\mathbf{#1}^{{\scriptstyle \text{#2}}}_{k}}} 					% J.B., 01/18/19
\newcommand{\bzex}[2]{\ensuremath{\mathbf{#1}^{\text{#2}}_{0}}} 					% J.B., 01/18/19
\newcommand{\bhk}[1]{\ensuremath{\widehat{\mathbf{#1}}_k}}
\newcommand{\bhz}[1]{\ensuremath{\widehat{\mathbf{#1}}_0}}

\newcommand{\bhkex}[2]{\ensuremath{\widehat{\mathbf{#1}}^{\text{#2}}_k}} 			% J.B., 06/06/19

\newcommand{\mgap}{\;\;}

\newcommand{\bz}[1]{\ensuremath{\mathbf{#1}_0}}
\newcommand{\biidx}[2]{\ensuremath{{\mathbf{#2}}_{#1}}}
\newcommand{\biidxex}[3]{\ensuremath{{\mathbf{#2}}_{#1}^{\text{#3}}}}				% J.B., 06/06/19
\newcommand{\In}{\ensuremath{\textbf{I}_n}}
\newcommand{\cref}[1]{\eqref{#1}}
\newcommand{\Cref}[1]{\eqref{#1}}
\newcommand{\bfz}{\ensuremath{\mathbf{z}}}
\newcommand{\bfx}{\ensuremath{\mathbf{x}}}
\newcommand{\bfy}{\ensuremath{\mathbf{y}}}
\newcommand{\bfd}{\ensuremath{\mathbf{d}}}
\newcommand{\bfQ}{\ensuremath{\mathbf{Q}}}
\newcommand{\bfepsilon}{\ensuremath{\boldsymbol{\epsilon}}}
\newcommand{\bbC}{\ensuremath{\mathbb{C}}}
\newcommand{\bbR}{\ensuremath{\mathbb{R}}}
\newcommand{\Fcurly}{\ensuremath{\mathcal{F}}}
\newcommand{\diag}[1]{\text{diag}(#1)}
\newcommand{\Refr}[1]{\mathfrak{Re}(#1)}
\renewcommand{\Im}[1]{\mathfrak{Im}(#1)}
\newcommand{\ds}{}
\newcommand{\hf}{}
\newcommand{\jjb}[1]{\textcolor{black}{#1}}
\newcommand{\jjbc}[1]{\textcolor{black}{#1}}
\renewcommand{\Re}{\mathbb{R}}
\journalname{COAP}
\begin{document}

% Adding the cover letter
\input{cover}

%\title{Insert your title here%\thanks{Grants or other notes
%%about the article that should go on the front page should be
%%placed here. General acknowledgments should be placed at the end of the article.}
%}
\title{Compact Representations of Structured BFGS Matrices\thanks{Submitted to the editors DATE.
This material was based upon work supported by the U.S. Department
of Energy, Office of Science, Office of Advanced Scientific Computing Research (ASCR)
under Contract DE-AC02-06CH11347.}}


\author{Johannes J. Brust         \and
	\jjbc{Zichao (Wendy) Di} \and 
        		Sven Leyffer \and
		Cosmin G. Petra  %etc.
}


\institute{J.J. Brust \at
              Argonne National Laboratory, Lemont, IL \\
%              Tel.: +123-45-678910\\
 %             Fax: +123-45-678910\\
              \email{jbrust@anl.gov}           %  \\
           \and
           W. Di \at
           Argonne National Laboratory, Lemont, IL \\
           \email{wendydi@mcs.anl.gov} 
           \and
           S. Leyffer \at
              Argonne National Laboratory, Lemont, IL \\
              \email{leyffer@mcs.anl.gov}
            \and
            C. G. Petra \at
            Lawrence Livermore National Laboratory, Livermore, CA \\
            \email{petra1@llnl.gov}            
}

\date{Received: date / Accepted: date}

\maketitle

\begin{abstract}
For general large-scale optimization problems compact representations
exist in which recursive quasi-Newton update formulas are represented as compact
matrix factorizations. For problems in which the objective function contains
additional structure, so-called structured quasi-Newton methods exploit
available second-derivative information and approximate unavailable 
second derivatives. This article develops the compact representations
of two structured Broyden-Fletcher-Goldfarb-Shanno update formulas.
The compact representations enable efficient limited memory and initialization strategies.
Two limited memory line search algorithms are described and tested on a collection
of problems, \jjbc{including a real world large scale imaging application}.
%\subclass{AMS 65K05 \and AMS 65F30 \and AMS 90C53 \and AMS 90C06}
\end{abstract}

%------------------------------------- Introduction ----------------------------------------------------%
\section{Introduction}
\label{sec:intro}

The unconstrained minimization problem is 
\begin{equation}
	\label{eq:intro_min}
	\underset{ \b{x} \in \Re^n }{ \text{ minimize } } f(\b{x}),
\end{equation}
where $ f:\Re^n \to \Re $ is assumed to be twice continuously differentiable.
If the Hessian matrix $ \nabla^2 f(\b{x}) \in \Re^{n \times n} $ is unavailable,
because it is unknown or difficult to compute, then quasi-Newton approaches are
effective methods, which approximate properties of the Hessian at each iteration,
 $ \nabla^2 f(\bko{x}) \approx \bko{B} $ \cite{DM77}. Arguably, the most
widely used quasi-Newton matrix is the Broyden-Fletcher-Goldfarb-Shanno (BFGS)
matrix \cite{Broyden70,Fletcher70,Goldfarb70,Shanno70}, because of its 
desirable results on many problems. Given $ \bk{s} \equiv \bko{x} - \bk{x} $ 
and $ \bk{y} \equiv \nabla f(\bko{x}) - \nabla f(\bk{x})  $ the BFGS recursive update 
formula is
\begin{equation}
	\label{eq:intro_bfgs}
	\bko{B} = \bk{B} - \frac{1}{\bk{s}^T \bk{B} \bk{s}} \bk{B}\bk{s} \bk{s}^T\bk{B} +
			\frac{1}{\bk{s}^T\bk{y}} \bk{y} \bk{y}^T.
\end{equation}
For a symmetric positive definite initialization $ \bz{B} \in \Re^{n \times n} $ 
\cref{eq:intro_bfgs} generates symmetric positive definite matrices
as long as $ \bk{s}^T \bk{y} > 0 $ for all $ k \ge 0 $ (see \cite[Section 2]{Fletcher70}).

\subsection{BFGS Compact Representation}
Byrd et al. \cite{ByrNS94} propose the compact representation of the recursive formula 
\cref{eq:intro_bfgs}. The compact representation has been successfully
used for large-scale unconstrained and constrained optimization \cite{ZhuByrdNocedal97}. Let the sequence of pairs $ \left \{ \biidx{i}{s}, \biidx{i}{y} \right \}_{i=0}^{k-1} $
be given, and let these vectors be collected in the matrices
$ \bk{S} = \left[ \bz{s}, \quad \cdots, \quad \biidx{k-1}{s} \right] \in \Re^{n \times k} $ and 
$ \bk{Y} = \left[ \bz{y}, \quad \cdots, \quad \biidx{k-1}{y} \right] \in \Re^{n \times k} $. Moreover, let
$ \bk{S}^T \bk{Y} = \bk{L} + \bk{R} $, where $ \bk{L} \in \Re^{k \times k} $
is the strictly lower triangular matrix, $ \bk{R} \in \Re^{k \times k} $ is the upper triangular matrix (including the diagonal),
and $ \jjbc{\bk{D} = \text{diag}(\biidx{0}{s}^T \biidx{0}{y},\cdots,\biidx{k-1}{s}^T \biidx{k-1}{y})} \in \Re^{k \times k} $ is the diagonal part of $ \bk{S}^T \bk{Y} $. The %\text{diag}(\bk{S}^T \bk{Y})
compact representation of the BFGS formula \cref{eq:intro_bfgs} is \cite[Theorem 2.3]{ByrNS94}:
\begin{equation}
	\label{eq:intro_cbfgs}
	\bk{B} = \bz{B} - \left[ \: \bz{B}\bk{S} \quad \bk{Y} \: \right]
	{\begin{bmatrix}
		\bk{S}^T \bz{B} \bk{S} 	& \bk{L} \\
		\bk{L}^T 				& - \bk{D}
	\end{bmatrix}}^{-1}
	\begin{bmatrix}
		\bk{S}^T \bz{B} \\
		\bk{Y}^T
	\end{bmatrix}.
\end{equation}
For large optimization problems limited memory versions of the compact representation in \cref{eq:intro_cbfgs}
are used. The limited memory versions typically store only the last $ m > 0 $ pairs
 $ \left \{ \biidx{i}{s}, \biidx{i}{y} \right \}_{i=k-m}^{k-1} $ when $ k \ge m$. 
In limited memory BFGS (L-BFGS) the dimensions of $ \bk{S} $ and $ \bk{Y} $
are consequently $ n \times m $. Usually the memory parameter is much smaller than the problem size,
namely, $ m \ll n $. A typical range for this parameter is $ 5 \le m \le 50 $ (see Boggs and Byrd in \cite{BoggsByrd19}). 
Moreover, in line search L-BFGS methods the initialization is frequently chosen as $ \bz{B} = \hat{\sigma}_k \In $, where 
$ \hat{\sigma}_k = \bkmo{y}^T \bkmo{y} / \bkmo{s}^T\bkmo{y} $. Such an initialization enables efficient
computations with the formula in \cref{eq:intro_cbfgs}, and adds extra information 
through the parameter $ \hat{\sigma}_k $, which depends on the iteration $ k $.

\subsection{Structured Problems}
When additional information about the structure of the objective function is 
known, it is desirable to include this information in a quasi-Newton update. Initial
research efforts on structured quasi-Newton methods were in the context of nonlinear least squares problems. These include the work of Gill and Murray \cite{GillMurray78},
Dennis et al. \cite{DennisGayWalsh81,DennisMartinezTapia89}, and Yabe and Takahashi \cite{YabeTakahashi91}. Recently, Petra et al. \cite{PetraChiangAnitescu19} formulated the general structured minimization problem as
\begin{equation}
	\label{eq:intro_strucmin}
	\underset{ \b{x} \in \Re^n }{ \text{ minimize } } f(\b{x}), \quad \quad 
	f(\b{x}) = \widehat{k}(\b{x}) + \widehat{u}(\b{x}), % \mathfrac from AMSfonts
\end{equation}
where $ \widehat{k}: \Re^n \to \Re $ has known gradients and known Hessians
and $ \widehat{u}: \Re^n \to \Re $ has known gradients but unknown 
Hessians. \jjb{For instance, objective functions composed of a general nonlinear function plus a regularizer or penalty term are described with \cref{eq:intro_strucmin}. 
Thus, applications such as regularized logistic regressions \cite{TeoVishwanthanSmolaLe10} or optimal control problems contain
structure that may be exploited, when we assume that the Hessian of the regularizer is known.}
\jjbc{We note that nonlinear least squares problems typically do not have the form as in \eqref{eq:intro_strucmin},
yet available second derivatives may also be used for this class of problems
after reformulating the quasi-Newton vectors. We will describe an image reconstruction application in the numerical experiments, Section 4.} 
Even though approximating the Hessian of the objective function 
in \cref{eq:intro_strucmin} by formula \cref{eq:intro_bfgs} or \cref{eq:intro_cbfgs} 
is possible, this would not exploit the known parts of the Hessian. 
Therefore in \cite{PetraChiangAnitescu19} structured BFGS (S-BFGS) 
updates are derived, which combine known Hessian information with
BFGS approximations for the unknown Hessian components. At each iteration the Hessian of the objective
 is approximated as $ \nabla^2 f(\bko{x}) \approx \nabla^2 \widehat{k}(\bko{x}) + \bko{A} $,
where $ \bko{A} $ approximates the unknown Hessian, that is, $ \bko{A} \approx \nabla^2 \widehat{u}(\bko{x}) $.
Given the known Hessian $ \nabla^2 \widehat{k}(\bko{x}) \equiv \bko{K} $ and the gradients of $ \widehat{u} $,
%$ \bhk{u} \equiv \nabla \widehat{u}(\bko{x})-\nabla \widehat{u}(\bk{x}) $ 
let $  \bk{u} \equiv  \bko{K}\bk{s} +(\nabla \widehat{u}(\bko{x})-\nabla \widehat{u}(\bk{x})) $.
One of two structured approximations from \cite{PetraChiangAnitescu19} is the 
structured BFGS-Minus (S-BFGS-M) update
\begin{equation}
	\label{eq:intro_sbfgsm}
	\bkoex{A}{M} = \bkex{B}{M} - \bko{K} -  \frac{1}{ \bk{s}^T\bkex{B}{M} \bk{s} }
	\bkex{B}{M}\bk{s} \bk{s}^T \bkex{B}{M} + 
	\frac{1}{\bk{s}^T\bk{u}}\bk{u}\bk{u}^T,
\end{equation}
%\begin{equation}
%	\label{eq:intro_sbfgsm}
%	\bkoex{A}{M} = \bkex{A}{M} + \bk{K} - \bko{K} -  \frac{1}{ \bk{s}^T( \bkex{A}{M} + \bk{K} ) \bk{s} }
%	( \bkex{A}{M} + \bk{K} )\bk{s} \bk{s}^T ( \bkex{A}{M} + \bk{K} )^T + 
%	\frac{(\bko{K}\bk{s} + \bhk{u})(\bko{K}\bk{s}+\bhk{u})^T}{\bk{s}^T\bhk{u}},
%\end{equation}
where $ \bkex{B}{M} = \bkex{A}{M} + \bk{K} $. By adding $\bko{K}$ to both sides, the update from \cref{eq:intro_sbfgsm} implies a 
formula for $ \bkoex{B}{M} $ that resembles \cref{eq:intro_bfgs}, in which $\bko{B}$, $\bk{B}$, and $ \bk{y} $ 
are replaced by $\bkoex{B}{M}$, $ \bkex{B}{M} $, and $ \bk{u} $, respectively. Consequently, $ \bkoex{B}{M} $
is symmetric positive definite given a symmetric positive definite initialization $ \b{B}_0^{\text{M}} $
as long as $ \bk{s}^T\bk{u} > 0 $ for $k \ge 0$.
A second formula is the structured BFGS-Plus (S-BFGS-P) update 
\begin{equation}
	\label{eq:intro_sbfgsp}
	\bkoex{A}{P} = \bkex{A}{P} - \frac{1}{\bk{s}^T\bhkex{B}{P}\bk{s}}
	 \bhkex{B}{P}\bk{s}\bk{s}^T\bhkex{B}{P}
				 + \frac{1}{\bk{s}^T \bk{u}}\bk{u}\bk{u}^T,
\end{equation}
where $ \bhkex{B}{P} = \bkex{A}{P} + \bko{K} $. \jjbc{After adding $ \biidx{k+1}{K} $
to both sides in \eqref{eq:intro_sbfgsp}, the left hand side ($\bkoex{B}{P} = \bkoex{A}{P} + \biidx{k+1}{K} $)
is positive definite if $ \bhkex{B}{P} $ is positive definite and $ \bk{s}^T\bk{u} > 0 $. However, in general,
$\bhkex{B}{P}$ does not have to be positive definite, because the known Hessian, $ \biidx{k+1}{K} $,
may not be positive definite. Similar to \eqref{eq:intro_bfgs} the update in \eqref{eq:intro_sbfgsp} 
is also rank 2.}
Both of the updates in eqs. \eqref{eq:intro_sbfgsm} and \eqref{eq:intro_sbfgsp} were implemented in a line search algorithm and 
compared with the BFGS formula \eqref{eq:intro_bfgs} in \cite{PetraChiangAnitescu19}. The structured %  unstructured
updates obtained better results in terms of iteration count and function evaluations than did the
unstructured counterparts. Unlike the BFGS formula from \cref{eq:intro_bfgs}, which 
recursively defines $ \bko{B} $ as a rank-2 update to $ \bk{B} $, the formulas 
for $\bkoex{A}{M}$ and $\bkoex{A}{P}$ in eqs. \cref{eq:intro_sbfgsm}, \eqref{eq:intro_sbfgsp} additionally depend
 on the known Hessians $ \bko{K} $ and $ \bk{K} $.
For this reason the compact representations of $\bkoex{A}{M}$ and $\bkoex{A}{P}$ are 
different from the one for $ \bko{B} $ in \cref{eq:intro_cbfgs} and have not yet been developed.
\jjbc{The updates in \eqref{eq:intro_sbfgsm} and \eqref{eq:intro_sbfgsp} are dense 
in general, and hence neither are suitable for large-scale optimization. Hence, here,
we develop first a compact representation of \eqref{eq:intro_sbfgsm} and 
\eqref{eq:intro_sbfgsp} and then show how to exploit them to develop structured
limited-memory quasi-Newton updates.}

%The inverse of \label{eq:intro_cbfgs} may be obtained by the Sherman-Morrison-Woodbury formula

\subsection{Article Contributions}
In this article we develop the compact representations of the structured BFGS updates $ \bkoex{A}{M} $ and
 $ \bkoex{A}{P} $ from eqs. \cref{eq:intro_sbfgsm} and \eqref{eq:intro_sbfgsp} 
 \jjbc{that lead to practical large-scale limited-memory implementations}. \jjbc{Unwinding the update formula in 
 \eqref{eq:intro_sbfgsp} is challenging, however by using an induction technique we are able to derive the explicit expression
 of the compact representation.} 
 We propose the limited memory versions of the 
compact structured BFGS (L-S-BFGS) matrices and describe line search algorithms 
\jjbc{(with slightly modified Wolfe conditions)} that implement
them. We exploit the compact representations in order to compute search directions by means of the
Sherman-Morrison-Woodbury formula and implement effective initialization strategies. 
Numerical experiments of the proposed L-S-BFGS methods on various problems are presented.

%---------------------------------------- Compact Representations ----------------------------------------------------%
\section{Compact Representations of Structured BFGS Updates}
\label{sec:compact}
To develop the compact representations of the structured BFGS formulas,
we define
\begin{equation}
	\label{eq:comp_usu}
	\bk{U} = \left[ \: \bz{u}, \quad \cdots ,\quad \bkmo{u} \: \right], \quad 
	\bk{S}^T\bk{U} = \bkex{L}{U} + \bk{R}^{\text{U}}, \quad
	\text{diag}(\bk{S}^T\bk{U}) = \bkex{D}{U},
\end{equation}
where $ \bk{U} \in \Re^{n \times k} $ collects all $ \bk{u} $ for $ k \ge 0 $ and
where $ \bkex{L}{U} \in \Re^{k \times k} $ is a strictly lower triangular matrix, 
$ \bk{R}^{\text{U}} \in \Re^{k \times k} $ is an upper triangular matrix (including the diagonal),
and $ \bkex{D}{U} \in \Re^{k \times k} $ is the diagonal part of $ \bk{S}^T \bk{U} $.

\subsection{Compact Representation of $ \bkex{A}{M} $}
Theorem \ref{thrm:comp_csbfgsm} contains the compact representation of
$ \bkex{A}{M} $.

\begin{theorem}
	\label{thrm:comp_csbfgsm}
	The compact representation of $ \bkex{A}{M} $ in the update formula
	\cref{eq:intro_sbfgsm} is 
	\begin{equation}
		\label{eq:comp_csbfgsm}
		\bkex{A}{M} = \biidxex{0}{B}{M} - \bk{K} - \left[ \: \biidxex{0}{B}{M}\bk{S} \quad \bk{U} \: \right]
		{\begin{bmatrix}
			\bk{S}^T \biidxex{0}{B}{M} \bk{S} 	& \bkex{L}{U} \\
			(\bkex{L}{U})^T 				& -\bkex{D}{U}
		\end{bmatrix}}^{-1}
		\begin{bmatrix}
			\bk{S}^T(\biidxex{0}{B}{M})^T \\
			\bk{U}^T
		\end{bmatrix},
	\end{equation}
	where $ \bk{S} $ is as defined in \cref{eq:intro_cbfgs}, $ \bk{U} $, $ \bkex{L}{U} $, and 
	$ \bkex{D}{U} $ are defined in \cref{eq:comp_usu}, and 
	%\biidxex{0}{B}{M} = \biidxex{0}{A}{M} + \bz{K}.
	\begin{equation*}
		\biidxex{0}{B}{M} = \biidxex{0}{A}{M} + \bz{K}.
	\end{equation*}
\end{theorem}
\begin{proof}\hspace{-0.15cm}.
Observe that by adding $ \bko{K} $ to both sides of \cref{eq:intro_sbfgsm} the update formula
of $ \bkoex{B}{M} $ becomes
\begin{equation*}
	\bkoex{B}{M} = \bkex{B}{M} - \frac{1}{\bk{s}^T \bkex{B}{M} \bk{s}} \bkex{B}{M}\bk{s} \bk{s}^T\bkex{B}{M} +
			\frac{1}{\bk{s}^T\bk{u}} \bk{u} \bk{u}^T.
\end{equation*}
This expression is the same as \cref{eq:intro_bfgs} when $ \bkoex{B}{M} $ is relabeled as $ \bko{B} $, 
$ \bkex{B}{M} $ is relabeled as $ \bk{B} $, and $ \bk{u} $ is relabeled as $ \bk{y} $. The compact representation
of \cref{eq:intro_bfgs} is given by \cref{eq:intro_cbfgs}, and therefore the compact representation of 
$ \bkex{B}{M} $ is given by \cref{eq:intro_cbfgs} with $ \bk{Y} $ replaced by $ \bk{U} $ and $ \bz{B} $
replaced by $ \biidxex{0}{B}{M} $. Then \cref{eq:comp_csbfgsm} is obtained by subtracting $ \bk{K}  $
from the compact representation of $ \bkex{B}{M} $, and noting that 
$ \biidxex{0}{B}{M} = \bzex{A}{M} + \bz{K} $. Since 
\biidxex{k}{B}{M} is symmetric positive definite as long as $ \biidxex{0}{B}{M} $
is symmetric positive definite and $ \bk{s}^T\bk{u} > 0  $ for $ k \ge 0 $, the inverse in the 
right-hand side of \cref{eq:comp_csbfgsm} is nonsingular as long as $ \biidxex{0}{B}{M} $
is symmetric positive definite and $ \bk{s}^T\bk{u} > 0  $ for $ k \ge 0$. \qed
\end{proof}
Corollary \ref{cor:comp_icsbfgsm} describes the compact representation of the inverse 
$ \bkex{H}{M} = ( \bk{K} +  \bkex{A}{M}  )^{-1} $, which is used to compute search directions
in a line search algorithm (e.g., $ \bkex{p}{M} = - \bkex{H}{M} \nabla f(\bk{x})) $.
\begin{corollary}
	\label{cor:comp_icsbfgsm}
	The inverse $ \bkex{H}{M} = \left( \bk{K} + \bkex{A}{M} \right)^{-1} $,
%	\begin{equation*}
%		\bkex{H}{M} = \left( \bk{K} + \bkex{A}{M} \right)^{-1},
%	\end{equation*}
	with the compact representation of \bkex{A}{M} from \cref{eq:comp_csbfgsm}, is given as
	\begin{align}
		\label{eq:comp_icsbfgsm}
		&\bkex{H}{M} = \biidxex{0}{H}{M} + \nonumber \\
		&\left[ \bk{S} \: \: \biidxex{0}{H}{M} \bk{U} \right] % \biidxex{0}{H}{M} +
			\begin{bmatrix}
				(\bkex{T}{U})^T\left( \bkex{D}{U} + \bk{U}^T \biidxex{0}{H}{M} \bk{U} \right) \bkex{T}{U} 	& \text{-}(\bkex{T}{U})^{T} \\
				\text{-}\bkex{T}{U} 																& \b{0}_{k \times k}
			\end{bmatrix}
			\begin{bmatrix}
				\bk{S}^T \\
				\bk{U}^T (\biidxex{0}{H}{M})^{T}
			\end{bmatrix},
	\end{align} 
%	\begin{equation}
%		\label{eq:comp_icsbfgsm}
%		\bkex{H}{M} = \biidxex{0}{H}{M} + \left[ \bk{S} \quad \biidxex{0}{H}{M} \bk{U} \right]
%			\begin{bmatrix}
%				(\bkex{T}{U})^T\left( \bkex{D}{U} + \bk{U}^T \biidxex{0}{H}{M} \bk{U} \right) \bkex{T}{U} 	& -(\bkex{T}{U})^{T} \\
%				-\bkex{T}{U} 																& \b{0}_{k \times k}
%			\end{bmatrix}
%			\begin{bmatrix}
%				\bk{S}^T \\
%				\bk{U}^T (\biidxex{0}{H}{M})^{T}
%			\end{bmatrix},
%	\end{equation} 
	where 
	\begin{equation*}
		\biidxex{0}{H}{M} = (\biidxex{0}{B}{M})^{-1} = (\biidxex{0}{A}{M} + \bz{K})^{-1},
	\end{equation*}	
	and where $ \bkex{T}{U} = (\bk{R}^U)^{-1} $ with $ \bk{S}$, $ \bk{U}$, $\bkex{D}{U} $, 
	and $ \bk{R}^U  $ defined in Theorem \ref{thrm:comp_csbfgsm} and \eqref{eq:comp_usu}. 
\end{corollary}
\begin{proof}\hspace{-0.15cm}.
	Define 
	\begin{equation*}
	 \bsk{\Xi} \equiv \left[ \: \biidxex{0}{B}{M} \bk{S} \quad \bk{U} \: \right], \quad \quad
		\bk{M} \equiv 
			\begin{bmatrix}
				\bk{S}^T \biidxex{0}{B}{M} \bk{S} 	& \bkex{L}{U} \\
				(\bkex{L}{U})^T 				& -\bkex{D}{U}
			\end{bmatrix},
	\end{equation*}
	for the compact representation of $ \bkex{A}{M} $ in \cref{eq:comp_csbfgsm}. 
	Let $ \biidxex{0}{H}{M} = (\biidxex{0}{B}{M})^{-1} $ then the expression of
	$ \bkex{H}{M} $ is obtained by the Sherman-Morrison-Woodbury identity:
	\begin{align*}
		\bkex{H}{M} &= \left( \bk{K} + \bkex{A}{M} \right)^{-1}	 \\
		&= \left( \biidxex{0}{B}{M} - \bsk{\Xi} \bk{M}^{-1} \bs{\Xi}^T_k \right)^{-1} \\
		&= (\biidxex{0}{B}{M})^{-1} + 
		(\biidxex{0}{B}{M})^{-1} \bsk{\Xi}\left[ \bk{M} - \bs{\Xi}^T_k (\biidxex{0}{B}{M})^{-1} \bsk{\Xi} \right]^{-1} \bs{\Xi}^T_k (\biidxex{0}{B}{M})^{-1} \\
		&= \biidxex{0}{H}{M} - 
		\biidxex{0}{H}{M} \bsk{\Xi} 
		{\begin{bmatrix}
			\b{0}_{k \times k} 	& \bkex{R}{U} \\
			(\bkex{R}{U})^T 			& \bkex{D}{U} + \bk{U}^T \biidxex{0}{H}{M} \bk{U}
		\end{bmatrix}}^{-1}
		\bs{\Xi}^T_k \biidxex{0}{H}{M} \\
		&=\biidxex{0}{H}{M} + 
		\biidxex{0}{H}{M} \bsk{\Xi} 
		{\begin{bmatrix}
			(\bkex{R}{U})^{-T} \left( \bkex{D}{U} + \bk{U}^T \biidxex{0}{H}{M} \bk{U} \right) (\bkex{R}{U})^{-1} 	& -(\bkex{R}{U})^{-T} \\
			-(\bkex{R}{U})^{-1} 															& \b{0}_{k \times k}
		\end{bmatrix}}
		\bs{\Xi}^T_k \biidxex{0}{H}{M}
%		&= \biidxex{0}{H}{M} - 
%		\biidxex{0}{H}{M} \bsk{\Xi} 
%		{\begin{bmatrix}
%			\bk{R}^{-T} \left( \bkex{D}{U} + \bk{U}^T \biidxex{0}{H}{M} \bk{U} \right) \bk{R}^{-1} 	& -\bk{R}^{-T} \\
%			-\bk{R}^{-1} 															& \b{0}_{k \times k}
%		\end{bmatrix}}
%		\bs{\Xi}^T_k \biidxex{0}{H}{M},
	\end{align*}
	where the third equality is obtained from applying the Sherman-Morrison-Woodbury inverse, the fourth equality
	uses the identity $ \bk{S}^T\bk{U} - \bkex{L}{U} = \bkex{R}{U} $, and the fifth equality is obtained 
	% by defining $ (\bsiidxex{0}{\Psi}{M})^{-1} \equiv \biidxex{0}{H}{M} $ and
	 by explicitly computing the inverse of the block matrix. Using $ (\bkex{R}{U})^{-1} = \bkex{T}{U} $ and 
	$ (\biidxex{0}{B}{M})^{-1} \bsk{\Xi} = (\biidxex{0}{B}{M})^{-1}[ \biidxex{0}{B}{M}\bk{S} \quad \bk{U} ] $ 
	yields the expression in \cref{eq:comp_icsbfgsm}. \qed
\end{proof}

\subsection{Compact Representation of $ \bkex{A}{P} $}
\jjbc{Developing the compact representation of \eqref{eq:intro_sbfgsp} is
more challenging and requires an inductive argument.}
Specifically, % To develop the compact representation of the structured BFGS matrix $ \bkex{A}{P} $
we define $ \bk{v} \equiv \bko{K} \bk{s} $ in addition to the expressions in \cref{eq:comp_usu} and
\begin{equation}
	\label{eq:comp_vsv}
	\bk{V} = \left[ \: \bz{v}, \quad \ldots ,\quad \bkmo{v} \: \right], \quad 
	\bk{S}^T\bk{V} = \bkex{L}{V} + \bk{R}^{\text{V}}, \quad \text{diag}(\bk{S}^T\bk{V}) = \bkex{D}{V},
\end{equation}
where $ \bk{V} \in \Re^{n \times k} $ collects all $ \bk{v} $ for $ k \ge 0 $ and
where $ \bkex{L}{V} \in \Re^{k \times k} $ is the strictly lower triangular matrix, 
$ \bk{R}^{\text{V}} \in \Re^{k \times k} $ is the upper triangular matrix (including the diagonal),
and $ \bkex{D}{V} \in \Re^{k \times k} $ is the diagonal part of $ \bk{S}^T \bk{V} $.
Theorem \ref{thrm:comp_csbfgsp} contains the compact representation of $ \bkex{A}{P} $.
%and it is proved by an induction.

\begin{theorem}
	\label{thrm:comp_csbfgsp}
	The compact representation of $ \bkex{A}{P} $ in the update formula
	\cref{eq:intro_sbfgsp} is 
	\begin{equation}
		\label{eq:comp_csbfgsp}
		\bkex{A}{P} = \biidxex{0}{A}{P} - \left[ \: \bk{Q} \quad \bk{U} \: \right]
		{\begin{bmatrix}
			\bkex{D}{V} + \bkex{L}{V} + (\bkex{L}{V})^T + \bk{S}^T \biidxex{0}{A}{P} \bk{S} 	& \bkex{L}{U} \\
			(\bkex{L}{U})^T 												& -\bkex{D}{U}
		\end{bmatrix}}^{-1}
		\begin{bmatrix}
			\bk{Q}^T \\
			\bk{U}^T
		\end{bmatrix},
	\end{equation}
	where 
	\begin{equation*}
		\bk{Q} \equiv \bk{V} + \biidxex{0}{A}{P} \bk{S}, 
	\end{equation*}
	and where $ \bk{S}, \bk{U}, \bkex{D}{U}$, and $\bkex{L}{U} $ are defined in \cref{thrm:comp_csbfgsp} and $ \bk{V} $, $ \bkex{L}{V} $, and
	$ \bkex{D}{V} $ are defined in \cref{eq:comp_vsv}.
\end{theorem}
\begin{proof}\hspace{-0.15cm}.
The proof of \cref{eq:comp_csbfgsp} is by induction. For $ k = 1 $ it follows that % in \cref{eq:comp_csbfgsp}
\begin{align*}
	\biidxex{1}{A}{P} &= \biidxex{0}{A}{P} - \left[ \bz{v} + \biidxex{0}{A}{P}\bz{s} \quad \bz{u} \right]
	{\begin{bmatrix}
		\bz{s}^T\bz{v} + \bz{s}^T\biidxex{0}{A}{P}\bz{s} 		& \\
												& -\bz{s}^T\bz{u}
	\end{bmatrix}}^{-1}
	\begin{bmatrix}
		(\bz{v} + \biidxex{0}{A}{P}\bz{s})^T \\
		\bz{u}^T
	\end{bmatrix} \\
	&= \biidxex{0}{A}{P} - \frac{1}{ \bz{s}^T( \biidx{1}{K} + \biidxex{0}{A}{P} )\bz{s} }( \biidx{1}{K} + \biidxex{0}{A}{P} )\bz{s} \bz{s}^T( \biidx{1}{K} + \biidxex{0}{A}{P} )^T
	+ \frac{1}{\bz{s}^T\bz{u}}\bz{u}\bz{u}^T,
\end{align*}
\jjbc{which shows that \eqref{eq:comp_csbfgsp} holds for $k=1$}.
This expression is the same as $ \biidxex{1}{A}{P} $ in \cref{eq:intro_sbfgsp}, and thus the compact representation holds for $ k =1$. 
% J.B., 06/07/19. The next part of the proof is taken from 'report_compact_SQN.tex'.
% This requires renaming some quantities.
Next assume that \cref{eq:comp_csbfgsp} is valid for $ k \ge 1 $, and in particular let it be represented as
\begin{equation}
	\label{eq:comp_csbfgsp_induc}
	\bkex{A}{P} = \bzex{A}{P} - \left[ \mgap \bk{Q}  \mgap \bk{U} \mgap \right]
				\left[
					\begin{array}{c c}
						(\bk{M})_{11} 		& (\bk{M})_{12} \\
						(\bk{M})^T_{12} 	& (\bk{M})_{22}
					\end{array}
				\right]^{-1}
				\left[
					\begin{array}{c}
						\bk{Q}^T \\
						\bk{U}^T
					\end{array}
				\right],
\end{equation}
where 
\begin{equation*}
	\mgap (\bk{M})_{11} = \bkex{D}{V} + \bkex{L}{V} + ( \bkex{L}{V} )^T + \bk{S}^T \bzex{A}{P} \bk{S}, 
	\mgap (\bk{M})_{12} = \bkex{L}{U}, 
	\mgap (\bk{M})_{22} = -\bkex{D}{U}.
\end{equation*}
We verify the validity of \cref{eq:comp_csbfgsp_induc} by substituting it in the update
formula \cref{eq:intro_sbfgsp}, and then seek the representation \jjbc{\eqref{eq:comp_csbfgsp_induc} for $k+1$:}
\begin{equation*}
	\bkoex{A}{P} = \bzex{A}{P} - \left[ \mgap \bko{Q}  \mgap \bko{U} \mgap \right]
				\left[
					\begin{array}{c c}
						(\bko{M})_{11} 		& (\bko{M})_{12} \\
						(\bko{M})^T_{12} 	& (\bko{M})_{22}
					\end{array}
				\right]^{-1}
				\left[
					\begin{array}{c}
						\bko{Q}^T \\
						\bko{U}^T
					\end{array}
				\right].
\end{equation*}
First let
\begin{equation*}
	\bk{q} =  \bk{v} + \bzex{A}{P}\bk{s}, \quad \bk{w} = \bk{Q}^T\bk{s}, \quad \bk{r} = \bk{U}^T \bk{s}, \quad \bs{\xi}_k = \left[
									\begin{array}{c}
										\bk{w}	\\
										\bk{r} 
									\end{array}
								\right],
\end{equation*}
and note that in \cref{eq:intro_sbfgsp} it holds that
\begin{align*}
	(\bkex{A}{P}+\bko{K})\bk{s} 	&= \bkex{A}{P}\bk{s}+\bk{v} \\
							&= \bzex{A}{P} \bk{s} - \left[ \mgap \bk{Q} \mgap \bk{U} \mgap \right] [\bk{M}]^{-1}
								\left[
									\begin{array}{c}
										\bk{Q}^T \bk{s} \\
										\bk{U}^T \bk{s} 
									\end{array}
								\right]
								+ \bk{v} \\
							&\equiv \bk{q} - \left[ \mgap \bk{Q} \mgap \bk{U} \mgap \right] [\bk{M}]^{-1}
								\left[
									\begin{array}{c}
										\bk{w}	\\
										\bk{r} 
									\end{array}
								\right] \\
							&\equiv \bk{q} - \left[ \mgap \bk{Q} \mgap \bk{U} \mgap \right] [\bk{M}]^{-1} \bs{\xi}_k.							
\end{align*}
Next we define $ \sigma_k^{\text{P}} = 1 / \bk{s}^T( \bkex{A}{P} + \bko{K} ) \bk{s} $ and obtain the following representation of $ \bkoex{A}{P} $
\jjbc{using \eqref{eq:intro_sbfgsp} and \eqref{eq:comp_csbfgsp_induc}:}
\begin{align*}
	\bkoex{A}{P} 	&= \bkex{A}{P} - \sigma_k^{\text{P}} ( \bkex{A}{P} \bk{s} + \bk{v} ) ( \bkex{A}{P} \bk{s} + \bk{v} )^T + \frac{1}{\bk{s}^T\bk{u}} \bk{u} \bk{u}^T \\
				&= \bzex{A}{P} - \sigma_k^{\text{P}}\left[ \mgap \bk{Q} \mgap \bk{U} \mgap \bk{q} \mgap \right]
					\left[
						\begin{array}{c c }
							\frac{\bk{M}^{-1}}{\sigma_k^{\text{P}}} +  \bk{M}^{-1} % (1/\sigma_k^{\text{P}})
								\bs{\xi}_k \bs{\xi}^T_k
								\bk{M}^{-1} & - \bk{M}^{-1} \bs{\xi}_k \\
								- \bs{\xi}^T_k \bk{M}^{-1} &
								1
						\end{array}
					\right]
					\left[
						\begin{array}{c}
							\bk{Q}^T \\
							\bk{U}^T \\
							\bk{q}^T
						\end{array}
					\right] \\
				&\quad + \frac{1}{\bk{s}^T\bk{u}} \bk{u} \bk{u}^T \\
				&= \bzex{A}{P} - \left[ \mgap \bk{Q} \mgap \bk{U} \mgap \bk{q} \mgap \right]
				\left[
						\begin{array}{c c }
								\bk{M} & 
								\left[
									\begin{array}{c}
										\bk{w} \\
										\bk{r}
									\end{array} 
								\right] \\
								\left[ \mgap \bk{w}^T \mgap \bk{r}^T \mgap \right]  & \bk{s}^T \bk{q}
						\end{array}
				\right]^{-1}
				\left[
						\begin{array}{c}
							\bk{Q}^T \\
							\bk{U}^T \\
							\bk{q}^T
						\end{array}
					\right]
					+ \frac{1}{\bk{s}^T\bk{u}} \bk{u} \bk{u}^T.
\end{align*} 
Using the permutation matrix $ \b{P} = \left[ \mgap \biidx{1}{e} \mgap \cdots \mgap \biidx{k}{e} \mgap \biidx{2k+1}{e} \mgap \cdots \mgap \biidx{2k}{e} \mgap \right] $,
we represent $ \bkoex{A}{P} $ as 
\begin{align*}
	\bkoex{A}{P} &= \bzex{A}{P} - \left[ \mgap \bk{Q} \mgap \bk{U} \mgap \bk{q} \mgap \right] \b{P} \b{P}^T
				\left[
						\begin{array}{c c }
								\bk{M} & 
								\left[
									\begin{array}{c}
										\bk{w} \\
										\bk{r}
									\end{array} 
								\right] \\
								\left[ \mgap \bk{w}^T \mgap \bk{r}^T \mgap \right]  & \bk{s}^T \bk{q}
						\end{array}
				\right]^{-1}
				\b{P} \b{P}^T
				\left[
						\begin{array}{c}
							\bk{Q}^T \\
							\bk{U}^T \\
							\bk{q}^T
						\end{array}
					\right] \\
			&\quad		+ \frac{1}{\bk{s}^T\bk{u}} \bk{u} \bk{u}^T \\
			&= \bzex{A}{P} - \left[ \mgap \bk{Q} \mgap \bk{q}  \mgap \bk{U} \mgap \bk{u} \mgap \right] 
				\left[
						\begin{array}{c c c c }
								(\bk{M})_{11} 		& \bk{w} & (\bk{M})_{12} 	& \b{0} \\
								\bk{w}^T 			& \bk{s}^T \bk{q} 		& \bk{r}^T & 0 \\
								(\bk{M})^T_{12} 	& \bk{r} & (\bk{M})_{22} & \b{0} \\
								\b{0} 			& 0 & \b{0} & -\bk{s}^T \bk{u} \\
						\end{array}
				\right]^{\text{-1}}%
				\left[
						\begin{array}{c}
							\bk{Q}^T \\
							\bk{q}^T \\
							\bk{U}^T \\
							\bk{u}^T
						\end{array}
					\right].
\end{align*}
Now we verify that the identities hold:
\begin{align*}
	\bko{Q} 		&= \left[ \mgap \bk{Q} \mgap \bk{q} \mgap \right] =  \left[ \mgap \bk{V} + \bzex{A}{P}\bk{S} \mgap \bk{v} + \bzex{A}{P}\bk{s} \mgap \right] = 
					\bko{V} + \bzex{A}{P}\bko{S},\\
	\bko{U} 		&= \left[ \mgap \bk{U} \mgap \bk{u} \mgap \right], \\
	(\bko{M})_{11} 	&= 	\left[
						\begin{array}{ c c }
							(\bk{M})_{11} 	& \bk{w} \\
							\bk{w}^T 		& \bk{s}^T \bk{q}
						\end{array}
					\right] = \bkoex{D}{V} + \bkoex{L}{V} + ( \bkoex{L}{V} )^T + \bko{S}^T \bzex{A}{P} \bko{S}, \\
	(\bko{M})_{12} 	&= 	\left[
						\begin{array}{ c c }
							(\bk{M})_{12} 	& \b{0} \\
							\bk{r}^T 		& 0
						\end{array}
					\right] = \bko{L}, \\
	(\bko{M})_{22} 	&= 	\left[
						\begin{array}{ c c }
							(\bk{M})_{22} 	& \b{0} \\
							\b{0} 		& -\bk{s}^T \bk{u}
						\end{array}
					\right] = -\bko{D}.
\end{align*}
Therefore we conclude that $\bkoex{A}{P}$ is of the form \eqref{eq:comp_csbfgsp} with $ k+1 $ replacing
the indices $ k $. \qed
\end{proof}

% Possibly add relation to SR1

\subsection{Limited Memory Compact Structured BFGS}
\label{subsec:limitedmemory}

The limited memory representations of Eqs. \eqref{eq:comp_csbfgsm} and \eqref{eq:comp_csbfgsp} are obtained
by storing only the last $ m \ge 1$ columns of $ \bk{S}, \bk{U} $ and $ \bk{V} $. \jjb{By setting $m \ll n$ limited memory strategies enable computational efficiencies and lower storage requirements, see e.g., \cite{Noc80}.} Updating $ \bk{S}, \bk{U} $ and $ \bk{V} $
requires replacing or inserting one column at each iteration. Let an underline below a matrix represent the 
matrix with its first column removed. That is, $ \underline{\b{S}}_k $ represents $ \bk{S} $ without its first column.
With this notation, a column update of a matrix, say $ \bk{S} $, by a vector $ \bk{s} $ is defined as follows.
\begin{equation*}
	\text{colUpdate}\left(\bk{S},\bk{s} \right) \equiv
	\begin{cases}
		[\: \bk{S} \: \bk{s}\:  ] 						& \text{ if } k < m \\
		[\: \underline{\b{S}}_k \: \bk{s}\:  ] 			& \text{ if } k \ge m. \\
	\end{cases}
\end{equation*}
Such a column update either directly appends a column to a matrix or first removes a column and then appends one. 
This column update will be used, for instance, to obtain $ \bko{S} $ from $ \bk{S} $ and $ \bk{s} $, i.e., $\bko{S}= \text{colUpdate}( \bk{S}, \bk{s} )$. Next, let an overline above a matrix represent
the matrix with its first row removed. That is, $ \overline{\b{S}^{T}_k \b{U}}_k $ represents $ \b{S}^{T}_k \bk{U} $ without its first row.
With this notation, a product update of, say $ \bk{S}^T\bk{U} $, by matrices $ \bk{S} $, \bk{U} 
and vectors $ \bk{s} $, $\bk{u}$ is defined as:
\begin{equation*}
	\text{prodUpdate} \left( \bk{S}^T\bk{U}, \bk{S}, \bk{U}, \bk{s}, \bk{u} \right) \equiv 
	\begin{cases}
		\left[
			\begin{array}{ c c }
				\bk{S}^T\bk{U} 		& \bk{S}^T\bk{u} \\
				\bk{s}^T\bk{U}	& \bk{s}^T \bk{u} 
			\end{array}
		\right] & \text{ if } k < m \vspace{0.1cm} \\		
		\left[
			\begin{array}{ c c }
				\left(\underline{\overline{\b{S}^{T}_k \b{U}_k}}\right) 			& 	\underline{\b{S}}_k^T\bk{u} \\
				\bk{s}^T \underline{\b{U}}_k					&	 \bk{s}^T \bk{u} 
			\end{array}
		\right] & \text{ if } k \ge m .\\
	\end{cases}
\end{equation*}
This product update is used to compute matrix products, such as, $ \bko{S}^T \bko{U} $, with 
$ \mathcal{O}(2mn) $ multiplications, instead of $ \mathcal{O}(m^2n) $ when the product 
$ \bk{S}^T \bk{U} $  had previously been stored. \jjbc{Note that a diagonal matrix can
be updated in this way by setting the rectangular matrices (e.g., $\bk{S}, \bk{U}$) to zero, such that e.g., 
$ \bkoex{D}{U} = \text{prodUpdate}(\bkex{D}{U},\b{0},\b{0},\bk{s},\bk{u})$. An upper triangular matrix
can be updated in a similar way, e.g., $ \bkoex{R}{U} = \text{prodUpdate}(\bkex{R}{U},\bk{S},\b{0},\bk{s},\bk{u}) $. 
To save computations, products with zeros matrices are never formed expliclty.}
%Moreover, we let ``$\text{diag}(\bk{S}^{T} \bk{U}) $" extract the diagonal elements of a matrix, say $\bk{S}^T\bk{U},$ while ``$  \text{tril}(\bk{S}^{T} \bk{U},-1) $" are the
%strictly lower triangular elements (elements below, excluding the main diagonal) and ``$ \text{triu}(\bk{S}^{T} \bk{U},0) $" the upper 
%triangular elements (elements above, including the main diagonal). 
\jjb{Section \ref{sec:alg} discusses computational and memory aspects in greater detail.}

\section{\jjbc{Limited-Memory Structured BFGS Line Search} Algorithms}
\label{sec:alg}
This section describes two line search algorithms with limited memory
 structured BFGS matrices. The compact representations
 enable efficient reinitialization strategies and search directions,
and we discuss these two components first, before presenting the overall algorithms. %are described in the next \cref{subsec:init,subsec:searchdir}.
 \subsection{Initializations}
 \label{subsec:init}
 For the limited memory BFGS matrix based on \cref{eq:intro_cbfgs} one commonly
 uses the initializations $ \bz{B}^{(k)} = \widehat{\sigma}_k \In  $, where
 $ \widehat{\sigma}_k = \bkmo{y}^T\bkmo{y} / \bkmo{s}^T \bkmo{y} $ (c.f. \cite{ByrNS94}).
 Choosing the initialization as a multiple of the identity matrix enables fast 
 computations with the matrix in \cref{eq:intro_cbfgs}. In particular, the
 inverse of this matrix may be computed efficiently by the Sherman-Morrison-Woodbury 
 identity. Because at the outset it is not necessarily obvious which
 initializations to use for the limited memory structured-BFGS (L-S-BFGS) matrices
 based on eqs. \eqref{eq:comp_csbfgsm} and \eqref{eq:comp_csbfgsp}, we investigate different approaches.
 We use the analysis in \cite{bb88}, which proposed formula $ \widehat{\sigma}_k $. 
 Additionally, in that work a second initialization $ \widehat{\sigma}^{(2)}_k = \bkmo{s}^T\bkmo{y} / \bkmo{s}^T \bkmo{s} $
 was proposed. Because in the S-BFGS methods the vectors $ \bhk{u} $ and $ \bk{u} $ are used
 instead of $ \bk{y} $ (unstructured BFGS), the initializations in this article are the below.
 \begin{equation}
 	\label{eq:initializations}
 	\sigma_{k+1} = 
	\begin{cases}
		\frac{\bk{u}^T\bk{u}}{\bk{s}^T\bk{u}} 		& \text{Init. 1} \\
		\frac{\bhk{u}^T\bhk{u}}{\bk{s}^T\bhk{u}} 	& \text{Init. 2} \\
		\frac{\bk{s}^T\bk{u}}{\bk{s}^T\bk{s}} 		& \text{Init. 3} \\
		\frac{\bk{s}^T\bhk{u}}{\bk{s}^T\bk{s}} 		& \text{Init. 4} \\
	\end{cases}
 \end{equation}
 Note that Init. 1 and Init. 2 are extensions of $ \widehat{\sigma}_k $ to structured methods. 
 Instead of using $ \bk{y} $ these initializations are defined by $ \bhk{u} $ and $ \bk{u} $. Init. 3
 and Init. 4 extend $ \widehat{\sigma}^{(2)}_k $. Observe that the vectors 
 $ \bhk{u} = \nabla \widehat{u}(\bko{x}) - \nabla \widehat{u}(\bk{x}) $ depend only on
 gradient information of $ \widehat{u}(\b{x}) $. In contrast, $ \bk{u} = \bko{K}\bk{s} + \bhk{u}  $
 depends on known second-derivative information, too. 
 Because the initial matrices 
$\b{A}^{\text{M}}_0$ and $ \b{A}^{\text{P}}_0$ affect the compact representations from 
Theorems \ref{thrm:comp_csbfgsm} and \ref{thrm:comp_csbfgsp} differently, we accordingly adjust
our initialization strategies for these two matrices. In particular, for L-S-BFGS-M the compact 
limited memory formula for $ \bk{B}^{\text{M}} $ simplifies if we take $ \b{B}^{\text{M}}_0 $ as
a multiple of the identity matrix: 
\begin{equation}
\label{eq:bminit}
	\b{B}^{\text{M}}_0 = \b{A}^{\text{M}}_0 + \b{K}_0 \equiv \sigma_k \b{I}.
\end{equation}
\jjb{The advantage of this choice is that it \jjbc{has similar} computational complexities to the L-BFGS formula from \cref{eq:intro_cbfgs}}. However by setting this default initialization for $ \b{B}^{\text{M}}_0  $ the corresponding limited memory matrices
$ \b{B}^{\text{M}}_k $ are not equivalent anymore to the full-memory matrices $ \b{B}^{\text{M}}_k $ defined by 
 \cref{eq:intro_sbfgsm}, even when $ k < m $. \jjb{In Section \ref{subsec:large_comp_m} computational techniques are discussed when $\b{B}^{\text{M}}_0$ is not taken as a multiple of the identity matrix.}
For L-S-BFGS-P we set $ \b{A}^{\text{P}}_0 = \sigma_k \b{I} $. This initialization, 
as long as $\sigma_k$ remains constant, implies that the limited memory compact representations from 
Theorem \ref{thrm:comp_csbfgsm} and the update formulas from \cref{eq:intro_sbfgsp} produce the same matrices when $ k < m $.

%(the limited memory matrix based on \cref{eq:comp_csbfgsm}),
% we use different initialization of the form
% \begin{equation}
% 	\label{eq:alg_csbfgsm_init}
%	\bsiidxex{0}{\Psi}{M}(\hat{\sigma},c) = \biidxex{0}{A}{M}(\widehat{\sigma}) + c\bz{K} = \widehat{\sigma} \In + c \bz{K}, 
% \end{equation} 
% for two scalars $ \widehat{\sigma} > 0 $ and $ c \in \{0,1\}  $. When the parameter $ c = 1 $, then $ \bz{K} $ is
% included as part of the initialization, otherwise, when $ c = 0 $, $ \bz{K} $ is not included. The initializations
% for the L-S-BFGS-P (the limited memory matrix based on \cref{eq:comp_csbfgsb}) are 
% \begin{equation}
% 	\label{eq:alg_csbfgsb_init}
%	\bsiidxex{0}{\Psi}{P}(\widehat{\sigma}) = \biidxex{0}{A}{P}(\widehat{\sigma}) = \widehat{\sigma} \In, 
% \end{equation}
% with $ \hat{\sigma} > 0 $. Similar to the value $ \bk{\sigma} $, we define the scalar
% \begin{equation}
% 	\label{eq:alg_hatsigma}
%	\widehat{\sigma}_k = \frac{ \bkmo{u}^T \bkmo{u} }{\bkmo{s}\bkmo{u}}.
% \end{equation}
 
% When, moreover, the initialization is $ \bsiidxex{0}{\Psi}{M}(\widehat{\sigma}_k,0) = \widehat{\sigma}_k\In $ is used,
 
 \subsection{Search Directions}
  \label{subsec:searchdir}
 The search directions for line search algorithms, with the structured BFGS approximations,
 are computed as
 \begin{equation}
 	\label{eq:alg_p}
		\bk{p} = - ( \bk{K} + \bk{A} )^{-1} \bk{g},
 \end{equation} 
 where $ \bk{g} = \nabla f(\bk{x}) $ and where $ \bk{A} $ is either the limited memory version of
$ \bkex{A}{M} $ from \cref{eq:comp_csbfgsm} or $ \bkex{A}{P} $ from \cref{eq:comp_csbfgsp}. 
When $ \bkex{A}{M} $ is used, we apply
the expression of the inverse from Corollary \ref{cor:comp_icsbfgsm}, in order to compute search directions.
In particular, with the initialization strategy $ \b{B}^{\text{M}}_0 = \sigma_k \b{I} $ from the preceding section, the
search directions  $ \cref{eq:alg_p} $ are computed efficiently by % 
  \begin{equation}
  \label{eq:pm}
  	\bkex{p}{M} = - \frac{\bk{g}}{\sigma_k}  -
	\left[ \bk{S} \quad  \bk{U} \right]
			\begin{bmatrix}
				(\bkex{T}{U})^T\left( \bkex{D}{U} + 1 / \sigma_k \bk{U}^T \bk{U} \right) \bkex{T}{U}  	& -\frac{(\bkex{T}{U})^T}{\sigma_k} \\
				-\frac{\bkex{T}{U}}{\sigma_k}											& \b{0}_{m \times m}
			\end{bmatrix}
			\left(
				\begin{bmatrix}
					\bk{S}^T\bk{g}  \\
					\bk{U}^T \bk{g} 
				\end{bmatrix}
			\right), % \bk{g} 
  \end{equation}
  where $ \bkex{T}{U} $ is defined in Corollary \ref{cor:comp_icsbfgsm}. \jjb{This computation is done efficiently assuming that all matrices have been updated before, such as $ \bk{U}^T \bk{U} $. Omitting terms of order $m$, the multiplication complexity for this search direction is $ \mathcal{O}(n(4m+1) + 3m^2)$. In particular, computing $ \bkex{p}{M} $
  can be done by: two vector multiplies with the $ n \times 2m $ matrix $ [ \: \bk{S} \: \bk{U} \: ] $ (order $ 4nm $), the scaling $ \frac{\bk{g}}{\sigma_k} $ (order $n$) and a matrix vector product with a structured $ 2m \times 2m  $ matrix. Since $\bkex{T}{U}$ represents a solve with an $ m \times m $ upper triangular matrix the vector product
  with the middle $ 2m \times 2m $ matrix is done in order $ 3m^2 $. 
  When $ \bkex{A}{P} $ is used, search directions are computed by solves of the linear system
  $ (\bk{K} + \bkex{A}{P}) \bkex{p}{P} = - \bk{g}  $. \jjbc{Because of the compact
  representation of $\bkex{A}{P}$ we can exploit structure in solving this system,
  as is described in Section \ref{subsec:large_comp_p} }}
  
  \subsection{Algorithms}
  \label{sec:algs}
  \jjbc{Similar to Petra et al. \cite{PetraChiangAnitescu19}, we use a strong Wolfe line search
  in our implementations of the new limited-memory compact representations.}
  %As in the work of Petra et al. \cite{PetraChiangAnitescu19}, the compact representations of the 
  %structured BFGS formulas are implemented in a strong Wolfe line search algorithm based on \cite{MoreThuente94}.
  For nonnegative constants $ 0 < c_1 \le c_2 $, the current iterate $ \bk{x} $ and search direction $ \bk{p} $,
   the strong Wolfe conditions define the step length parameter
  $ \alpha $ by two inequalities
  \begin{align}
  	\label{eq:wolfecond}
	&f(\bk{x} + \alpha \bk{p}) \le f(\bk{x}) + c_1 \alpha (\bk{p}^T \nabla f(\bk{x})), \text{ and } \\
	&\left| \bk{p}^T \nabla f(\bk{x} + \alpha \bk{p}) \right| \le c_2 \left| \bk{p}^T \nabla f(\bk{x}) \right|. \nonumber
  \end{align}  
   Because the S-BFGS-M
  matrix from \cref{eq:comp_csbfgsm} is positive definite as long as $ \bk{s}^T \bk{u} > 0 $ for $ k \ge 0 $ 
  \jjbc{(rather than $ \bk{s}^T\bk{y} >0 $ for $ k \ge 0$ for L-BFGS)},
  the line searches in our algorithms include this condition. \jjbc{As in \cite[Appendix A]{PetraChiangAnitescu19}
  a variant of the Mor\'{e}-Thuente \cite{MoreThuente94} line search is used. This line search is identical
  to the one of Mor\'{e}-Thuente, except for one condition. Specifically, given a trial step length $ \alpha_t $ and
  trial $ \b{u}_t $,
  our line search terminates when the conditions in \eqref{eq:wolfecond} and additionally $  \bk{s}^T\b{u}_t > 0  $ holds.
  \cite[Proposition 17]{PetraChiangAnitescu19} ensures the existence of a step length $ \alpha $ that satisfies all of the above conditions.
  Such a line search variant is straight forward to implement, by adding one additional condition
  to a Mor\'{e}-Thuente line search.} 
  Moreover, when S-BFGS-M is used,
  new search directions are computed by using the inverse from Corollary \ref{cor:comp_icsbfgsm}. In contrast,
  because the S-BFGS-P matrix from \cref{eq:comp_csbfgsp} is not necessarily positive definite even
  if $ \bk{s}^T \bk{u} > 0 $ for $ k \ge 0 $ (see \cite{PetraChiangAnitescu19}), our implementation checks
  whether $ \bk{K} + \bkex{A}{P} $ is positive definite, before computing a new search direction
  \jjbc{However, if it is known that $ \bk{K} $ is positive definite for all $ k \ge 0 $ (which is often the case
  in applications) than ensuring that $ \bk{s}^T\bk{u} > 0 $ for $ k \ge 0 $ ensure positive definiteness,
  in this case too}.
  If this matrix is positive definite, then a new search direction is computed by solving the linear system
  $ (\bk{K} + \bkex{A}{P}) \bkex{p}{P} = - \bk{g} $. Otherwise the search direction is computed by
  solving the system $ (\bk{K} + \bkex{A}{P} + \delta \In) \bkex{p}{P} = - \bk{g} $, where the scalar $ \delta > 0 $
  ensures that $  (\bk{K} + \bkex{A}{P} + \delta \In) \succ 0 $ (Here $ \delta $ is chosen as the the
  first $ \delta =10^j, j = 0,1,\ldots$ that yields a positive definite matrix). The proposed limited memory line search
  algorithms are listed in Algorithm \ref{alg:alg_csbfgsm} and Algorithm \ref{alg:alg_csbfgsp}. % \cref{alg:alg_csbfgsm,alg:alg_csbfgsp}. 
%  As in \cite{PetraChiangAnitescu19},
%  both algorithms use a strong Wolfe line search based on \cite{MoreThuente94}
 
%Our analysis leads to the algorithm in \cref{alg:buildtree}.
% TODO: This algorithm uses the notation [M]_{i \le j} to define the elements of a upper triangular matrix.
% This may be somewhat inaccurate.
\begin{algorithm}[tbhb]
\caption{Limited Memory Structured-BFGS-Minus (L-S-BFGS-M)}
\label{alg:alg_csbfgsm}
\begin{algorithmic}[1]
\STATE{Initialize: $ k = 0$, $ m > 0 $, $ \epsilon > 0 $, $ \sigma_k > 0 $, $ 0 < c_1 \le c_2 $, $\bk{x}$, $ \bk{g} = \nabla f(\bk{x}) = \nabla \hat{k}(\bk{x}) + \nabla \hat{u}(\bk{x}) $, 
$ \bk{S} = 0, \bk{U} = 0, \bkex{D}{U} = 0, (\bkex{R}{U})^{-1} = 0 $, $\bk{U}^T\bk{U}=0$, $ \b{H}_0 = (1/ \sigma_k) \b{I}$, $ \bsk{\Theta} = \left[ \bk{S} \quad  \b{H}_0\bk{U}  \right] $} % $\bk{S}^T\bk{U}=0$, 
\WHILE{$ \| \bk{g} \|_{\infty} > \epsilon $}
\STATE{Compute: \label{alg:csbfgsm_p}
$$ 
	\bk{p} = - \b{H}_0\bk{g} + \bsk{\Theta} \bk{M} (\bskt{\Theta} \bk{g}),
$$
where 
$$ 
	\bk{M} =  
	\biggl[
		\begin{smallmatrix}
			(\bk{R}^{\scriptsize \textnormal{U}})^{-T}\left( \bk{D}^{\scriptsize \textnormal{U}} + \bk{U}^T \b{H}_0 \bk{U} \right) (\bk{R}^{\scriptsize \textnormal{U}})^{-1}  	& -(\bk{R}^{\scriptsize \textnormal{U}})^{-T} \\
				-(\bk{R}^{\scriptsize \textnormal{U}})^{-1}														& \b{0}
		\end{smallmatrix}
	\biggr].
$$}
\STATE{Strong Wolfe line search: 
$$
	\bko{x} = \bk{x} + \alpha \bk{p}, 
$$ 
where $ \alpha > 0 $, $ \bko{x} $ satisfies strong Wolfe conditions (cf. \cite{PetraChiangAnitescu19} and \cref{eq:wolfecond}), 
$\bk{s} = \bko{x} - \bk{x} $, $ \bk{s}^T\bk{u} > 0 $.}
\STATE{Updates: $\bko{g} = \nabla f(\bko{x})$, $\bk{u} = \nabla^2 \hat{k}(\bko{x}) \bk{s} +(  \nabla \hat{u}(\bko{x}) - \nabla \hat{u}(\bk{x}))$}
\STATE{ $\bko{S} = \text{colUpdate}(\bk{S}, \bk{s}) $  }
\STATE{ $\bko{U} = \text{colUpdate}(\bk{U}, \bk{u}) $ } 
\jjbc{\STATE{ $ \bkoex{R}{U} = \text{prodUpdate}( \bkoex{R}{U},\bk{S},\b{0},\bk{s}, \bk{u}) $ }}
%\STATE{ $\bko{S}^T\bko{U} = \text{prodUpdate}( \bk{S}^T\bk{U},\bk{S},\bk{U},\bk{s}, \bk{u}) $ }
\STATE{ $\bko{U}^T\bko{U} = \text{prodUpdate}( \bk{U}^T\bk{U},\bk{U},\bk{U},\bk{u}, \bk{u}) $ }
%$$
%	\bko{S} =
%	\begin{cases}
%		\left[ \bk{S} \quad \bk{s} \right], 				& k \le m \\
%		\left[ \underbar{\b{S}}_k \quad \bk{s} \right], 	& k > m,
%	\end{cases} \quad 
%	\bko{U} =
%	\begin{cases}
%		\left[ \bk{U} \quad \bk{u} \right], 				& k \le m \\
%		\left[ \underbar{\b{U}}_k \quad \bk{u} \right] 	& k > m,
%	\end{cases}
%$$
%$$
%	\bko{S}^T\bko{U} =
%		\begin{cases}
%			\bko{S}^T\bko{U}, 													& k \le m \\
%			\begin{bmatrix}
%				\underbar{\b{S}}^T_k\underbar{\b{U}}_k 	& \underbar{\b{S}}^T_k\bk{u} \\
%				\bk{s}^T\underbar{\b{U}}_k 			& \bk{s}^T \bk{u}
%			\end{bmatrix}, 														& k > m.		
%		\end{cases}
%$$
\jjbc{\STATE{$ \bkoex{D}{U} = \text{prodUpdate}( \bkex{D}{U},\b{0},\b{0},\bk{s}, \bk{u}) $}}
%\STATE{$ \bkoex{D}{U} = \text{diag}(\bko{S}^T\bko{U}) $}
%\STATE{$ \bkoex{R}{U} = \text{triu}(\bko{S}^T\bko{U},0) $}
%\STATE{$ [\bkoex{R}{U}]_{ij} = [ \bko{S}^T\bko{U} ]_{i \le j} $, $ [\bkoex{D}{U}]_{ii} = [ \bko{S}^T\bko{U} ]_{ii} $}
%$$
%	\bkoex{D}{U} =
%	\begin{cases}
%		\begin{bmatrix}
%			\bkex{D}{U} 	& \\
%					& \bk{s}^T \bk{u}
%		\end{bmatrix}, 							& k \le m \\
%		
%		\begin{bmatrix}
%			\overline{\underbar{\b{D}}}_k^{U} 	& \\
%					& \bk{s}^T \bk{u}
%		\end{bmatrix}, 							& k > m
%	\end{cases}, \quad 
%	\bkoex{R}{U} =
%	\begin{cases}
%		\begin{bmatrix}
%			\bkex{R}{U} 		& \bk{S}^T\bk{u} \\
%							& \bk{s}^T \bk{u}
%		\end{bmatrix}, 							& k \le m \\
%		
%		\begin{bmatrix}
%			\overline{\underbar{\b{R}}}_k^{U} 		& \underbar{\b{S}}^T_k \bk{u} \\
%											& \bk{s}^T \bk{u}
%		\end{bmatrix}, 							& k > m
%	\end{cases},
%$$
\STATE{Compute: $ \sigma_{k+1}$}
\STATE{ $ \b{H}_0 = (1/ \sigma_{k+1})\b{I}$, update $ \bko{M} $, $ \bsko{\Theta} $ using Theorem \ref{thrm:comp_csbfgsm}, $ k = k+1 $}
%$$
%	\bko{M} =
%	\begin{bmatrix}
%	(\bkoex{R}{U})^{-T}\left( \bkoex{D}{U} + \bko{U}^T \bs{\Psi}^{-1}_0 \bko{U} \right) (\bkoex{R}{U})^{-1}  	& -(\bkoex{R}{U})^{-T} \\
%				-(\bkoex{R}{U})^{-1}
%	\end{bmatrix},
%$$
%$  \bsko{\Theta} = \left[ \bko{S} \quad  \bs{\Psi}^{-1}_0\bko{U}  \right] $, $ k = k+1 $. 
%}
%\bkex{D}{U} = \diag(\bko{S}^T \bko{U}) $, \bkex{T}{U} = (\text{triu}(\bko{S}^T \bko{U}))^{-1} }
\ENDWHILE
\RETURN $\bk{x}$
\end{algorithmic}
\end{algorithm}

Note that $ \bskt{\Theta} \bk{g} $ on Line \ref{alg:csbfgsm_p} in Algorithm \ref{alg:alg_csbfgsm} is computed as 
$ 
\begin{bmatrix}
	\bk{S}^T\bk{g} \\
	\b{H}_0(\bk{U}^T\bk{g})
\end{bmatrix} 
$
so that only one linear solve with $ \b{H}_0 = (\b{K}_0 + \b{A}^M_0)^{-1}  $ is needed, when the algorithm does
not use a multiple of the identity as the initialization.

\begin{algorithm}[tbhb]
\caption{Limited Memory Structured-BFGS-Plus (L-S-BFGS-P)}
\label{alg:alg_csbfgsp}
\begin{algorithmic}[1]
\STATE{Initialize: $ k = 0$, $ m > 0 $, $ \epsilon > 0 $, $ \sigma_k > 0 $, $ 0 < c_1 \le c_2 $, $\bk{x}$,  $ \bk{g} = \nabla f(\bk{x}) = \nabla \hat{k}(\bk{x}) + \nabla \hat{u}(\bk{x}) $, 
$ \bk{K} = \nabla^2 \hat{k}(\bk{x}) $, $ \bk{S} = 0, \bk{U} = 0, \bk{V} = 0, \bkex{D}{U} = 0, \bkex{L}{U} = 0, \bkex{D}{V} = 0, \bkex{L}{V} = 0 $, $\bsk{\Omega} = 0$, $\bk{S}^T\bk{S} = 0$, $\bk{A} = \sigma_k \b{I}$} % $\bk{S}^T\bk{U} = 0$, $\bk{S}^T\bk{V} = 0$, 
% $ \bsk{\Omega} = 0 $, 
%$ \bk{N} =  
%	\biggl[
%		\begin{smallmatrix}
%			\bkex{L}{V} + \bkex{D}{V} + (\bkex{L}{V})^T + \bk{S}^T \bk{A} \bk{S}   					& \bkex{L}{U} \\
%				(\bkex{L}{U})^T															& \bkex{D}{U}
%		\end{smallmatrix}
%	\biggr]^{-1} = 0.
%$
\WHILE{$ \| \bk{g} \|_{\infty} > \epsilon $ }
\IF{ \label{alg:csbfgsp_pdck} $ ( \bk{K} + \bk{A} ) \not\succ 0 $}
\STATE{Find  $ \delta > 0 $ such that $ ( \bk{K} + \bk{A} + \delta \In ) \succ 0 $}
\ENDIF
\STATE{Solve: \label{alg:csbfgsp_p}
$$ 
	(\bk{K} + \bk{A})\bk{p} = -\bk{g}
$$}
\STATE{Strong Wolfe line search: 
$$
	\bko{x} = \bk{x} + \alpha \bk{p}, 
$$ 
where $ \alpha > 0 $, $ \bko{x} $ satisfies strong Wolfe conditions (cf. \cite{PetraChiangAnitescu19} and \cref{eq:wolfecond}), 
$ \bk{s} = \bko{x} - \bk{x} $.}
\STATE{Updates: $\bko{g} = \nabla f(\bko{x})$, $ \bko{K} = \nabla^2 \hat{k}(\bko{x}) $, $ \bk{v} = \bko{K}\bk{s} $, $\bk{u} = \bk{v} +(  \nabla \hat{u}(\bko{x}) - \nabla \hat{u}(\bk{x}))$}
\STATE{ $\bko{S} = \text{colUpdate}(\bk{S}, \bk{s}) $  }
\STATE{ $\bko{U} = \text{colUpdate}(\bk{U}, \bk{u}) $ }
\STATE{ $\bko{V} = \text{colUpdate}(\bk{V}, \bk{v}) $ }
\jjbc{\STATE{ $\bkoex{L}{U} = \text{prodUpdate}( \bkex{L}{U},\b{0},\bk{U},\bk{s}, \b{0}) $ }}
\jjbc{\STATE{ $\bkoex{L}{V} = \text{prodUpdate}( \bkex{L}{V},\b{0},\bk{V},\bk{s}, \b{0}) $ }}
%\STATE{ $\bko{S}^T\bko{U} = \text{prodUpdate}( \bk{S}^T\bk{U},\bk{S},\bk{U},\bk{s}, \bk{u}) $ }
%\STATE{ $\bko{S}^T\bko{V} = \text{prodUpdate}( \bk{S}^T\bk{V},\bk{S},\bk{V},\bk{s}, \bk{v}) $ }
\STATE{ $\bko{S}^T\bko{S} = \text{prodUpdate}( \bk{S}^T\bk{S},\bk{S},\bk{S},\bk{s}, \bk{s}) $ }
%$$
%	\bko{S} =
%	\begin{cases}
%		\left[ \bk{S} \quad \bk{s} \right], 				& k \le m  \\
%		\left[ \underbar{\b{S}}_k \mgap \bk{s} \right], 	& k > m,
%	\end{cases} \quad
%	\bko{U} =
%	\begin{cases}
%		\left[ \bk{U} \quad \bk{u} \right], 			 	& k \le m \\
%		\left[ \underbar{\b{U}}_k \quad \bk{u} \right], 	& k > m,
%	\end{cases} 
%$$
%$$
%		\bko{V} =
%	\begin{cases}
%		\left[ \bk{V} \quad \bk{v} \right], 				& k \le m  \\
%		\left[ \underbar{\b{V}}_k \quad \bk{v} \right], 	& k > m,
%	\end{cases}
%$$
%$$
%	\bko{S}^T\bko{U} =
%		\begin{cases}
%			\bko{S}^T\bko{U}, 													& k \le m \\
%			\begin{bmatrix}
%				\underbar{\b{S}}^T_k\underbar{\b{U}}_k 	& \underbar{\b{S}}^T_k\bk{u} \\
%				\bk{s}^T\underbar{\b{U}}_k 			& \bk{s}^T \bk{u}
%			\end{bmatrix}, 														& k > m,		
%		\end{cases}
%$$
%$$
%	\bko{S}^T\bko{V} =
%		\begin{cases}
%			\bko{S}^T\bko{V}, 													& k \le m \\
%			\begin{bmatrix}
%				\underbar{\b{S}}^T_k\underbar{\b{V}}_k 	& \underbar{\b{S}}^T_k\bk{v} \\
%				\bk{s}^T\underbar{\b{V}}_k 			& \bk{s}^T \bk{v}
%			\end{bmatrix}, 														& k > m,		
%		\end{cases}
%$$
%\STATE{$\bkoex{L}{U} = \text{tril}(\bko{S}^T\bko{U},-1)$}
%\STATE{$\bkoex{L}{V} = \text{tril}(\bko{S}^T\bko{V},-1)$}
%\STATE{$\bkoex{D}{U} = \text{diag}(\bko{S}^T\bko{U})$}
%\STATE{$\bkoex{D}{V} = \text{diag}(\bko{S}^T\bko{U})$}
\jjbc{\STATE{ $\bkoex{D}{U} = \text{prodUpdate}( \bkex{D}{U},\b{0},\b{0},\bk{s}, \bk{u}) $ }}
\jjbc{\STATE{ $\bkoex{D}{V} = \text{prodUpdate}( \bkex{D}{V},\b{0},\b{0},\bk{s}, \bk{v}) $ }}
%\STATE{$ [\bkoex{L}{U}]_{ij} = [ \bko{S}^T\bko{U} ]_{j<i} $, $ [\bkoex{D}{U}]_{ii} = [ \bko{S}^T\bko{U} ]_{ii} $,
%$ [\bkoex{L}{V}]_{ij} = [ \bko{S}^T\bko{V} ]_{j<i} $, $ [\bkoex{D}{V}]_{ii} = [ \bko{S}^T\bko{V} ]_{ii} $
%}
\STATE{Compute: $ \sigma_{k+1}$}
\STATE{ $ \b{A}_0 = (1/ \sigma_{k+1})\b{I}$, update $ \bsko{\Omega} = \left[ \: \bko{V} + \b{A}_0\bko{S} \: \bk{U} \: \right] $ }
\STATE{
$$
\bko{A} = \b{A}_0 -  \bsko{\Omega}
		{\begin{bmatrix}
			\bkoex{D}{V} + \bkoex{L}{V} + (\bkoex{L}{V})^T + \bko{S}^T \b{A}_0 \bko{S} 	& \bkoex{L}{U} \\
			(\bkoex{L}{U})^T 												& -\bkoex{D}{U}
		\end{bmatrix}}^{-1}
		\bskot{\Omega}	
$$
\STATE{$ k = k+1 $}
}
\ENDWHILE
\RETURN $\bk{x}$
\end{algorithmic}
\end{algorithm}
%Comparing \cref{alg:alg_csbfgsm,alg:alg_csbfgsp},
Algorithm \ref{alg:alg_csbfgsp} is expected to be
computationally more expensive than Algorithm \ref{alg:alg_csbfgsm} because it tests for the positive definiteness
of $ \bk{K} + \bk{A} $ in Line \ref{alg:csbfgsp_pdck} and it computes search directions by the solve
in Line \ref{alg:csbfgsp_p}. However, the structured quasi-Newton approximation in Algorithm \ref{alg:alg_csbfgsp}
may be a more accurate approximation of the true Hessian (see \cite{PetraChiangAnitescu19}), which
may result in fewer iterations or better convergence properties. \jjbc{Note that as in 
\cite[Section 3.1.2]{PetraChiangAnitescu19} computational efforts for ensuring
positive definiteness may largely reduced by e.g., defining $ \delta = \text{max}(0,(\varepsilon-
(\bk{u} + \bk{v})^T\bk{s})/\|\bk{s}\|^2 $), for $0 < \varepsilon$.} Unlike Algorithm \ref{alg:alg_csbfgsp}, 
Algorithm \ref{alg:alg_csbfgsm} does not require solves involving large linear systems. 

 \subsection{Large-Scale Computation Considerations}
 \label{subsec:large_comp}
 \jjb{This section discusses computational complexity and memory requirements of the
structured Hessian approximations when the problems \jjbc{are} large. In particular, if $n$ is large 
 the Hessian matrices $\bk{K}$ typically exhibit additional structure, such as being diagonal or
 sparse. When $ \bk{K} $ is sparse and solves with it can be done efficiently,
 the compact representation of $ \bkex{A}{M} $ and $ \bkex{A}{P} $ can be exploited
 to compute inverses of $ \bk{K} + \bk{A} $ efficiently. \jjbc{Note that Algorithm \ref{alg:alg_csbfgsm} % \eqref{eq:bminit}
 is directly applicable to large problems, because the formula in \eqref{eq:pm} does not use solves with \bk{K} }. \jjbc{Nevertheless, observe that} the matrices 
 $ \bk{K} + \bk{A} $, (with limited memory $\bk{A}$  from Theorem \ref{thrm:comp_csbfgsm} or Theorem \ref{thrm:comp_csbfgsp}, respectively), have the form with $m \ll n$:
 \begin{equation}	
 	\label{eq:compform}
	 \bk{K} + \bk{A} \equiv \bh{K}_0 -
	 	\left[
	 		\begin{array}{c}
	 			\vline \\
	 			\vline \\
	 			\bsk{\Xi} \\
	 			\vline \\
	 			\vline \\
	 		\end{array}
	 	\right]
	 		\big[
	 			\bk{M} 
	 		\big]^{-1}
	 		\big[ \rule[.5ex]{3.5em}{0.4pt} \mgap \bs{\Xi}_k^T \mgap \rule[.5ex]{3.5em}{0.4pt} \big],
 \end{equation}
for some $ \bh{K}_0 $. If $\bkex{A}{M}$ is used in \cref{eq:compform} then $ \bh{K}_0 = \bz{K} + \bzex{A}{M} $
 and $ \bsk{\Xi} $, $ \bk{M} $ correspond to the remaining terms in Theorem \ref{thrm:comp_csbfgsm}. Using $\bkex{A}{P}$ in \cref{eq:compform} then $ \bh{K}_0 = \bk{K} + \bzex{A}{P} $ and $ \bsk{\Xi} $, $ \bk{M} $ correspond to the remaining terms in Theorem \ref{thrm:comp_csbfgsp}. Because of its structure the matrix in \cref{eq:compform}
 can be inverted efficiently by the Sherman-Morrison-Woodbury formula as long as solves with $ \bh{K}_0 $ can be done efficiently.
 Next, L-S-BFGS-M and L-S-BFGS-P are discussed in the situation when solves with 
 $ \bh{K}_0 $ are done efficiently. Afterwards we relate these methods to S-BFGS-M, S-BFGS-P and BFGS, L-BFGS.}
 
 \subsubsection{Computations for L-S-BFGS-M}
 \label{subsec:large_comp_m}
 \jjb{The most efficient computations are achieved when $ \bh{K}_0 $ is set as a multiple of the identity matrix
  $ \sigma_k \b{I} $ (cf. \Cref{subsec:searchdir} with $\mathcal{O}(n(4m+1)+3m^2)$ multiplications).
  This approach however omits the $ \bz{K} $ term. Nevertheless, when $ \bz{K} $ has additional structure such that factorizations and solves with it can be done
  in, say $ nl $ multiplications, search directions can be computed efficiently in this case, without omitting $ \bz{K} $. In particular, the search direction is computed as 
  $ \bkex{p}{M} = -(\bk{K} + \bkex{A}{M})^{-1}\bk{g}= -\bkex{H}{M}\bk{g} $ where $ \bkex{H}{M} $ is the inverse from 
  \cref{cor:comp_icsbfgsm}. The initialization matrix is $ \bzex{H}{M} = ( \sigma_k \b{I} + \bz{K})^{-1} $. To determine the search direction two matrix vector products with
  the $ n \times 2m $ matrices $ [\: \bk{S} \quad \bzex{H}{M} \bk{U} \:] $ are required, at
  complexity $ \mathcal{O}(4nm + 2nl) $. The product with the $ 2m \times 2m $ middle
  matrix is done at $\mathcal{O}(2nm + nl + 2m^2)$. Subsequently, $ -\bzex{H}{M}\bk{g} $
  is obtained at $nl$ multiplications. The total complexity is thus 
  $ \mathcal{O}(n(6m+4l)+2m^2) $. Note that if $ \sigma_k $ is set to a constant value, say $ \sigma_k = \bar{\sigma} $, then the complexity can be further reduced by storing 
  the matrix $ \bar{\b{U}}_k = [\: \bar{\b{u}}_{k-m} \ldots \bar{\b{u}}_{k-1} \: ] $, where $ \bar{\b{u}}_{i} = (\bz{K}+\bar{\sigma} \b{I})^{-1}\b{u}_i $. The computational cost in this situation is $\mathcal{O}(n(4m+l)+3m^2)$, excluding the
  updating cost of the vector $\bar{\b{u}}_{i}$ at order $nl$. With a constant $ \sigma_k $ only one factorization of $ (\bz{K}+\bar{\sigma} \b{I}) $ is required.}
  
  \subsubsection{Computations for L-S-BFGS-P}
  \label{subsec:large_comp_p}
  \jjb{When $ \bkex{A}{P} $ is used in \cref{eq:compform} with $ \bhz{K} = (\bk{K} + \bzex{A}{P}) $ and $ \bhk{Q} = \bhz{K}^{-1} \bk{Q} $, $ \bhk{U} = \bhz{K}^{-1} \bk{U} $ the inverse has the form
  \begin{equation*}
	  (\bk{K} + \bkex{A}{P})^{-1} = \bhz{K}^{-1}\left( \b{I}_n + \bsk{\Xi}\left( \bk{M} - \bs{\Xi}_k^T\bhz{K}^{-1} \bs{\Xi}_k \right)^{-1} \bs{\Xi}^T_k \bhz{K}^{-1}\right),
  \end{equation*}
  where
    $
	   \bs{\Xi}_k^T\bhz{K}^{-1} \bs{\Xi}_k = 
	   		\bigg[
			\begin{smallmatrix}
	   			\bk{Q}^T \bhk{Q} 	& \bk{Q}^T\bhk{U}  \\
	   			\bk{U}^T\bhk{Q} 												
	   			& \bk{U}^T\bhk{U}
	   		\end{smallmatrix}
			\bigg]
  $
%  \begin{equation*}
%	   \bs{\Xi}_k^T\bhz{K}^{-1} \bs{\Xi}_k = 
%	   		\begin{bmatrix}
%	   			\bk{Q}^T \bhk{Q} 	& \bk{Q}^T\bhk{U}  \\
%	   			\bk{U}^T\bhk{Q} 												
%	   			& \bk{U}^T\bhk{U}
%	   		\end{bmatrix}
%  \end{equation*}
  and $ \bsk{\Xi} $, $ \bk{M} $ are defined in \cref{thrm:comp_csbfgsp}. Assuming
  that $ \bk{M} $, $\bk{Q}$, $\bk{U}$ had previously been updated, computing the search direction $ \bkex{p}{P} = -(\bk{K} + \bkex{A}{P})^{-1}\bk{g} $ may be done as follows; First, 
  $ \bhk{Q} $, $ \bhk{U} $ are computed in $ \mathcal{O}(2nlm) $ multiplications. Then
  the $ 2m \times 2m $ matrix $ \bs{\Xi}_k^T\bhz{K}^{-1} \bs{\Xi}_k $ is formed in
  $ \mathcal{O}(3nm^2) $ multiplications. Combining the former terms and solving with the (small) $2m \times 2m $ matrix explicitly, the direction $ \bkex{p}{P} $ is computed
  in $ \mathcal{O}(n(2lm + 3m^2 + 4m + 1) + m^3) $ multiplications. Note that this approach requires an additional $ 2nm $ storage locations for the matrices $ \bhk{Q} $, $ \bhk{U} $. Two additional remarks; first, since $ \bk{Q} = \bk{V} + \bzex{A}{P} \bk{S} $, the update of $ \bk{Q} $ uses $ \mathcal{O}(nl) $ multiplications to form a new $ \b{v}_k $ and additional $ nm $ multiplications if  $ \bzex{A}{P} = \sigma_k \b{I} $. If 
  $ \sigma_k $ remains constant, say $ \sigma_k = \bar{\sigma} $, then the update 
  of $ \bk{Q} $ is done at only $ \mathcal{O}(nl) $ multiplications, because $\bzex{A}{P} \bk{S}$ does not need to be recomputed each iteration. Second, if $ \bk{K} = \bz{K} $, in other words if $ \bk{K} $ is a constant matrix then Theorems \ref{thrm:comp_csbfgsm} and \ref{thrm:comp_csbfgsp} reduce to the same expressions yielding the same computational complexities.}
 
   \subsubsection{Memory Usage and Comparison}
   \label{subsec:memory}
   \jjb{This section addresses the memory usage of the proposed representations
   and relates their computational complexities to existing methods. As an overall guideline, the representations from \cref{eq:compform} use $ 2nm + 4m^2 $ storage locations, excluding the $ \bhz{K} $ term. This estimate is refined if the particular structure of the matrix $ \bk{M} $ is taken into consideration. For example, the
   matrices $ \bkex{T}{U} $ and $ \bkex{D}{U} $ from Theorem \ref{thrm:comp_csbfgsm} are 
   upper triangular and diagonal, respectively. Thus, when $ \bzex{H}{M} = \sigma_k \b{I} $, and when the matrix $ \bk{U}^T\bk{U} \in \mathbb{R}^{m \times m} $ is stored and updated, the memory requirement for the limited memory version of $\bkex{H}{M}$
   in \cref{thrm:comp_csbfgsm} are $ \mathcal{O}(2nm + \frac{3}{2}m^2 + m) $ locations. We summarize the computational demands of the different methods in a table. % \ref{tbl:complx_comp}   
   \begin{table}
   \label{tbl:complx_comp}
   \caption{Comparison of computational demands for BFGS,L-BFGS,S-BFGS-M,S-BFGS-P,L-S-BFGS-M,L-S-BFGS-P, excluding storage of $ \bk{K} $ and where \jjbc{solves with $ \bk{K} $ are assumed to cost $ \mathcal{O}(nl) $ multiplications and vector multiplies cost $ \mathcal{O}(l) $.} Terms of order $\mathcal{O}(m)$ or lower are omitted. ($^{\dagger}$ The search direction cost for L-S-BFGS-P do not include the identity regularization with $ \delta $ from Sec. \ref{sec:algs}, as other techniques are possible.) }
   \centering
   	\begin{tabular}{l l l l}
   		\text{Method} 	& \text{Search Direction} 	& \text{Memory} 		& \text{Update} \\
   		\hline BFGS 	& $\mathcal{O}(n^2)$ 		& $\mathcal{O}(n^2)$ &  $\mathcal{O}(n^2)$ \\
   		L-BFGS (\cref{eq:intro_cbfgs}, \cite{ByrNS94}) & $\mathcal{O}(n(4m+1)+m^2)$ & $ \mathcal{O}(2nm + \frac{3}{2}m^2) $ & $ \mathcal{O}(2nm) $ \\
   		S-BFGS-M (\cref{eq:intro_sbfgsm},\cite{PetraChiangAnitescu19}) & $\mathcal{O}(n^2)$ & $\mathcal{O}(n^2)$ & $\mathcal{O}(n^2)$ \\
   		S-BFGS-P (\cref{eq:intro_sbfgsp},\cite{PetraChiangAnitescu19}) & $\mathcal{O}(n^2)^{\dagger}$ & $\mathcal{O}(n^2)$ & $\mathcal{O}(n^2)$ \\
   		L-S-BFGS-M (\cref{eq:pm}) & $\mathcal{O}(n(4m+1)+m^2)$ & $\mathcal{O}(2nm + \frac{3}{2}m^2)$ & $\mathcal{O}(2nm + l)$ \\
   		L-S-BFGS-M (\cref{eq:compform}) & $ \mathcal{O}(n(6m+4l)+m^2)$ & $\mathcal{O}(2nm + \frac{3}{2}m^2)$ & $\mathcal{O}(n(m + l))$ \\
   		L-S-BFGS-P (\cref{eq:compform}) & $ \mathcal{O}(n(2lm + 3m^2 + $ & $\mathcal{O}(4nm + 3m^2)$ & $\mathcal{O}(n(3m + l))$ \\
   		& $ 4m + 1) + m^3)^{\dagger}$ & &
   	\end{tabular}
   	\end{table}
   	Note that when $ m \ll n$ and $ l \ll n $ L-BFGS, L-S-BFGS-M and L-S-BFGS-P
   	enable computations with complexity lower than $ n^2 $ and therefore allow for large values of $n$.}
	\jjbc{Moreover, Table 1 shows that the proposed limited-memory BFGS methods have similar search direction
	complexity to unstructured L-BFGS, but higher update cost.}
	
\section{Numerical Experiments}
\label{sec:numex}
This section describes the numerical experiments for the proposed methods in
Section \ref{sec:alg}. The numerical experiments are carried out in MATLAB 2016a on a 
MacBook Pro @2.6 GHz Intel Core i7, with 32 GB of memory. The experiments are 
divided into \jjbc{five} parts.
%With the first experiment,
In Experiment I, we investigate initialization strategies. Experiment II
compares the limited memory methods with the full-memory methods.
The tests in this experiment are on
the same 61 CUTEst \cite{GouOT03} problems as in \cite{PetraChiangAnitescu19}, unless
otherwise noted. \jjbc{In Experiment III, we use classification data from LIBSVM 
(a library for support vector machines \cite{Chang11}) in order to solve regularized logistic regression 
problems with the proposed methods. We also include an application to PDE constrained optimization.}
%In Experiment III, the proposed methods are used in \jjbc{one} structured problem applications. 
%(\jjbc{)  %The second application is an optimal control problem from PDE-constrained optimization.
In Experiment IV, the proposed methods and L-BFGS with   %{Accessed 21 June 2018}. % \cite{ZhuByrdNocedal97}
the IPOPT \cite{WaechterBiegler06} solver are compared. \jjbc{Experiment V, describes 
a real world application from image reconstruction.}
%In summary, the experiments are as follows:
%\begin{description}
%	\item[Experiment I:] Initialization strategies 
%		\begin{enumerate}
%			\renewcommand{\theenumi}{\Alph{enumi}}
%			\item L-S-BFGS-M
%			\item L-S-BFGS-P
%			\item Eigenvalue estimation of $ \sigma_k $
%		\end{enumerate}
%	\item[Experiment II:] Memory choice
%		\begin{enumerate}
%			\renewcommand{\theenumi}{\Alph{enumi}}
%			\item L-S-BFGS-M
%			\item L-S-BFGS-P
%		\end{enumerate}	
%	\item[Experiment III:] Applications
%		\begin{enumerate}
%			\renewcommand{\theenumi}{\Alph{enumi}}
%			\item Logistic regressions
%			\item Optimal control
%		\end{enumerate}
%	\item[Experiment IV:] Comparisons with L-BFGS-B, IPOPT
%\end{description}

Extended performance profiles as in \cite{MahajanLeyfferKirches11} are provided. These profiles are an extension
of the well known profiles of Dolan and Mor\'{e} \cite{DolanMore02}. 
 We compare the number of iterations and the total computational time for each solver on the test set of problems, 
 \jjbc{unless otherwise stated}. 
 The performance metric $ \rho_s(\tau) $ with a given number of test problems $ n_p $ is
\begin{equation*}
		\rho_s(\tau) = \frac{\text{card}\left\{ p : \pi_{p,s} \le \tau \right\}}{n_p} \quad \text{and} \quad \pi_{p,s} = \frac{t_{p,s}}{ \underset{1\le i \le S, i \ne s}{\text{ min } t_{p,i}} },
\end{equation*} 
where $ t_{p,s}$ is the ``output'' (i.e., iterations or time) of
``solver'' $s$ on problem $p$. Here $ S $ denotes the total number of solvers for a given comparison. This metric measures
the proportion of how close a given solver is to the best result. The extended performance profiles are the same as the classical ones for $ \tau \ge 1 $.
In the profiles we include a dashed vertical grey line, to indicate $ \tau = 1 $. In all experiments the line search parameters are
set to $ c_1 = 1\times10^{-4} $ and $ c_2 = 0.9 $.

%The convergence tolerance is set to $ \epsilon = 10^{-6} $.

\subsection{Experiment I}
\label{subsec:EX_I}
This experiment investigates the initialization strategies from Section \ref{sec:alg}. To this end, the problems in
this experiment are not meant to be overly challenging, yet they are meant to enable some variations. Therefore,
we define the quadratic functions
\begin{equation*}
	%\label{eq:quad_func}
	Q_i(\b{x}; \phi, r) \equiv \frac{1}{2}\b{x}^T( \phi \cdot \b{I} + \b{Q}_i \b{D}_i \b{Q}_i^T )\b{x},
\end{equation*} 
with scalar parameters $  0 < \phi  $, $ 1 \le r \le n $ and where $ \b{D}_i \in \mathbb{R}^{r \times r} $ is a diagonal
matrix and $ \b{Q}_i \in \mathbb{R}^{n \times r} $ has orthonormal gaussian columns. Note that $ r $ eigenvalues
of the Hessian $ \nabla^2 Q_i$ are the diagonal elements of $ \phi \cdot \b{I} + \b{D}_i $, while the remaining $ (n-r) $ 
eigenvalues are $ \phi  $. Therefore, by varying $ \phi, r $, and the elements of $ \b{D}_i $, Hessian matrices with
different spectral properties are formed. In particular, when $ r \ll n $, the eigenvalues are clustered around $ \phi $.
In the experiments of this section we investigate \jjbc{ $ \phi = 1 $. In Appendix A, we include tests
when $ \phi = 1000 $. } 
%two values are investigated; specifically, $ \phi = \{1,1000 \} $.  
The structured objective functions from 
\cref{eq:intro_strucmin} are defined by
\begin{equation}
\label{eq:struc_quad}
	\widehat{k}(\b{x}) = \b{x}^T \b{g} + Q_1(\b{x};\phi,r), \quad \quad \widehat{u}(\b{x}) = Q_2(\b{x};\phi,r).
\end{equation}
We refer to the objective functions $ f(\b{x}) = \widehat{k}(\b{x}) + \widehat{u}(\b{x})  $ defined by \cref{eq:struc_quad}
as \emph{structured quadratics}. The problems in this experiment have dimensions $ n = j \cdot 100 $ with corresponding
$ r = j \cdot 10 $ for $ 1 \le j \le 7 $. Since some of the problem data in this experiment is randomly generated (e.g., the 
orthonormal matrices $ \b{Q}_i $), the experiments are repeated five times for each $n$. 
The reported results are of the average values of the five individual runs. For all solvers we set $ m = 8 $ (memory parameter),
$ \epsilon = 5\times 10^{-6} $ ($ \|\bk{g}\|_{\infty} \le \epsilon $), and maximum iterations to 10,000 \jjbc{This limit was not 
reached in the experiments}.
  
\subsubsection{Experiment I.A: L-S-BFGS-M}
\label{subsub:EX_I_LSBFGSM}
\begin{figure*}[t!]
	\begin{minipage}{0.48\textwidth}
		\includegraphics[trim=0 0 20 15,clip,width=\textwidth]{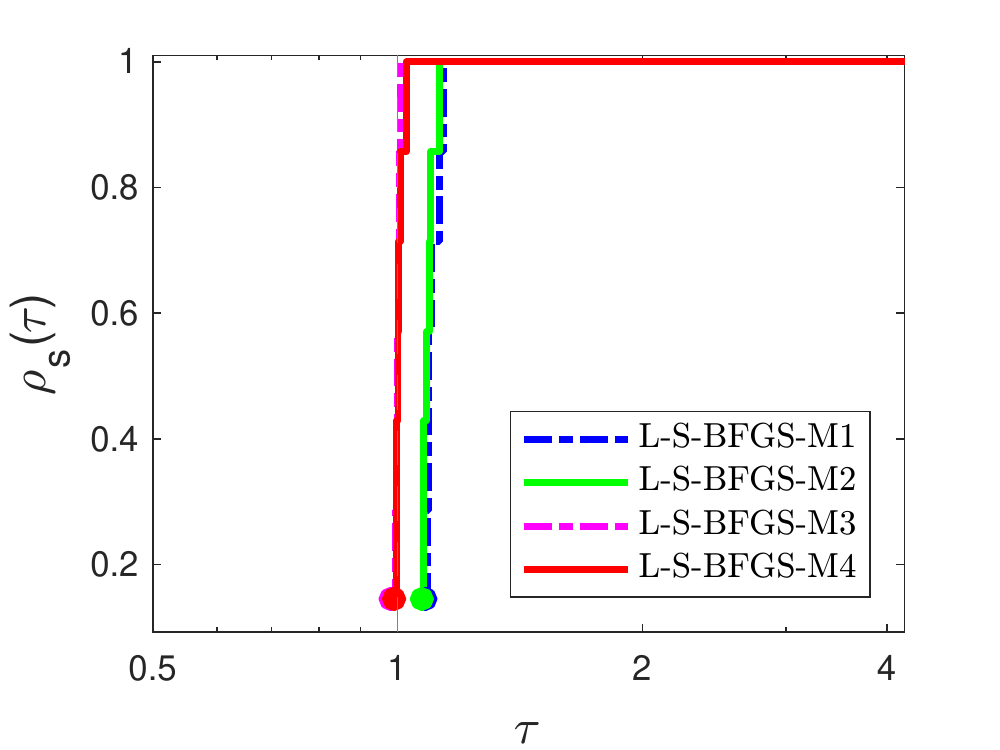}
	\end{minipage}
		\hfill
	\begin{minipage}{0.48\textwidth}
		\includegraphics[trim=0 0 20 15,clip,width=\textwidth]{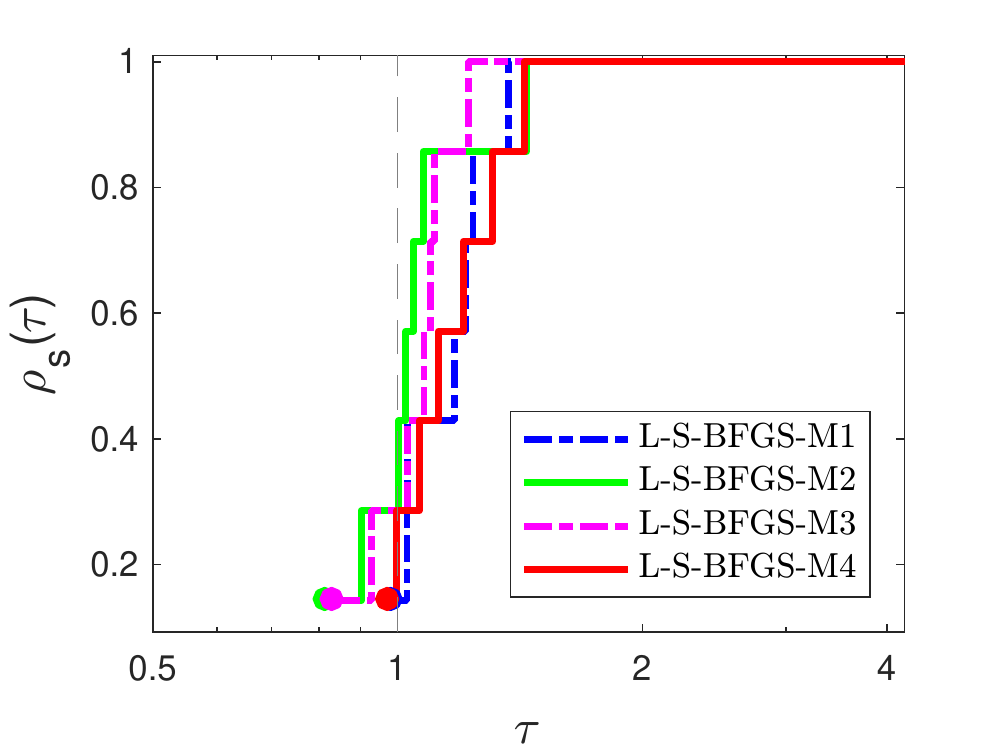}
	\end{minipage}
	\caption{Comparison of initialization strategies for L-S-BFGS-M on problems with eigenvalues clustered around 1 with 
	$1 \le \lambda_r \le 1000 $ and $ \lambda_{r+1} = \cdots = \lambda_{n} = 1 $.  Left: number of iterations; right: time.}
		\label{fig:EX_IA_SM}       
\end{figure*}

Experiment I.A compares the four L-S-BFGS-M initializations on the structured quadratic objective functions with eigenvalues
clustered around 1. In particular, $ \phi = 1 $, and the elements of $ \b{D}_i $ are uniformly distributed  in the interval $[0,999]$.
The results are displayed in Fig. \ref{fig:EX_IA_SM}.
We observe that in terms of number of iterations, Init. 4 (red) and Init. 3 (purple) perform similarly
and that also Init. 2 (green) and Init. 1 (blue) perform similarly. Overall, Init. 4 and 
Init. 3 requirer fewer iterations on the structured quadratics. Moreover, the solid
lines are above the dashed ones for both pairs. This indicates that including only gradient information in $ \bhk{u} $ and in the initialization strategy, as opposed to also including 2$^{\text{nd}}$ derivative information from $ \bk{u} $, may be desirable for this problem. Init. 1 and Init. 2 are fastest on these problems. Even though these initializations require a larger number of iterations, they can be faster because the line searches terminate more quickly. 

%\subsubsection{Experiment I.B: L-S-BFGS-P}
%\label{subsub:EX_I_LSBFGSB}
\jjbc{Next, we compare the four  L-S-BFGS-P initializations. As before, experiments on problems with eigenvalues clustered at 1 are done.}
%Experiment I.B compares the four L-S-BFGS-P initializations. As in Section \ref{subsub:EX_I_LSBFGSM} 
%experiments on problems with eigenvalues clustered at 1 are done. 
\jjbc{Experiments with eigenvalues clustered around 1000 are included in Appendix A}. The respective outcomes 
are in Figure \ref{fig:EX_IA_SM_P}. % and \ref{fig:EX_IA_LRG_P}.

\begin{figure*}[t!]
	\begin{minipage}{0.48\textwidth}
		\includegraphics[trim=0 0 20 15,clip,width=\textwidth]{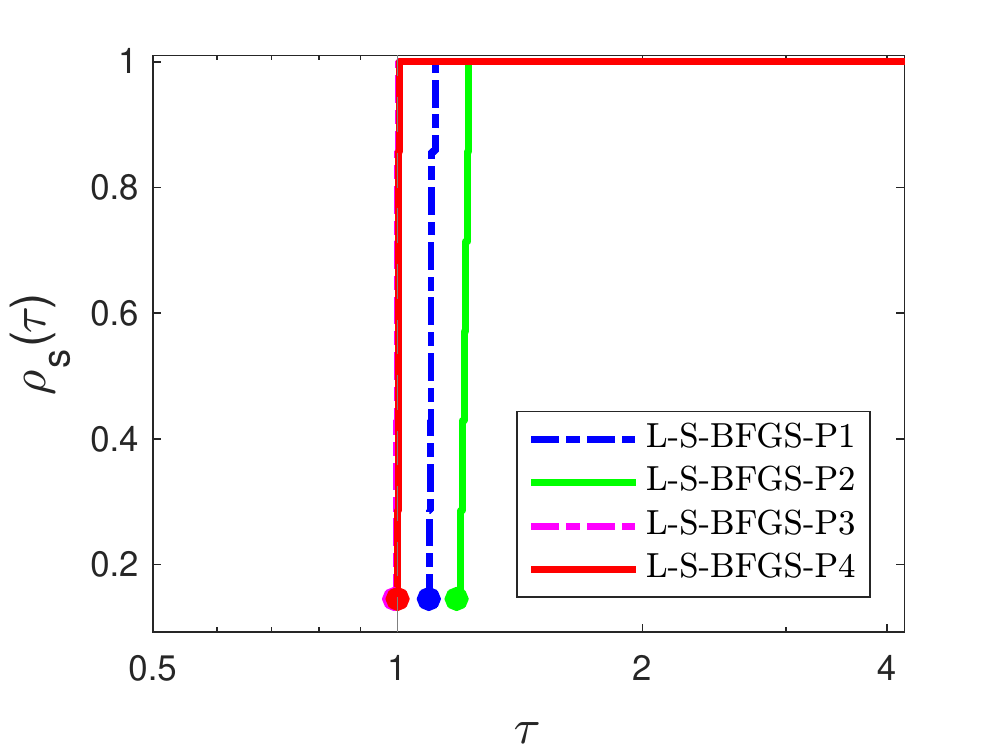}
	\end{minipage}
		\hfill
	\begin{minipage}{0.48\textwidth}
		\includegraphics[trim=0 0 20 15,clip,width=\textwidth]{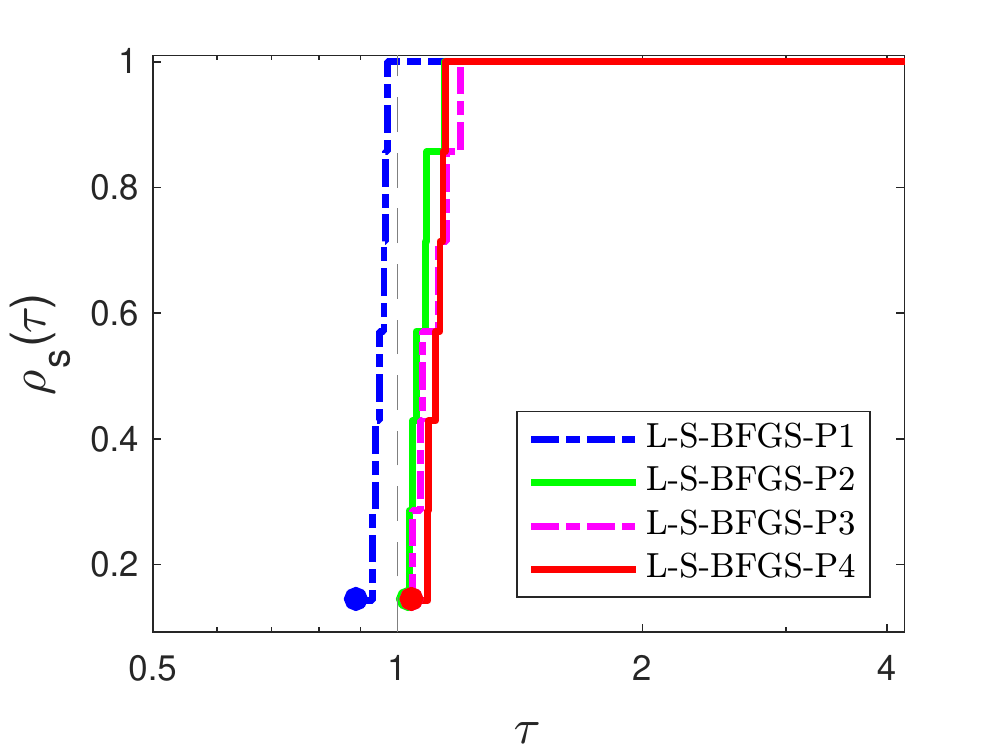}
	\end{minipage}
	\caption{Comparison of initialization strategies for L-S-BFGS-P on problems with eigenvalues clustered around 1 with 
	$1 \le \lambda_r \le 1000 $ and $ \lambda_{r+1} = \cdots = \lambda_{n} = 1 $.  Left: number of iterations; right: time.}
		\label{fig:EX_IA_SM_P}       
\end{figure*}
We observe that, similar to Figure \ref{fig:EX_IA_SM}, Init. 3 and Init. 4 do best in iterations, while Init. 1 does best in time. %In Figure \ref{fig:EX_IA_SM_P} 
%Experiment on large clustered eigenvalues.

To analyze the properties of the scaling factor $ \sigma_k $ in greater detail, Section \ref{subsub:EX_I_EIG} describes experiments
that relate $ \sigma_k $ to eigenvalues.

\subsubsection{Experiment I.B: Eigenvalue Estimation}
\begin{figure*}[t!]
	\begin{minipage}{0.48\textwidth}
		\includegraphics[trim=0 0 30 0,clip,width=\textwidth]{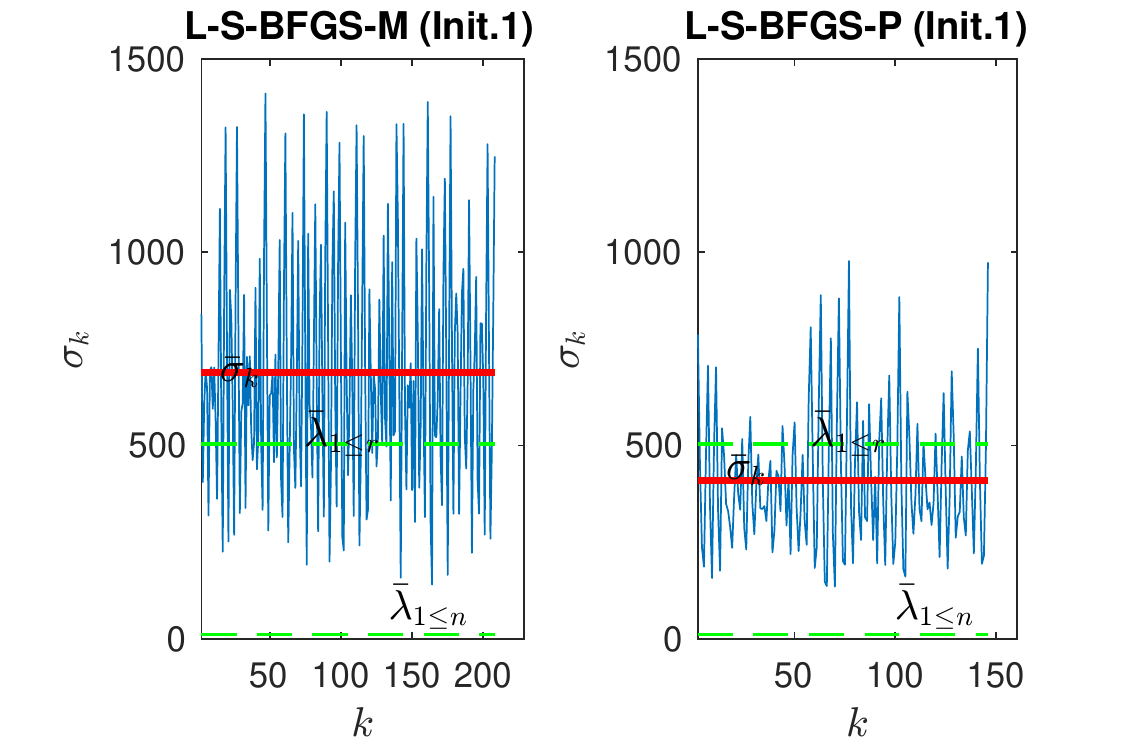}
	\end{minipage}
		\hfill
	\begin{minipage}{0.48\textwidth}
		\includegraphics[trim=0 0 20 0,clip,width=\textwidth]{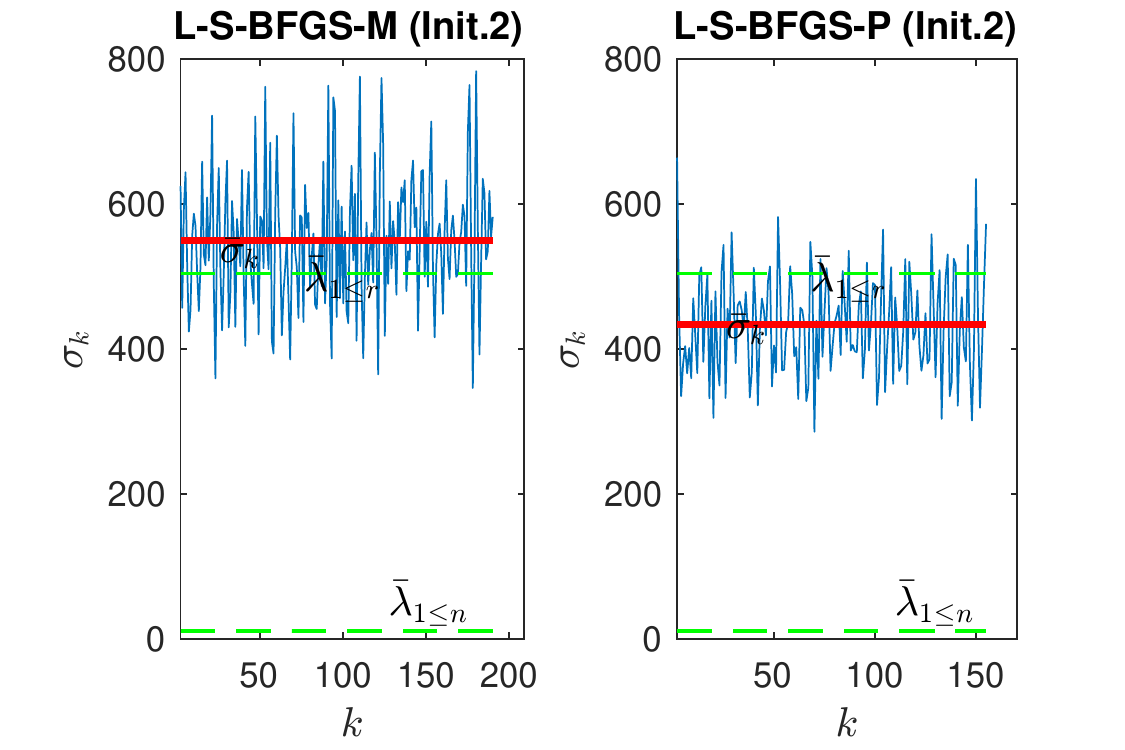}
	\end{minipage}
	\begin{minipage}{0.48\textwidth}
		\includegraphics[trim=0 0 30 0,clip,width=\textwidth]{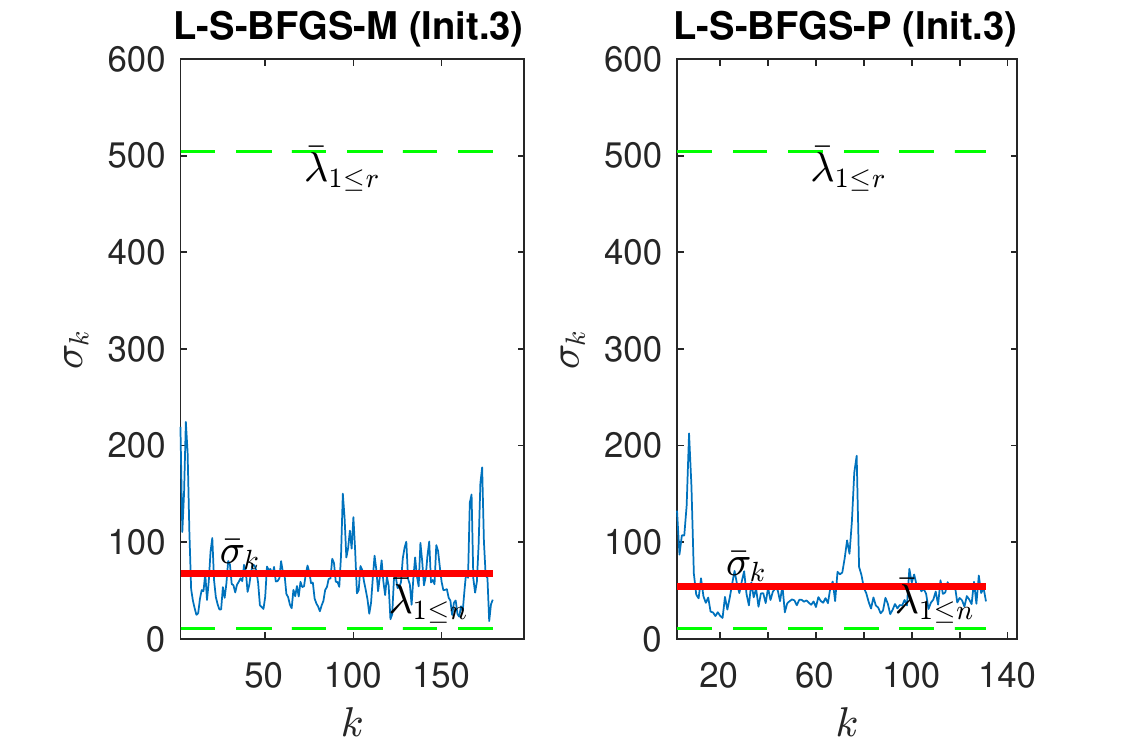}
	\end{minipage}
		\hfill
	\begin{minipage}{0.48\textwidth}
		\includegraphics[trim=0 0 20 0,clip,width=\textwidth]{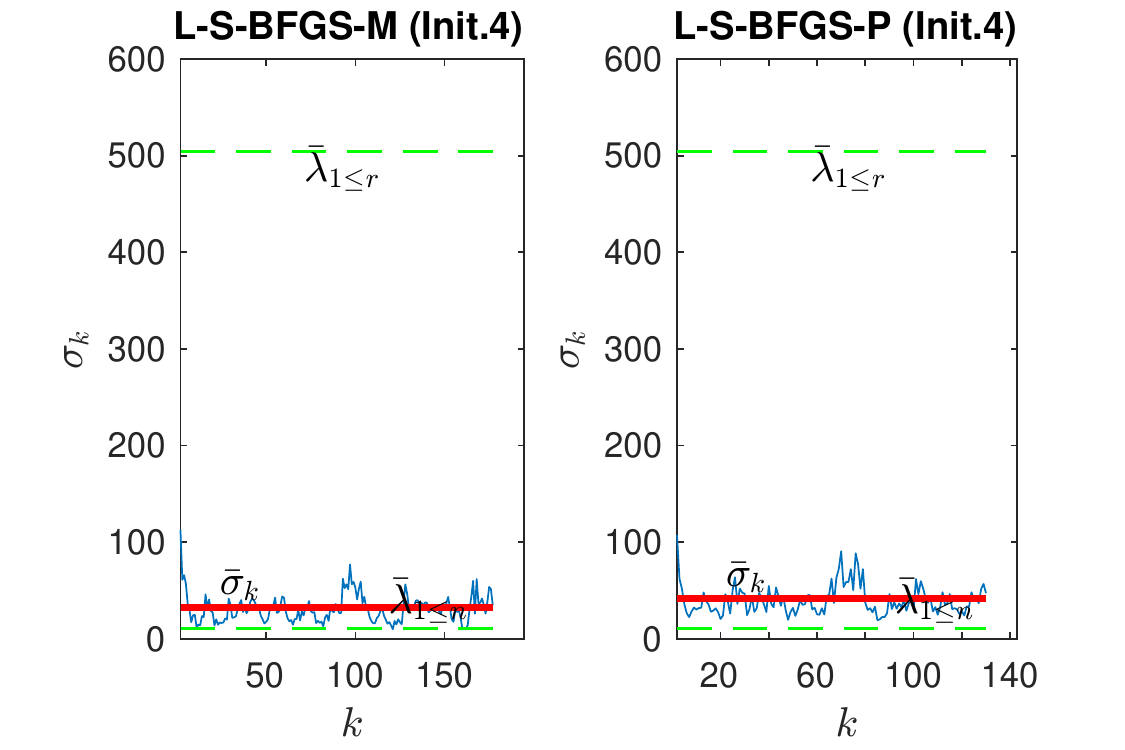}
	\end{minipage}
	\caption{Eigenvalue estimation with initialization parameter $ \sigma_k $. The eigenvalues are clustered around 1 with 
	$1 \le \lambda_r \le 1000 $ and $ \lambda_{r+1} = \cdots = \lambda_{n} = 1 $. }
		\label{fig:EX_IC_SM}       
\end{figure*}

\begin{figure*}[t!]
	\begin{minipage}{0.48\textwidth}
		\includegraphics[trim=0 0 30 0,clip,width=\textwidth]{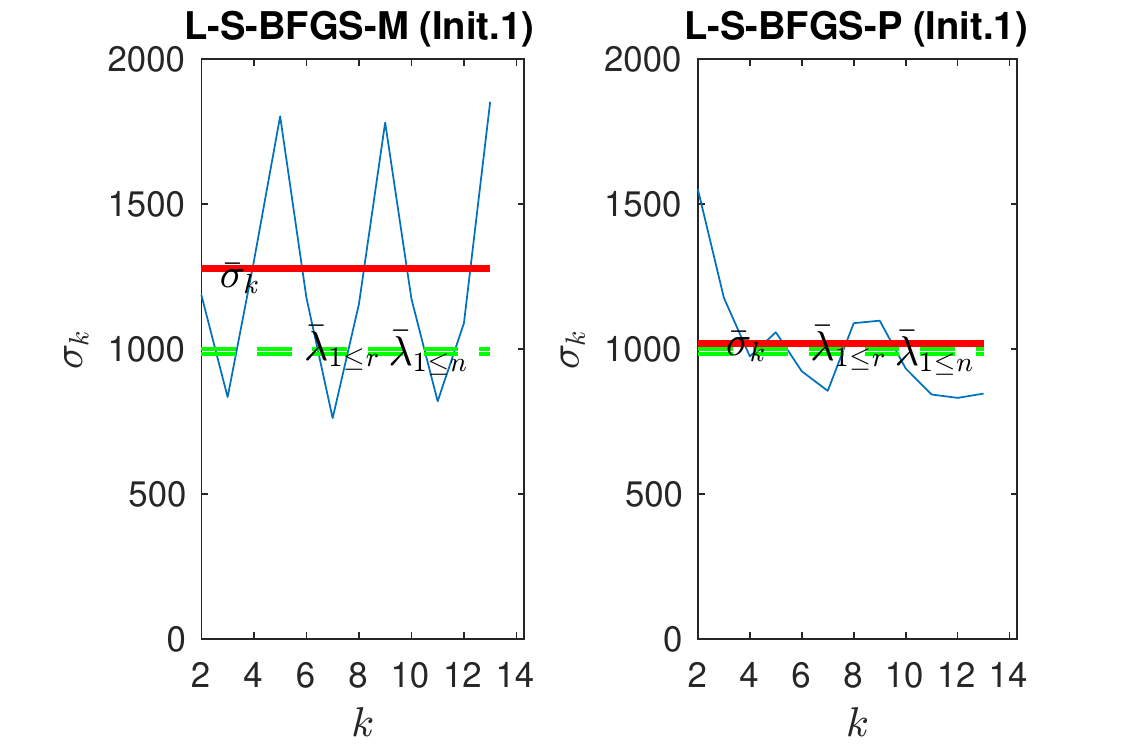}
	\end{minipage}
		\hfill
	\begin{minipage}{0.48\textwidth}
		\includegraphics[trim=0 0 20 0,clip,width=\textwidth]{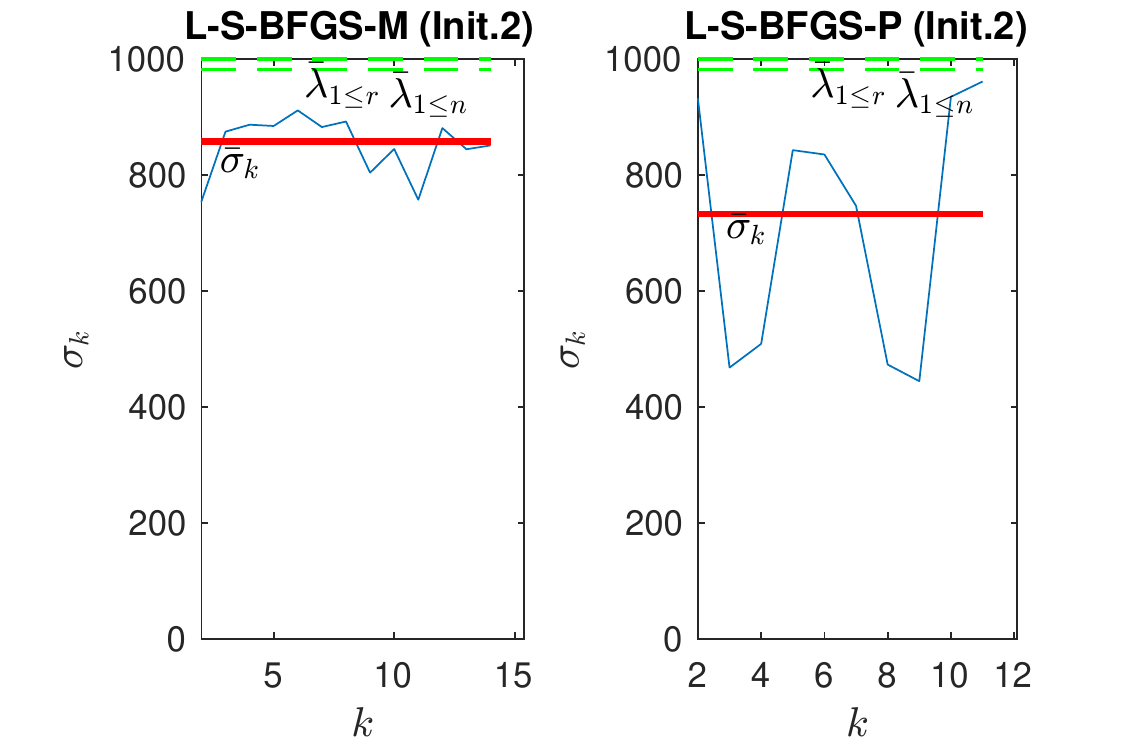}
	\end{minipage}
	\begin{minipage}{0.48\textwidth}
		\includegraphics[trim=0 0 30 0,clip,width=\textwidth]{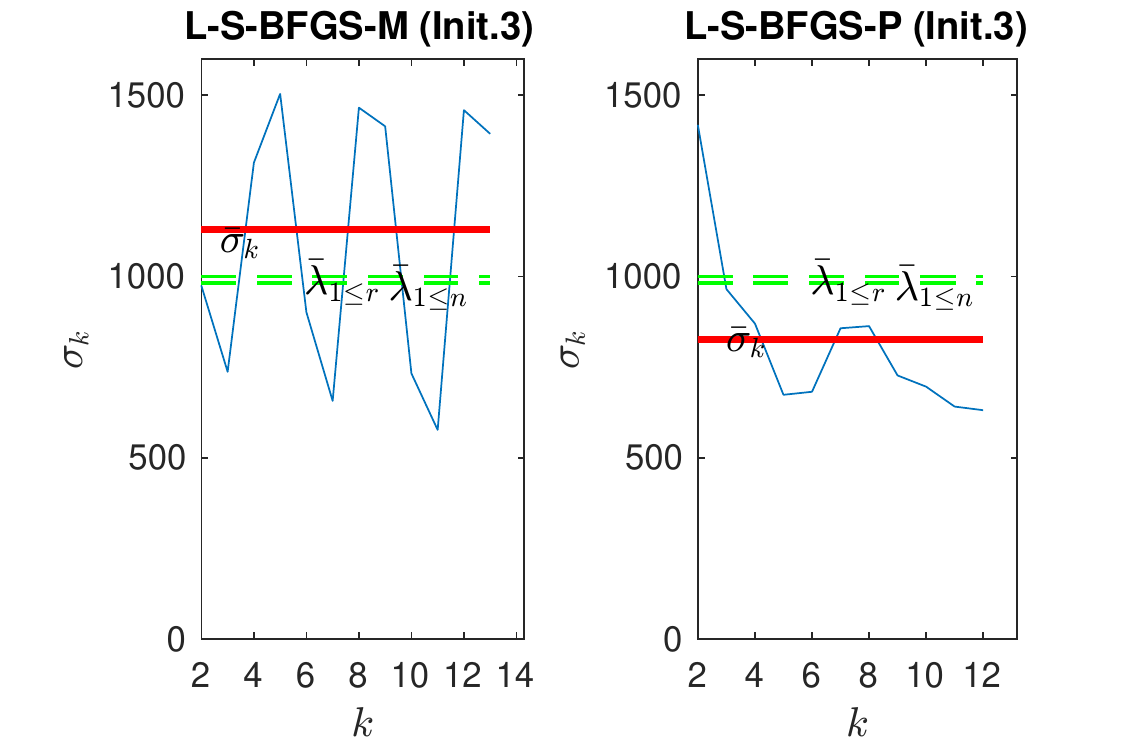}
	\end{minipage}
		\hfill
	\begin{minipage}{0.48\textwidth}
		\includegraphics[trim=0 0 20 0,clip,width=\textwidth]{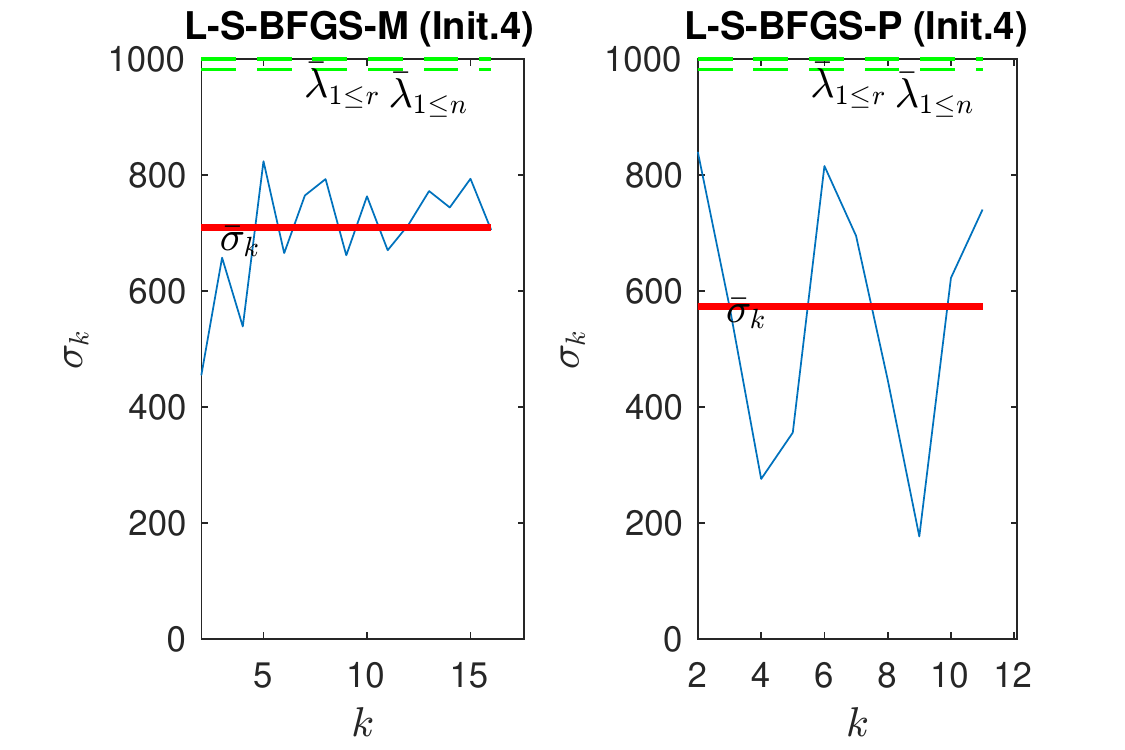}
	\end{minipage}
	\caption{Eigenvalue estimation with scaling parameter. The eigenvalues are clustered around 1,000 with 
	$1 \le \lambda_r \le 1000 $ and $ \lambda_{r+1} = \cdots = \lambda_{n} = 1000 $. }
		\label{fig:EX_IC_LRG}       
\end{figure*}

In Experiment I.B we investigate the dynamics of $ \sigma_k $ in the four initialization
strategies from \cref{eq:initializations} on a fixed problem as the iteration count $k$ increases. In particular, we use one representative run from the average results of the preceding two subsections, where $ n =100 $ and $ r = 10 $. In Figure \ref{fig:EX_IC_SM} the evolution of $ \sigma_k $ of all four initializations for both;
L-S-BFGS-M and L-S-BFGS-P is displayed on a structured quadratic problem with eigenvalues clustered at 1. In Figure \ref{fig:EX_IC_LRG} the same quantities are displayed for structured quadratic problems with eigenvalues clustered at 1000. In green $ \bar{\lambda}_{1 \le n} $ and $ \bar{\lambda}_{1 \le r} $ are displayed, which correspond to the median taken over the first $ 1,2,\cdots, n $ (all) and the first $ 1,2,\cdots, r $ eigenvalues, respectively. Because in Figure \ref{fig:EX_IC_SM} the eigenvalues are clustered around 1, $ \bar{\lambda}_{1 \le n} = 1 $. In Figure \ref{fig:EX_IC_LRG} the eigenvalues are clustered around 1000 and $ \bar{\lambda}_{1 \le r} = 1000 $. In red $ \bar{\sigma}_k $ is the average $ \sigma_k $ value over all iterations.

\label{subsub:EX_I_EIG}

Across all plots in Figures \ref{fig:EX_IC_SM} and \ref{fig:EX_IC_LRG} we observe that the dynamics of
$\sigma_k$ for L-S-BFGS-M and L-S-BFGS-P are similar. Moreover, the average $ \bar{\sigma}_k $ is higher for Init. 1 and Init. 2 than for Init. 3 and Init. 4. The variability of Init. 2 appears less than that of Init. 1, while the variability of Init. 4 appears less than that of Init. 3. We observe that Init. 1 and 2 approximate a large eigenvalue well, whereas Init. 3 and Init. 4 approximate smaller eigenvalues better (cf. Figure \ref{fig:EX_IC_SM} lower half). Since large $ \sigma_k $ values typically result in shorter step lengths (step computations use $ 1 / \sigma_k $), choosing Init. 1 and Init. 2 result in shorter step lengths on average. \jjb{Taking shorter average steps can be a desirable conservative strategy when the approximation to
the full Hessian matrix is not very accurate. Therefore as a general guideline, Init. 1 and Init. 2 appear more suited for problems in which
it is difficult to approximate the Hessian accurately, and Init. 1 and Init. 2 are more suited for problems in which larger step sizes are desirable.}

%\begin{figure}[tbhp]
%	\centering
%	\includegraphics{figs/SM_SIM1_n3_Ini1}
%	\includegraphics{figs/SM_SIM1_n3_Ini2}
%	\includegraphics{figs/SM_SIM1_n3_Ini3}
%	\includegraphics{figs/SM_SIM1_n3_Ini4}
%	%\subfloat[Iterations]{\label{fig:num_csbfgsm8its}\includegraphics{figs/its_EXM_m8}}
%	%\subfloat[Times]{\label{fig:num_csbfgsm8times}\includegraphics{figs/times_EXM_m8}}
%	\caption{Eigenvalue estimation with scaling parameter.}
%	\label{fig:scale_sm}
%\end{figure}

\subsection{Experiment II}
\label{subsec:EX_II}

Experiment II compares the limited memory structured formulas with the full-memory update formulas from Petra et al. 
\jjbc{\cite{PetraChiangAnitescu19} on the CUTEst problems from \cite{PetraChiangAnitescu19}}. The full-memory algorithms from \cite{PetraChiangAnitescu19},
which use Eqs. \cref{eq:intro_sbfgsm} and \eqref{eq:intro_sbfgsp}, are called S-BFGS-M and S-BFGS-P, 
respectively. The line search procedures of the limited memory structured BFGS algorithms
(Algorithms \ref{alg:alg_csbfgsm} and \ref{alg:alg_csbfgsp}) are the same as for the full memory algorithms.
Moreover, the initializations in the full memory algorithms are set as $ \b{A}^{\text{M}}_0 = \bar{\sigma} \In $
for S-BFGS-M, and $ \b{A}^{\text{P}}_0 = \bar{\sigma} \In $ for S-BFGS-P, where
$ \bar{\sigma} = 10^{i} $ for the first $ i \ge 0 $ that satisfies $ (10^{i} \In + \bz{K}) \succ 0 $ (usually $i=0$). 
%These initializations
%ensure that, when $ k \le m $, the algorithms L-BFGS-M1 and L-BFGS-P1 generate the 
%same iterates as S-BFGS-M and S-BFGS-P. Only after $ k > 0 $ do iterates of these
%limited memory versions and the full-recursive updates differ. The other 
%versions from \cref{tbl:num_versions}, however, generally produce iterates
%that are different from the full-recursive algorithms.
The experiments are divided into two main parts. Experiment II.A. tests the
limited memory structured BFGS-Minus versions corresponding to
Algorithm \ref{alg:alg_csbfgsm}. Experiment II.A. is further 
subdivided into the cases in which the memory parameters are $ m = 8 $
and $ m = 50 $. \jjb{These values represent a typical value ($ m = 8$) and a relatively large value ($ m = 50$), cf. e.g., \cite{BoggsByrd19}}. Experiment II.B. tests the
limited memory structured BFGS-Plus versions corresponding to
Algorithm \ref{alg:alg_csbfgsp}. As before, Experiment II.B. is further 
subdivided into the cases in which the memory parameters are $ m = 8 $
and $ m = 50 $. For all the solvers, we set
$ \epsilon = 1\times 10^{-6} $ ($ \|\bk{g}\|_{\infty} \le \epsilon $) and maximum iterations to 1,000.

\subsubsection{Experiment II.A: L-S-BFGS-M}
\label{subsubsec:EX_II_A}

\begin{figure*}[t!]
	\begin{minipage}{0.48\textwidth}
		\includegraphics[trim=0 0 20 15,clip,width=\textwidth]{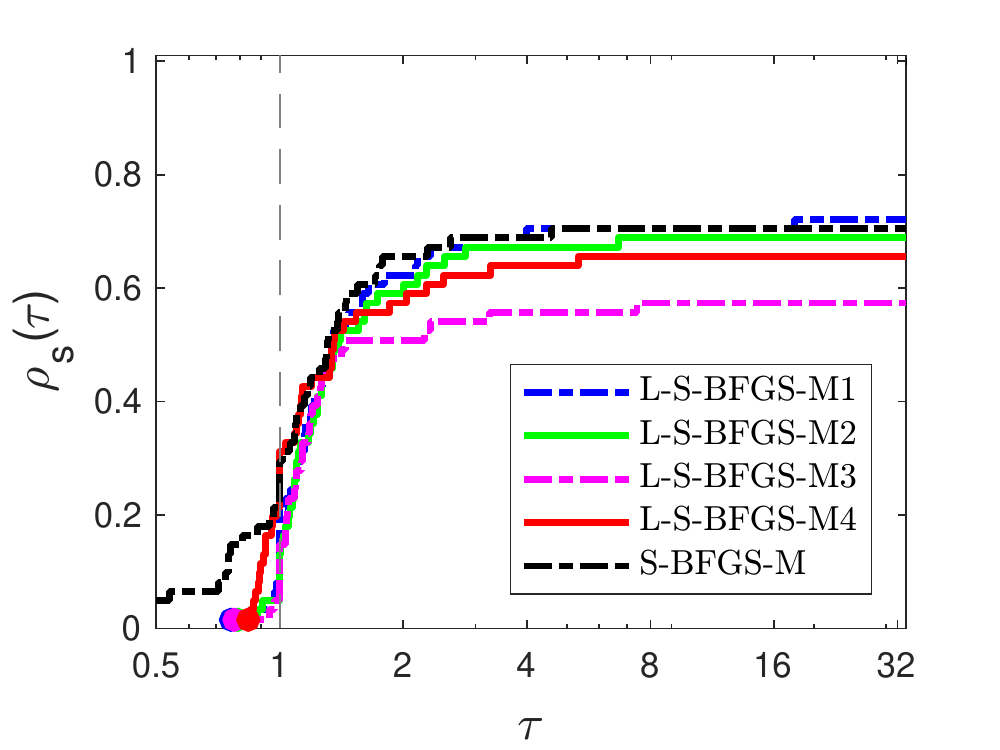}
	\end{minipage}
		\hfill
	\begin{minipage}{0.48\textwidth}
		\includegraphics[trim=0 0 20 15,clip,width=\textwidth]{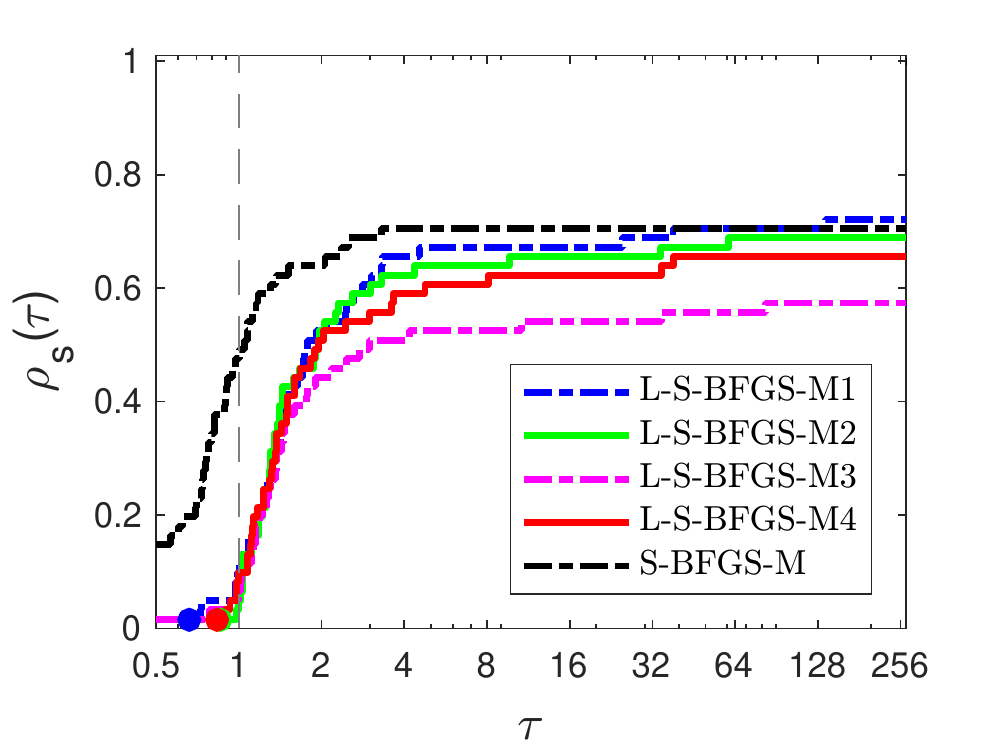}
	\end{minipage}
	\caption{Comparison of four initialization strategies of L-S-BFGS-M from \cref{eq:initializations} to the full-recursive method S-BFGS-M (corresponding to \cref{eq:intro_sbfgsm}) on all 62 CUTEst problems from \cite{PetraChiangAnitescu19}. The limited memory parameter is $ m = 8 $. Left:
	number of iterations; right: time.}
		\label{fig:EX_II_A_m8}       
\end{figure*}

In Experiment II.A we compare the limited memory implementations of Algorithm \ref{alg:alg_csbfgsm} with initialization strategies in
\cref{eq:initializations} with the full-recursive S-BFGS-M method from \cref{eq:intro_sbfgsm}. 
The solvers are tested on all 62 CUTEst problems from \cite{PetraChiangAnitescu19}. Figure \ref{fig:EX_II_A_m8} contains
the results for the limited memory parameter $ m =8 $.

%\begin{figure}[tbhp]
%	\centering
%	\includegraphics{figs/its_EXM_m8}
%	%\subfloat[Iterations]{\label{fig:num_csbfgsm8its}\includegraphics{figs/its_EXM_m8}}
%	%\subfloat[Times]{\label{fig:num_csbfgsm8times}\includegraphics{figs/times_EXM_m8}}
%	\caption{Comparison of three versions of L-S-BFGS-M from \cref{tbl:num_versions} to the full-recursive method
%	 S-BFGS-M (corresponding to \cref{eq:intro_sbfgsm}). The limited memory parameter is $ m =8 $.}
%	\label{fig:num_csbfgsm8}
%\end{figure}

We observe that the full-memory S-BFGS-M (black) does well in terms of number of iterations
and execution time. However, L-S-BFGS-M1 (Init. 1, blue), a limited memory version with memory
of only $ m =8 $, does comparatively well. In particular, this strategy is able
to solve one more problem, as indicated by the stair step at the right end of the plot.

%L-S-BFGS-M3 (green) obtains the overall best results in this experiment. 
%L-S-BFGS-M3 uses the nonconstant initialization 
%$ \bs{\Psi}_0(\widehat{\sigma}_k,0) = \widehat{\sigma}_k \In $, whereas the full-recursive
%method (blue) is based on a  constant initialization $ \bs{\Psi}_0 = \bar{\sigma} \In + \bz{K} $.
%Version L-S-BFGS-M1 (dashed-red), also uses a constant initialization 
%$ \bs{\Psi}_0(\bar{\sigma},1) = \bar{\sigma} \In + \bz{K} $, and the iterations it generates 
%are theoretically the same as those of S-BFGS-M for $ k \le m $. For $ k > m $ the limited memory
%updating begins, and subsequently the quasi-Newton approximation in 
%L-S-BFGS-M1 are expected to become less accurate. Overall, L-S-BFGS-M2, which
%uses an initialization that combines a constant and an updated component, i.e., 
%$ \bs{\Psi}_0(\widehat{\sigma}_k,1) = \widehat{\sigma}_k \In + \bz{K} $, requires the most iterations. 

\begin{figure*}[t!]
	\begin{minipage}{0.48\textwidth}
		\includegraphics[trim=0 0 20 15,clip,width=\textwidth]{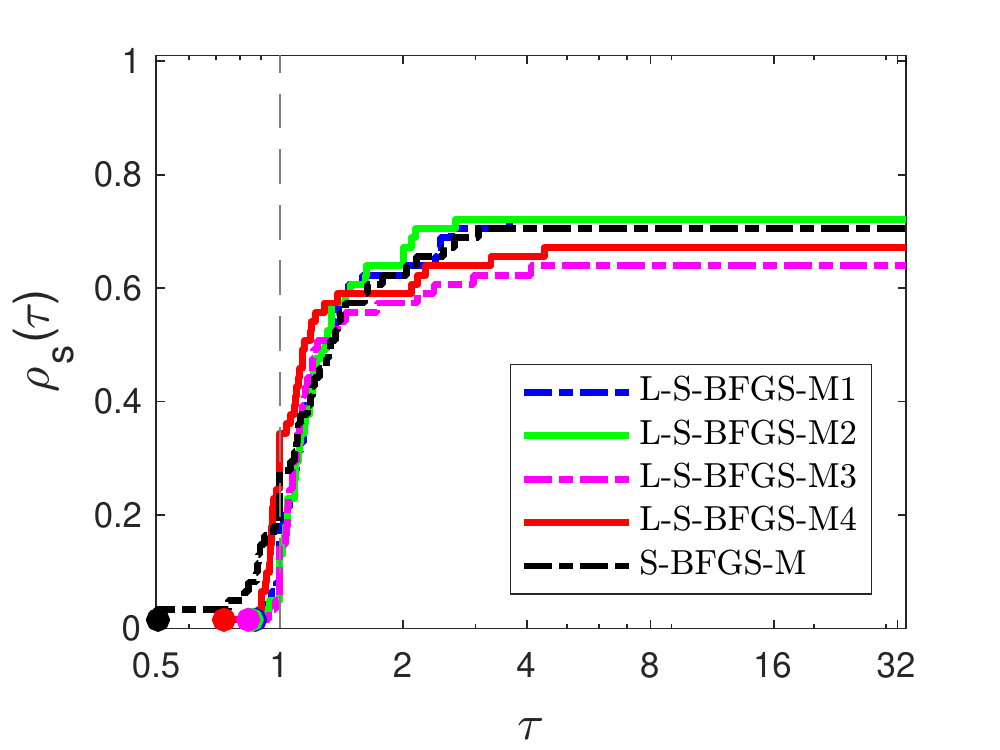}
	\end{minipage}
		\hfill
	\begin{minipage}{0.48\textwidth}
		\includegraphics[trim=0 0 20 15,clip,width=\textwidth]{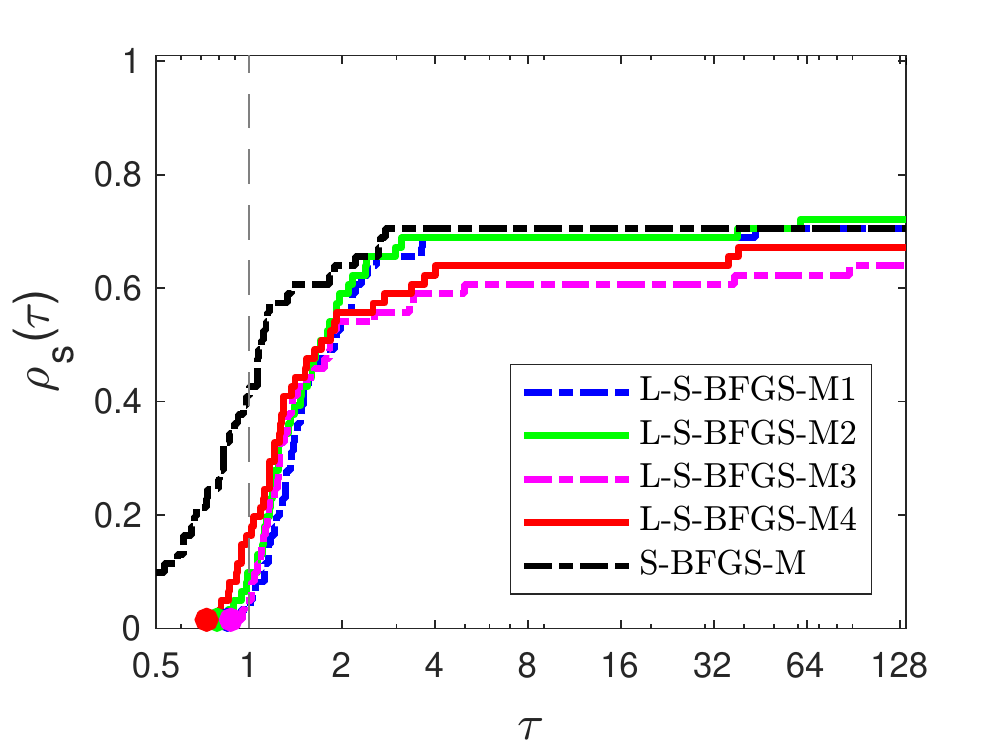}
	\end{minipage}
	\caption{Comparison of four initialization strategies of L-S-BFGS-M from \cref{eq:initializations} with the full-recursive method S-BFGS-M (corresponding to \cref{eq:intro_sbfgsm}) on all 62 CUTEst problems from \cite{PetraChiangAnitescu19}. The limited memory parameter is $ m = 50 $. Left: number of iterations; right: time.}
		\label{fig:EX_II_A_m50}       
\end{figure*}

Figure \ref{fig:EX_II_A_m50} shows
the results for the limited memory parameter $ m = 50 $. A larger
limited memory parameter makes using limited memory
structured matrices more computationally expensive but is also expected
to increase the accuracy of the quasi-Newton approximations.

%\begin{figure}

%\begin{figure}[tbhp]
%	\centering
%	\includegraphics{figs/its_EXM_m50}
%	%\subfloat[Iterations]{\label{fig:num_csbfgsm8its}\includegraphics{figs/its_EXM_m8}}
%	%\subfloat[Times]{\label{fig:num_csbfgsm8times}\includegraphics{figs/times_EXM_m8}}
%	\caption{Comparison of three versions of L-S-BFGS-M from \cref{tbl:num_versions} to the full-recursive method
%	 S-BFGS-M (corresponding to \cref{eq:intro_sbfgsm}). The limited memory parameter is $ m = 50 $.}
%	\label{fig:num_csbfgsm50}
%\end{figure}

Note that the outcomes of S-BFGS-M (black) in Figure \ref{fig:EX_II_A_m50}
are the same as those in Figure \ref{fig:EX_II_A_m8}, because it does not depend on  
the memory parameter. For the limited memory versions we observe that the 
outcomes of L-S-BFGS-M2 (green) improve notably, whereas the other limited memory
versions remain roughly unchanged. Using the initialization strategies (Init. 1 or Init. 2), limited memory solvers are able to solve one more problem than the full-memory method can,
as indicated by the highest ending lines in the plot. 
We suggest that Init. 1 and Init. 2 (see Section \ref{subsub:EX_I_EIG}) generate initialization parameters $ \sigma_k  $ that are on average larger than those generated by Init. 3 or Init. 4. These larger values in turn result in shorter average step sizes, which appears advantageous on general
nonlinear problems.
%\jjbc{Remarkably, we conclude that even though the computational costs with the limited memory compact representations are}
% The improvements are possibly a result of more accurate quasi-Newton
%approximations. Overall, however, L-S-BFGS-M3 (green) still appears the 
%most robust.

\subsubsection{Experiment II.B: L-S-BFGS-P}
\label{subsubsec:EX_II_B}

In Experiment II.B we compare the versions of Algorithm \ref{alg:alg_csbfgsp} using the initialization strategies from \cref{eq:initializations} with the 
full memory recursive S-BFGS-P method \cref{eq:intro_sbfgsp}. The solvers
are run on 55 of the 62 CUTEst problems from \cite{PetraChiangAnitescu19} for which $ n \le 2500 $. Figure \ref{fig:EX_II_B_m8} contains
the results for the limited memory parameter $ m =8 $:

% trim left bottom right top
\begin{figure*}[t!]
	\begin{minipage}{0.48\textwidth}
		\includegraphics[trim=0 0 20 15,clip,width=\textwidth]{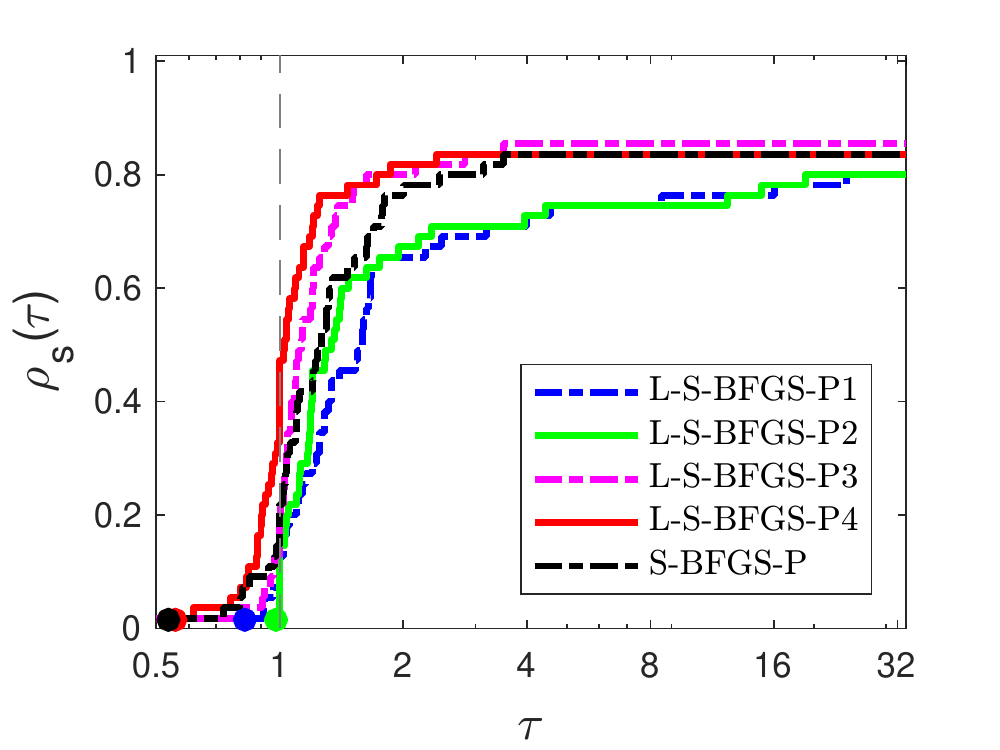}
	\end{minipage}
		\hfill
	\begin{minipage}{0.48\textwidth}
		\includegraphics[trim=0 0 20 15,clip,width=\textwidth]{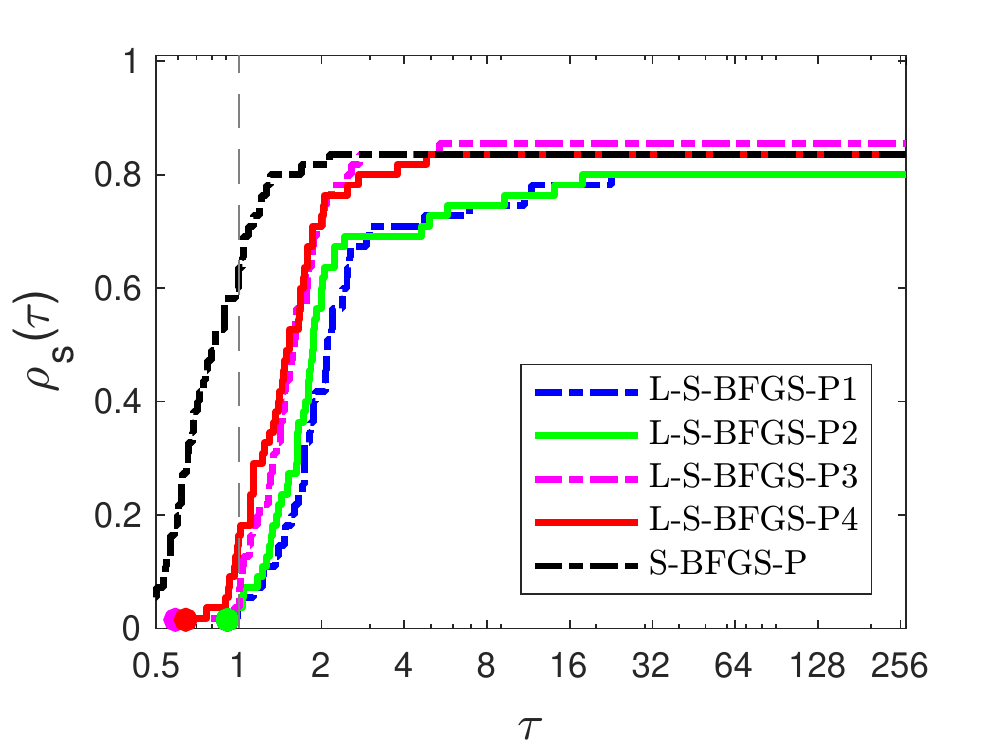}
	\end{minipage}
	\caption{Comparison of four initialization strategies of L-S-BFGS-P from \cref{eq:initializations} to the full-recursive method S-BFGS-P (corresponding to \cref{eq:intro_sbfgsp}) on 55 CUTEst problems from \cite{PetraChiangAnitescu19}. The limited memory parameter is $ m = 8 $. Left: number of iterations; right: time.}
		\label{fig:EX_II_B_m8}       
\end{figure*}

%\begin{figure}[tbhp]
%	\centering
%	\includegraphics{figs/its_EXP_m8_RED}
%	%\subfloat[Iterations]{\label{fig:num_csbfgsm8its}\includegraphics{figs/its_EXM_m8}}
%	%\subfloat[Times]{\label{fig:num_csbfgsm8times}\includegraphics{figs/times_EXM_m8}}
%	\caption{Comparison of two versions of L-S-BFGS-P from \cref{tbl:num_versions} to the full-recursive method
%	 S-BFGS-P (corresponding to \cref{eq:intro_sbfgsp}). The limited memory parameter is $ m = 8 $.}
%	\label{fig:num_csbfgsp8}
%\end{figure}

We observe that for a relatively small memory parameter $ m = 8 $, L-S-BFGS-M3 (Init. 3, purple) solves the most problems. L-S-BFGS-M4 (Init. 4, red) requires the fewest iterations, as indicated by the highest circle on the y-axis in the left panel of Figure \ref{fig:EX_II_B_m8}. 

%the limited memory structured 
%BFGS Plus versions (red and purple) obtain worse outcomes than the full memory
%structured BFGS-Plus method (blue). The L-S-BFGS-P1 version uses the 
%same initialization as the full memory method. Thus the first $ k \le m $ iterates
%of L-S-BFGS-P1 are theoretically equivalent to the one from S-BFGS-P. However
%when $ k > m $ the iterates differ, because of the limited memory updating. 

\begin{figure*}[t!]
	\begin{minipage}{0.48\textwidth}
		\includegraphics[trim=0 0 20 15,clip,width=\textwidth]{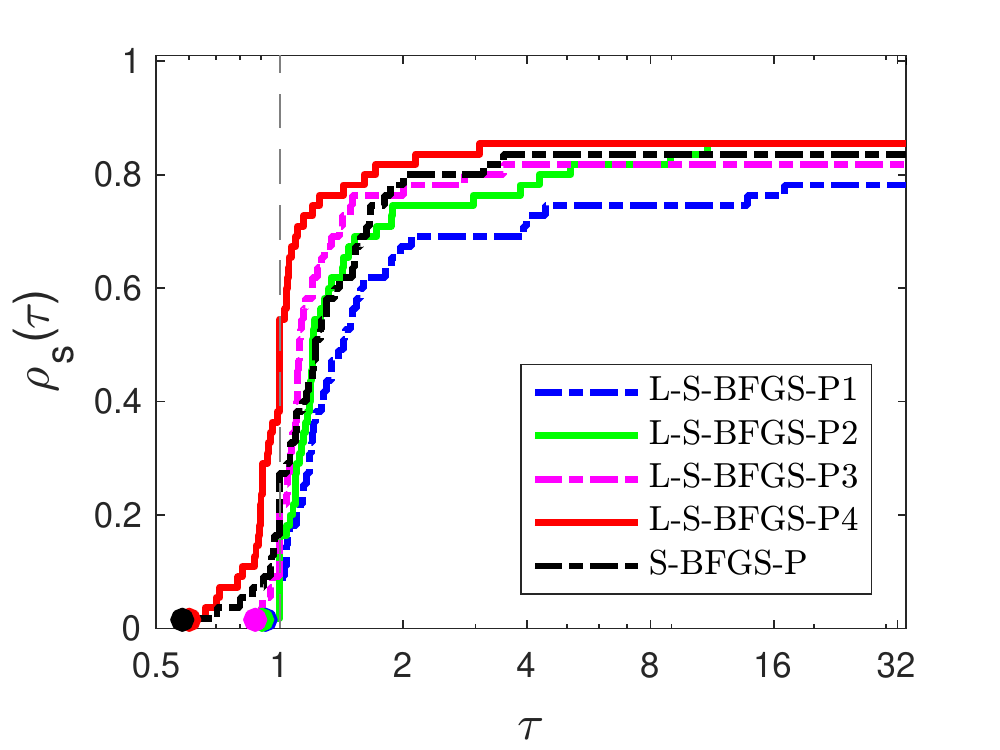}
	\end{minipage}
		\hfill
	\begin{minipage}{0.48\textwidth}
		\includegraphics[trim=0 0 20 15,clip,width=\textwidth]{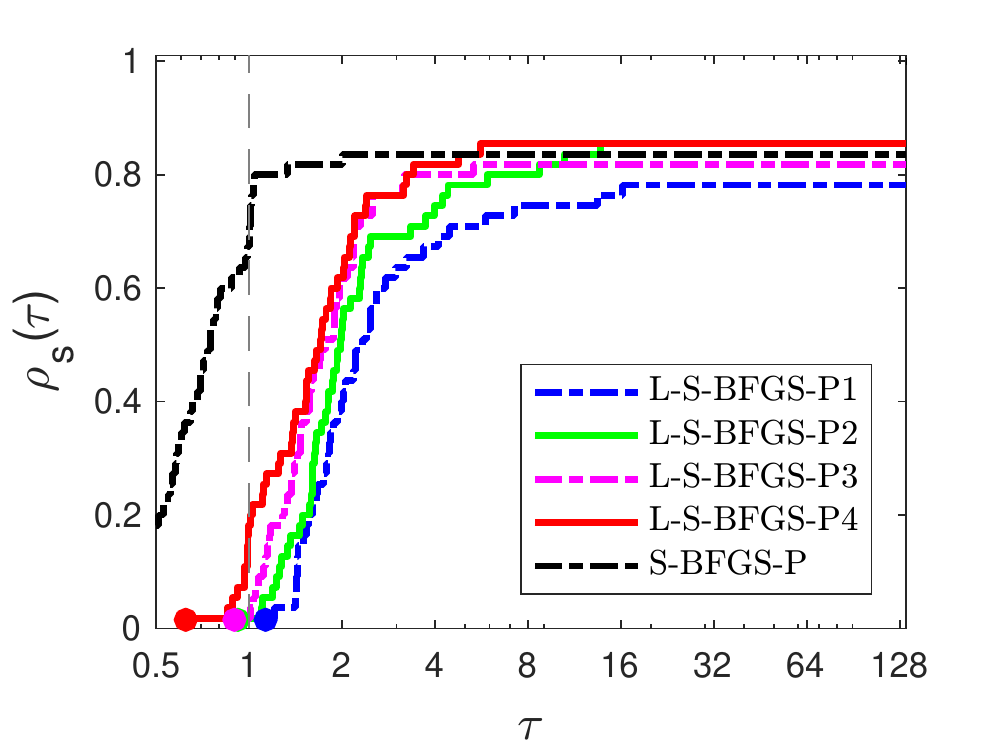}
	\end{minipage}
	\caption{Comparison of four initialization strategies of L-S-BFGS-P from \cref{eq:initializations} to the full-recursive method S-BFGS-P (corresponding to \cref{eq:intro_sbfgsp}) on 55 CUTEst problems from \cite{PetraChiangAnitescu19}. The limited memory parameter is $ m = 50 $. Left: number of iterations; right: time.}
		\label{fig:EX_II_B_m50}       
\end{figure*}

Figure \ref{fig:EX_II_B_m50} shows
the results for the limited memory parameter $ m = 50 $. A larger parameter makes using limited memory
structured matrices more computationally expensive but is also expected
to increase the accuracy of the quasi-Newton approximations.

%\begin{figure}[tbhp]
%	\centering
%	\includegraphics{figs/its_EXP_m50_RED}
%	%\subfloat[Iterations]{\label{fig:num_csbfgsm8its}\includegraphics{figs/its_EXM_m8}}
%	%\subfloat[Times]{\label{fig:num_csbfgsm8times}\includegraphics{figs/times_EXM_m8}}
%	\caption{Comparison of two versions of L-S-BFGS-P from \cref{tbl:num_versions} to the full-recursive method
%	 S-BFGS-P (corresponding to \cref{eq:intro_sbfgsp}). The limited memory parameter is $ m = 50 $.}
%	\label{fig:num_csbfgsp50}
%\end{figure}

Note that the outcomes of S-BFGS-P in
Figure \ref{fig:EX_II_B_m50}
%GAIL - note I changed this -- is this ok?\cref{fig:EX_II_B_m8}
are the same as in Figure \ref{fig:EX_II_B_m8}, because the
full-memory solver does not depend on the memory parameter. For a larger memory $ m = 50 $, the outcomes of L-S-BFGS-P2 (green) and L-S-BFGS-P4 (red) improve notably. Overall,
L-S-BFGS-P4 solves the most problems.

From the experiments in this section, we find that initialization strategies Init.1 and
Init. 2 appear most desirable for L-S-BFGS-M, whereas Init. 4 and Init. 2 appear most desirable for L-S-BFGS-P.
%the limited memory structured 
%BFGS Plus versions (red and purple) improve notably. Overall,
%L-S-BFGS-P2 (purple), which uses a non-constant initialization, solves the most problems. 
%The outcomes of L-S-BFGS-P1 (purple) and S-BFGS-P (blue) become more similar, partially
%a consequence of the fact that the first $ k \le m $ iterates of both methods are theoretically
%the same.

\subsection{Experiment III}
This section describes \jjbc{one} application of the methods in the context of machine learning. \jjbc{A $2^{\text{nd}}$ similar application to PDE constrained optimization is included, too}. For all solvers we set $ m = 8 $ (memory parameter), % in Appendix B
$ \epsilon = 1\times 10^{-6} $ ($ \|\bk{g}\|_{\infty} \le \epsilon $) and maximum iterations to 10,000. \jjb{Since some of the problems
in this section are large we use the techniques described in Section \ref{subsec:large_comp} throughout the experiments.}  
\jjbc{Because some of the problems in this experiment are very large, the recursive formulas from \eqref{eq:intro_sbfgsm} and \eqref{eq:intro_sbfgsp} (with $ m = \infty $) cannot be directly used 
on these problems. However, the limited memory compact representations use the memory parameter $ m $ to threshold the 
computational and memory cost and are therefore applicable.}

%\subsubsection{Experiment III: Logistic Regressions} % .A
%label{subsub:EX_III_A}
\subsubsection{Experiment III.A: Logistic Regressions} % .A
\begin{figure*}[t!]
	\begin{minipage}{0.48\textwidth}
		\includegraphics[trim=0 0 20 15,clip,width=\textwidth]{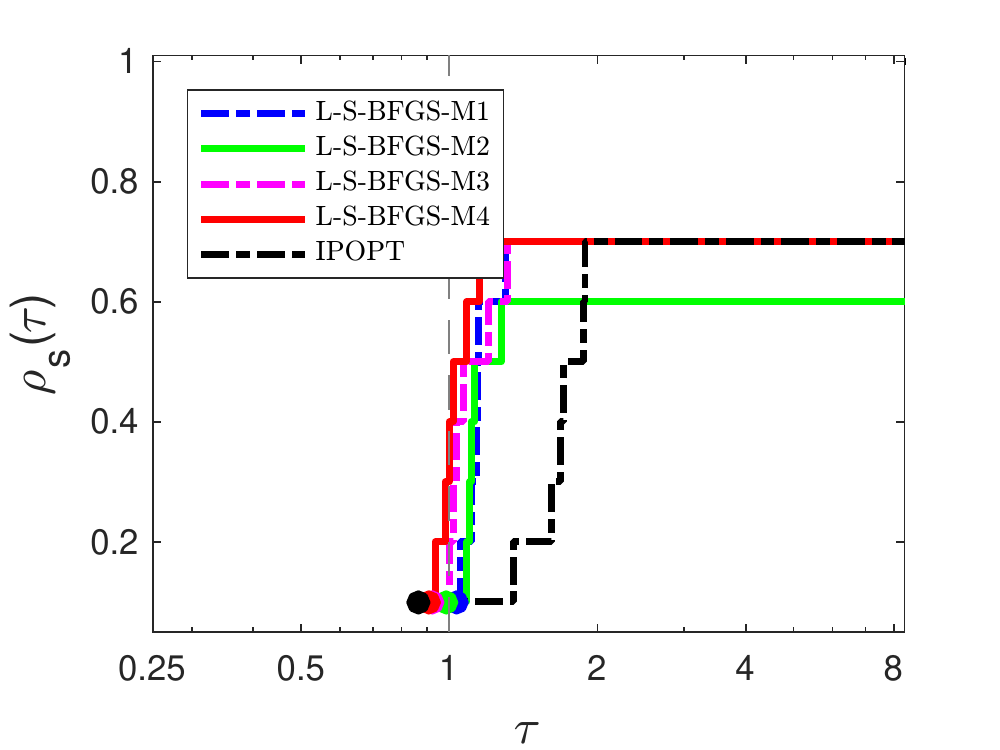}
	\end{minipage}
		\hfill
	\begin{minipage}{0.48\textwidth}
		\includegraphics[trim=0 0 20 15,clip,width=\textwidth]{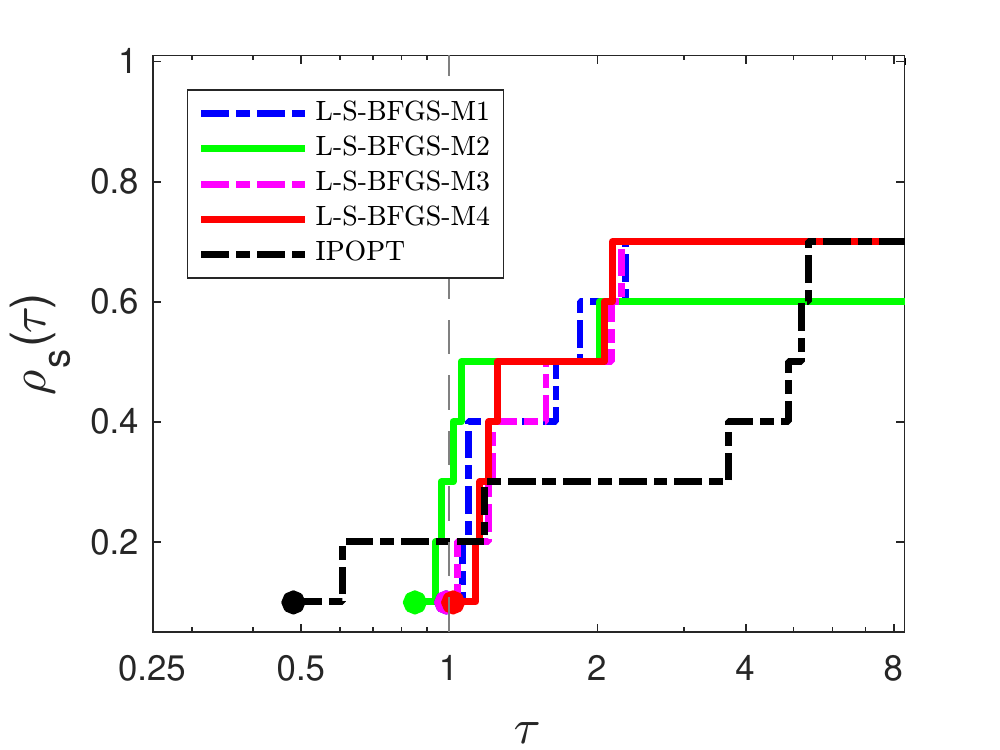}
	\end{minipage}
	\caption{Comparison of L-S-BFGS-M solvers on 10 logistic regression
	classification problems using data from LIBSVM. Left: number of iterations,
	right: time.}
		\label{fig:EX_III_A}       
\end{figure*}

%Experiment III tests the proposed methods on
The problems in this section are defined by smooth-structured objective functions % .A
from machine learning, as described, for example, in \cite{SraNowozinWright11}. In particular,
logistic regression problems use smooth objective functions for classification tasks (for instance, \cite{ByrdChinNeveittNocedal11}), which
often depend on a large number of data points and many variables. The classification problems are defined by the data pairs $ \{ \b{d}_i, y_i \}_{i=1}^{D} $,
where the so-called feature vectors $ \b{d}_{i} \in \mathbb{R}^n $ may be large, 
and the so-called labels $ y_i \in \{-1,1\} $ are scalars. In \cite{TeoVishwanthanSmolaLe10} regularized logistic regression problems are described
in which the objective function is composed of two terms. The optimization problems are
formulated as
\begin{equation*}
	\underset{\b{x} \in \mathbb{R}^n}{ \text{ minimize } } 
	\frac{\lambda}{2} \| \b{x} \|^2_2 +
	\sum_{i=1}^{D} \log\left(1+ \exp( - y_i \b{x}^T \b{d}_i ) \right), 
\end{equation*}
where $ \lambda > 0 $. The regularization term,  $ \frac{\lambda}{2} \| \b{x} \|^2 $,
has a second derivative,$\lambda \b{I}$, that is readily available. Therefore, we define the known and unknown components for this problem as
\begin{equation}
	\label{eq:ex_III_A_ku}
	\widehat{k}(\b{x}) = \frac{\lambda}{2} \| \b{x} \|^2_2, 
	\quad \quad \widehat{u}(\b{x}) = \sum_{i=1}^{D} \log\left(1+ \exp( - y_i \b{x}^T \b{d}_i ) \right).
\end{equation}
This data was obtained from \texttt{www.csie.ntu.edu.tw/\textasciitilde{}cjlin/libsvm/} (retrieved on 10/03/19). Ten problems were used, with problem dimensions listed in Table 2. % \ref{tbl:libsvm} 
% experiment's 
\begin{table}[h!]
\caption{List of dimensions for 10 LIBSVM logistic regression problems. Here $D$ denotes the number of training
pairs $\{\b{d}_i, y_i \}_{i=1}^D$, and $ n $ denotes the number of variables/feature weights (the size of the problem).}
\centering
\begin{tabular}{ | c | c | c | }
\hline \textbf{Problem} 		& $ D $ 		& $ n $ \\
\hline \text{rcv1} 				& 20242 		 	 		& 47236 \\
\text{duke} 					& 34 				 			 & 7129 \\
\text{gisette} 					& 6000 				 			& 5000 \\
\text{colon\_cancer} 			& 62 				  			& 2000 \\
\text{leukemia} 				& 38 								& 7129 \\
\text{real\_sim} 				& 72309 							& 20958 \\
\text{madelon} 					& 2000 								& 500 \\
\text{w8a} 						& 49749 			 		& 300 \\
\text{mushrooms} 				& 2000 						& 500 \\
\text{a9a} 						& 32561 					& 123 \\
 \hline  
\end{tabular}
\label{tbl:libsvm}
\end{table}
Some of the problems are large, with $ n \ge 5000 $ and \jjbc{thus we focus on the computations as described in Section \ref{subsec:large_comp_m}}. The regularization parameter is set as $ \lambda = 10^{-3} $. \jjbc{For comparison, we include
IPOPT \cite{WaechterBiegler06} with a L-BFGS quasi-Newton matrix % L-BFGS and 
(we use a precompiled Mex file with IPOPT 3.12.12, MUMPS and MA57). %{Accessed 21 June 2018}.  \cite{ZhuByrdNocedal97}  
We specify the limited memory BFGS option for IPOPT using the setting 
\texttt{hessian\_approximation= `limited memory'} and tolerances by 
\texttt{tol=9.5e-10} and 
\texttt{acceptable\_tol = 9.5e-10}.}
The results of the experiments are shown in Figure \ref{fig:EX_III_A}. \jjbc{We observe that all solvers, except for L-S-BFGS-M2,
solve the same total number of problems. Moreover, the structured L-BFGS solvers tend to use fewer
iterations and overall less computational time than IPOPT's L-BFGS method.}
%\jjb{Because the Hessian of $ \widehat{k}(\b{x}) $ is constant, the search directions of L-S-BFGS-M and L-S-BFGS-P are the same, when the computations with L-S-BFGS-M are done as described in \Cref{subsec:large_comp_m}. Thus we focus on presenting the results with this method.} 

%\begin{table}
%\label{tbl:libsvm}
%\begin{tabular}{ | c | c | c | c | }
%\hline \textbf{Problem} 		& \multicolumn{2}{c}{$ D $} 		& $ n $ \\
%\hline							& \text{train} 	& \text{test} 	& \\
%\hline \text{rcv1} 				& 20242 		& 677399 		& 47236 \\
%\text{duke} 					& 34 			& 4 			& 7129 \\
%\text{gisette} 					& 6000 			& 1000 			& 5000 \\
%\text{colon\_cancer} 			& 62 			& NA 			& 2000 \\
%\text{leukemia} 				& 38 			& 34 			& 7129 \\
%\text{real\_sim} 				& 72309 		& NA 			& 20958 \\
%\text{madelon} 					& 2000 			& 600 			& 500 \\
%\text{w8a} 						& 49749 		& 14951 		& 300 \\
%\text{mushrooms} 				& 2000 			&	 600 		& 500 \\
%\text{a9a} 						& 32561 		& 16281 		& 123 \\
% \hline  
%\end{tabular}
%\end{table}

\jjbc{Next, we describe experiments for optimal control problems
with similar structures.}

%In this experiment L-S-BFGS-M1 (blue) appears the most robust. This
%solver is based on Init. 1, which typically yields the largest values 
%of $ \sigma_k $ (see \cref{subsub:EX_I_EIG}). Overall, we observe that
%the dashed lines ,
%%GAIL - should you say what the dashed lines represent
%which typically produce larger values of $ \sigma_k  $,
%do best on these problems. Larger values of $ \sigma_k $ typically result
%in shorter steps (see e.g., Line \ref{alg:csbfgsm_p} in \cref{alg:alg_csbfgsm}
%where $ \b{H}_0 = (1/\sigma_k) \b{I} $). Therefore when, as in this case,
%the Hessian of $ \widehat{u}(\b{x}) $ may be difficult to accurately approximate,
%because it depends on data, it appears advantageous to use initialization
%strategies, which produce relatively large $ \sigma_k $.

\subsubsection{Experiment III.B: Optimal Control Problems}
\label{subsub:EX_III_B}
This experiment describes a typical situation in PDE constrained optimization.
In particular, if the PDE is nonlinear, then we can compute gradients efficiently
using the adjoint equation, but Hessians of the unknown part cannot be 
computed efficiently.
%Experiment III.B tests the L-S-BFGS-M methods on 
%optimal control problems with a PDE constraint. In particular,
Denoting $ u $ as the horizontal axis and $ v $ as the vertical axis,
then 2D Poisson problems, with an unknown control $ x(u,v) $, are defined
by the differential equation: $ y_{uu} + y_{vv} = x $.
%\begin{equation*}
%	\frac{\partial^2 y }{\partial u^2} + \frac{\partial^2 y }{\partial v^2} = x(u,v).
%\end{equation*}
The solution $ y(u,v) $ has known boundary values on a box $(u,v) \in [0,1]^2 $; in other words,
$ y(0,v)$, $y(1,v)$, $y(u,0)$, and $y(u,1) $ are known. Discretizing the domain and splitting it into
an interior and boundary part, we get for the optimal control problem 
\begin{equation*}
	\underset{\b{x} \in \mathbb{R}^n}{\text{ minimize }}
	\frac{1}{2} \left \{ \| \b{x} \|^2_2  +   \left\|  \b{y}(\b{x}) - \b{y}^* \right\|^2_2 \right \}   \quad
	\text{ subject to } \quad
	\b{A} \b{y} = \b{x} + \b{g},
\end{equation*}

%\begin{alignat*}{2}
%	& \underset{\b{x} \in \mathbb{R}^n}{\text{ minimize }}
%	&& \frac{1}{2} \| \b{x} \|^2_2  +  \frac{1}{2} \left\|  \b{y}(\b{x}) - \b{y}^* \right\|^2_2    \\
%	& \text{ subject to } && \\
%	&  && \b{A} \b{y} = \b{x} 					\quad \text{(Interior)} \\
%	&  && \textcolor{white}{\b{A}}\b{y} = \b{g} 		\quad \text{(Boundary)},
%\end{alignat*}
\noindent where $ \b{g} \in \mathbb{R}^n $ represents a vector with boundary information, 
$ \b{A} \in \mathbb{R}^{n \times n} $ is a matrix resulting from a 5-point stencil finite difference discretization 
of the partial derivatives, and $ \b{y}^* $ are fixed data values. Because the Hessian of the regularization
term,  $ \frac{1}{2} \| \b{x} \|^2_2 $, is straightforward to compute, we define 
the structured objective function by 
\begin{equation}
	\label{eq:ex_III_B_ku}
	\widehat{k}(\b{x}) = \frac{1}{2} \| \b{x} \|^2_2, 
	\quad \quad \widehat{u}(\b{x}) = \frac{1}{2} \left\|  \b{y}(\b{x}) - \b{y}^* \right\|^2_2,
\end{equation}
using $ \b{y}(\b{x}) = \b{A}^{-1} (\b{x} + \b{g}) $. The number of variables is defined
by the formula $ n = (10\times j - 2)^2 $, where $ j = 2,3,\cdots,10 $, which corresponds
to discretizations with $20,30,\cdots,100$ mesh points in one direction. The largest
problem has $ n = 9604 $ variables. For comparison we also include the implementation
of a ``standard" BFGS method from \cite{PetraChiangAnitescu19},
which uses the same line search as do the limited memory structured methods and
IPOPT's L-BFGS method.

\begin{figure*}[t!]
	\begin{minipage}{0.48\textwidth}
		\includegraphics[trim=0 0 20 15,clip,width=\textwidth]{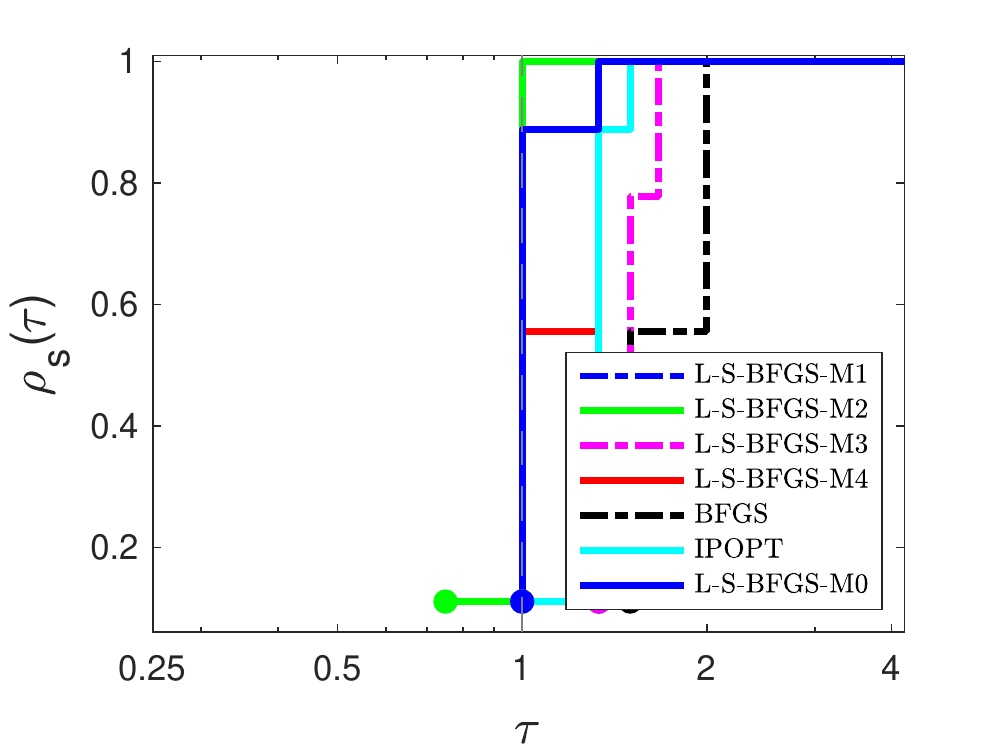}
	\end{minipage}
		\hfill
	\begin{minipage}{0.48\textwidth}
		\includegraphics[trim=0 0 20 15,clip,width=\textwidth]{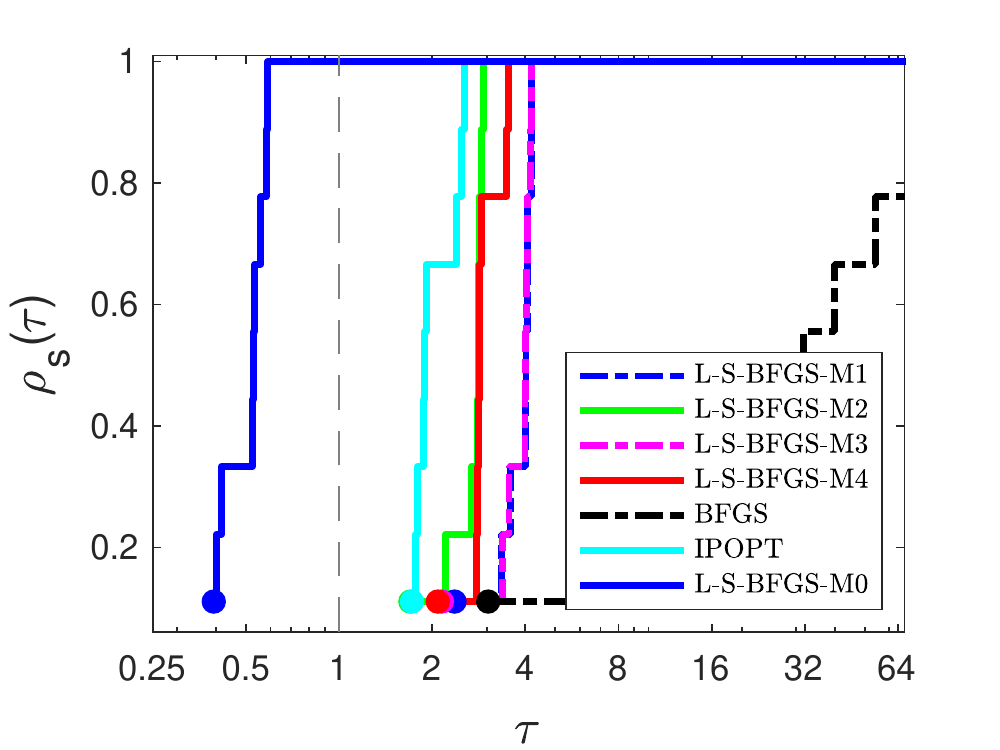}
	\end{minipage}
	\caption{Comparison of L-S-BFGS-M solvers on PDE constrained optimal control problems. The 
	dimensions of the problems are $n = (10\times j - 2)^2 $ for $ j = 2,3,\cdots,10 $. Left: number of iterations,
	right: time. \jjbc{L-S-BFGS-M0 represents an efficient implementation of Algorithm 1 in which the initialization 
	is the constant $ \sigma_k = 0 $.}}
		\label{fig:EX_III_B}       
\end{figure*}

%In this experiment, L-S-BFGS-M3 and L-S-BFGS-M1 (dashed purple and blue) appear to do best overall. Both solvers
%typically generate $ \sigma_k $ values, which are larger than the ones from the solvers with solid lines 
%(cf. observations from \cref{subsub:EX_I_EIG}). Larger values of $ \sigma_k $ yield shorter step lengths,
%which may be regarded as a more conservative strategy. The effect of computing search directions
%by the Sherman-Morrison-Woodbury formula in Line \ref{alg:csbfgsm_p} from \cref{alg:alg_csbfgsm},
%appears in the right-hand plot of \cref{fig:EX_III_B}. In particular, the limited memory versions do not 
%require solves with linear systems that depend on the dimension $ n $, whereas the BFGS implementation
%from \cite{PetraChiangAnitescu19} is based on solves with $ n \times n $ dimensional systems.

\subsection{Experiment IV}
\label{subsec:EX_IV}
\begin{figure*}[t!]
	\begin{minipage}{0.48\textwidth}
		\includegraphics[trim=0 0 20 15,clip,width=\textwidth]{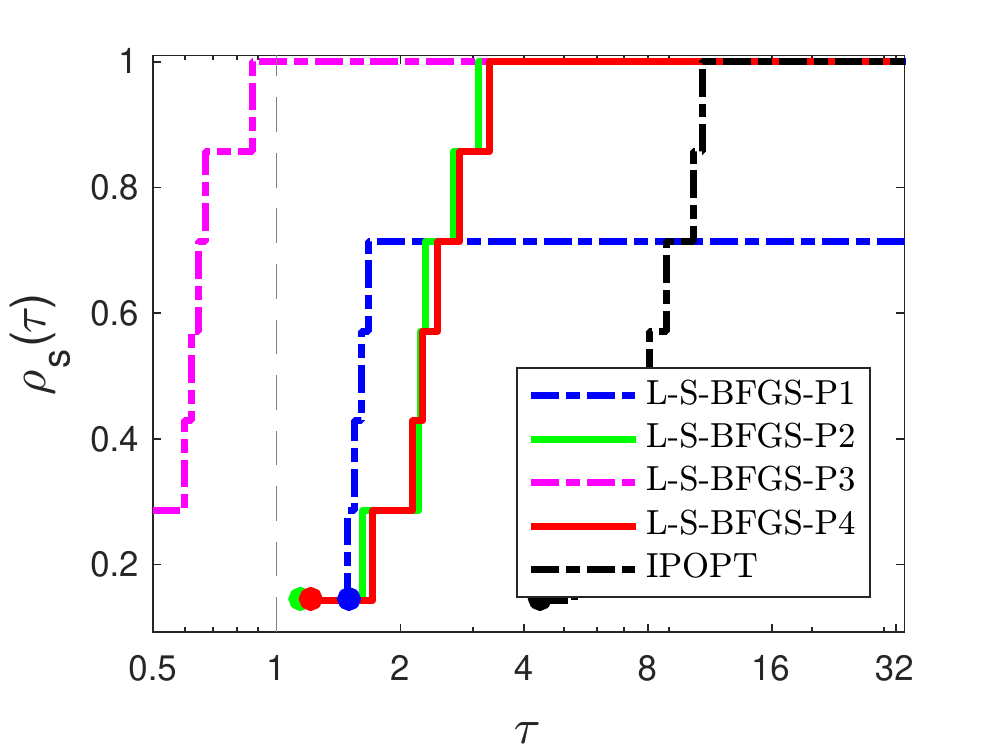}
	\end{minipage}
		\hfill
	\begin{minipage}{0.48\textwidth}
		\includegraphics[trim=0 0 20 15,clip,width=\textwidth]{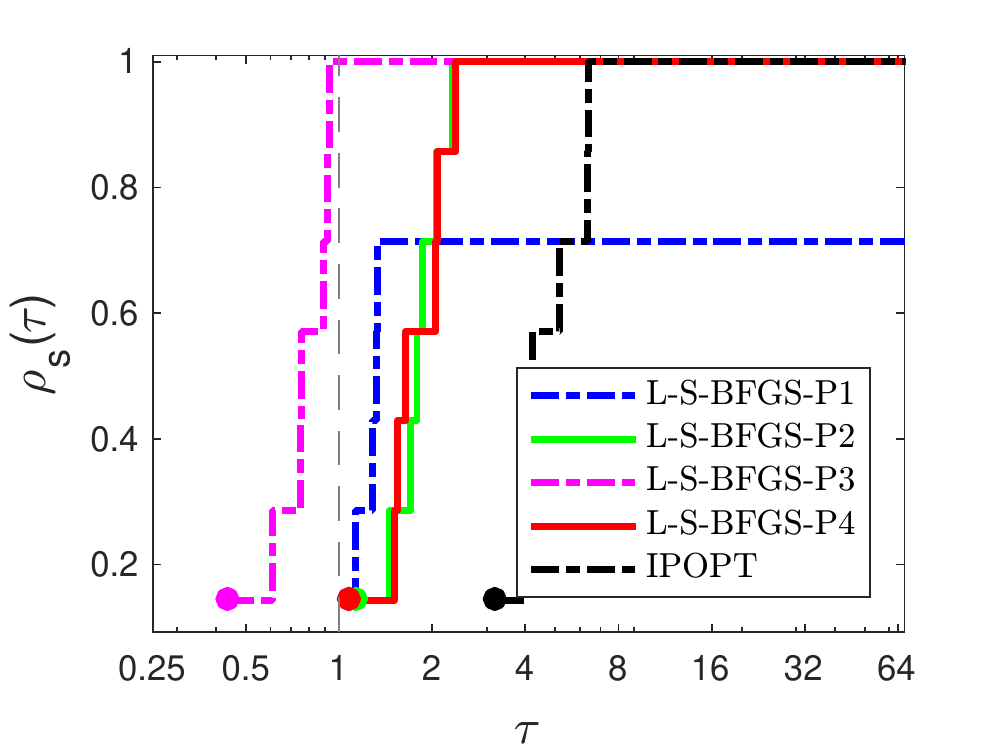}
	\end{minipage}
	\caption{Comparison of L-S-BFGS-P on structured objective functions to IPOPT and L-BFGS-B. Left: number of iterations, right: time.}
		\label{fig:EX_IV_SM_P}       
\end{figure*}

In this experiment the structured solvers are compared to 
IPOPT \cite{WaechterBiegler06} with an L-BFGS quasi-Newton matrix % L-BFGS and 
(we use a precompiled Mex file with IPOPT 3.12.12 that includes MUMPS and MA57 libraries). %{Accessed 21 June 2018}.  \cite{ZhuByrdNocedal97}  
The objective function is a structured quartic function 
\begin{equation}
\label{eq:struquartic}
	f(\b{x}) = \widehat{k}(\b{x}) + \widehat{u}(\b{x}), \quad \widehat{k}(\b{x}) =   \frac{1}{12} \sum_{i=1}^{n} (a^2_i x^4_i + 12 x_i g_i ), \quad
	\widehat{u}(\b{x}) = \frac{1}{2} \sum_{i=1}^n q_i x^2_i,
\end{equation}
where the data $ a_i, g_i $ and $ q_i $ are random normal variables with $ n = j \times 100, 1 \le j \le 7  $.  The
starting values are all ones, i.e., $ \bz{x} = \b{1} $. 
%As in \cref{subsec:EX_I}, the problem objective functions are structured quadratics defined by
%the expressions \cref{eq:struc_quad}. 
We specify the limited memory BFGS option for IPOPT using the setting 
\texttt{hessian\_approximation= `limited memory'} and tolerances by 
\texttt{tol=9.5e-10} and 
\texttt{acceptable\_tol = 9.5e-10}. For all solvers we set $ m = 8 $ (memory parameter),
and maximum iterations to 10,000. A solver is regarded to have converged when $ \|\bk{g}\|_{\infty} \le 9.5\times 10^{-5} $. 
The average outcomes of 5 runs of the experiments are in Figure \ref{fig:EX_IV_SM_P}. %\cref{fig:EX_IV_SM_M,fig:EX_IV_SM_P}.
%The L-BFGS solver uses the same strong Wolfe line search as  the algorithms from \cref{sec:alg}.
%\begin{figure*}[t!]
%	\begin{minipage}{0.48\textwidth}
%		\includegraphics[width=\textwidth]{figs/EXPERIMENT_IV_SM_ITER}
%	\end{minipage}
%		\hfill
%	\begin{minipage}{0.48\textwidth}
%		\includegraphics[width=\textwidth]{figs/EXPERIMENT_IV_SM_TIME}
%	\end{minipage}
%	\caption{Comparison of L-S-BFGS-M on structured objective functions with IPOPT and L-BFGS-B. 
%	Left: number of iterations, right: time.}
%		\label{fig:EX_IV_SM_M}       
%\end{figure*}
%
%In \cref{fig:EX_IV_SM_M} the proposed structured methods do well in number of iterations, where
%L-S-BFGS-M4 (red) requires the fewest iterations overall. L-BFGS-B (yellow)
%requires the least amount of time. Note that L-BFGS-B is based on low-level languages
%such as Fortran, which are inherently faster than MATLAB. IPOPT appears to require more iterations
%and time in relation to the other solvers. \jjb{In this regard, note that IPOPT is intended as a general purpose nonlinear programming
%solver that also handles general constraints. Therefore IPOPT may make additional computations to solve the problems.}
%Experiment on large clustered eigenvalues.

IPOPT and the L-S-BFGS-P solvers converge to the specified tolerances on all problems. The outcomes of the number of iterations (left plot) and computational times (right plot) in Figure \ref{fig:EX_IV_SM_P} are consistent. % L-BFGS converges on about 28$\%$  of the problems. 
In particular, we observe that the differences in the number of iterations are roughly reflected in the difference in the computational
times. In this problem the known Hessian is sensitive to changes in $ \b{x} $, and including second-order information in the quasi-Newton approximations yields outcomes with fewer iterations.

%objects
\subsection{Experiment V}
\label{subsec:EX_V}
\begin{figure*}[t!]	
		\includegraphics[width=\textwidth]{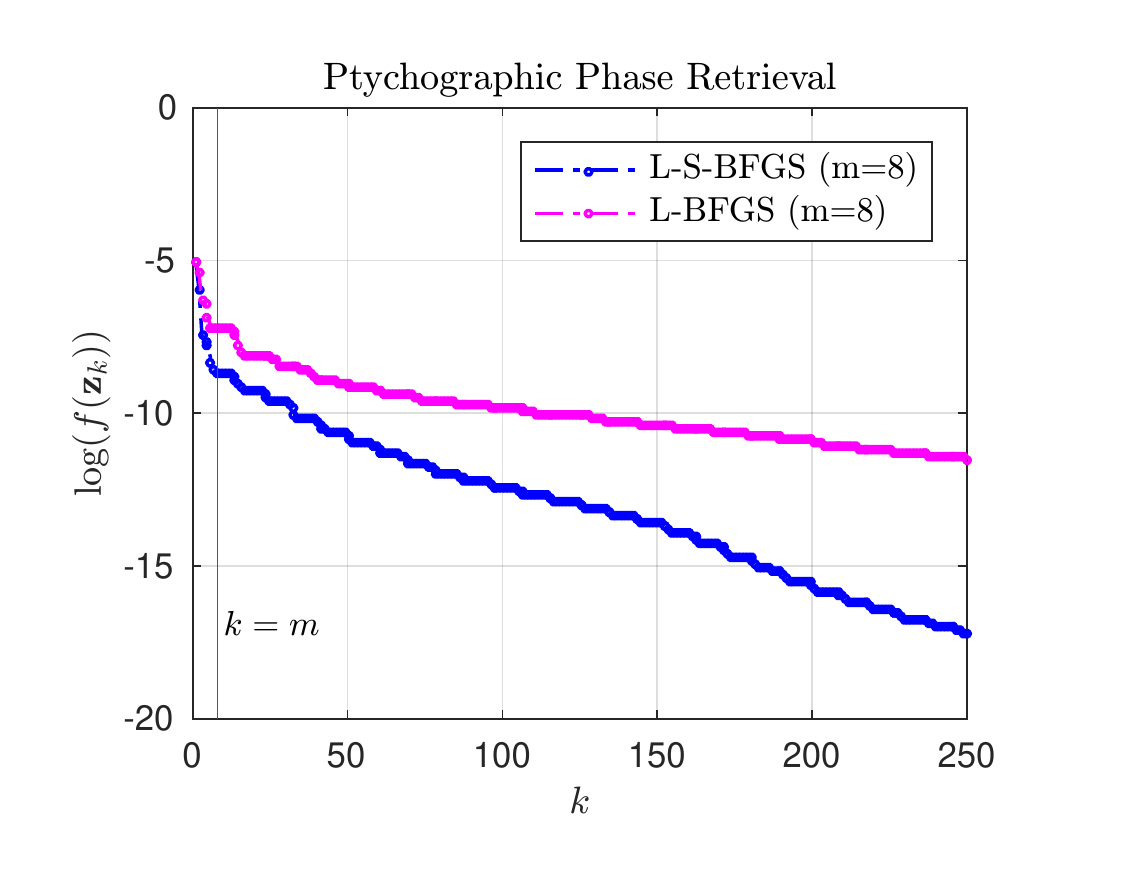}
		\caption{\jjbc{Comparison of log objective function values on a Ptychographic Phase Retrieval problem for the 
		first 250 iterations of 2 L-BFGS solvers. The L-S-BFGS method, using the known derivatives in $ \b{K} $, converges faster
		than the classical L-BFGS method with an identity initialization. The computational cost for the search directions in this problem scales according to $n$ with limited memory.
		}}
		\label{fig:EX_V}       
\end{figure*}

\jjbc{This experiment describes the application of the structured compact BFGS methods on an imaging problem. 
Large-scale imaging problems are challenging, because they involve large amounts of 
data and high-dimensional parameter space. Typically, image reconstruction problems are formulated
as optimization problems. In \cite{WuDi20}, efficient gradient-based quasi-Newton techniques for large-scale
ptychographic phase retrieval are described. However, even if the objective function is not directly
formulated as in problem \eqref{eq:intro_strucmin} it may still be possible to exploit  known 2$^{\text{nd}}$ derivatives. 
%---------------------- Updated problem formulation ----------------------%
Let $\bfz=\bfx+\bfy i \in \bbC^{n^2}$ be the object of interest, and $\bfd_j \in \bbR^{m^2}$ be the observed data (or intensities) measured from the $j^{th}$ probe, where $n^2$ and $m^2$ are the dimensions of the vectorized object and data resolution images, respectively. A ptychography experiment is modeled by
\begin{align} \label{eq:forwardProblem}
    \bfd_j = |\mathcal{F}(\bfQ_j\bfz)|^2 + \bfepsilon_j, \quad j=1,\ldots,N,
\end{align}
where $N$ is the total number of probes (or scanning positions), $\mathcal{F}\colon \mathbb{C}^{m^2} \mapsto \mathbb{C}^{m^2}$ is the two-dimensional discrete Fourier operator, $\bfQ_j \in \bbC^{m^2 \times n^2}$ is the $k^{th}$ probe (a diagonal \textit{illumination} matrix), and $\bfepsilon_j \in \bbR^{m^2}$ is the noise corresponding to the $k^{th}$ measurement error. The diagonal elements of $\bfQ_j$ are nonzero in the columns corresponding to the pixels being illuminated in the object at scanning step $j$.
There are different ways for formulating the reconstruction problem. One such formulation is the amplitude-based error metric
\begin{align}\label{eq:misfit1}
  \ds \underset{\b{z}}{ \textnormal{ minimize } } f(\bfz) = \frac{1}{2}\hf \sum_{j=1}^N \left\| \left|  \mathcal{F}(\bfQ_j \bfz)\right| - \sqrt{\bfd_j} \right\|_2^2 = \frac{1}{2}\hf \sum_{j=1}^N r_j^T r_j,
\end{align}
where $r_j=|\mathcal{F}(\bfQ_j \bfz)| - \sqrt{\bfd_j}$. Let $d_j=\sqrt{\bfd_j}$. Here, $ f \colon \mathbb{C}^{n^2} \mapsto \bbR$ is a real-valued cost function defined on the complex domain, and is therefore not complex-differentiable \cite{remmert2012theory}. To overcome the lack of complex-differentiability, it is common to employ the notion of $\bbC\bbR$ (Wirtinger) Calculus, where the derivatives of the real and imaginary parts of $\bfz$ are computed independently
\cite{remmert2012theory,sorber2012unconstrained}. For these real-valued functions, the mere existence of these Wirtinger derivatives is necessary and sufficient for the existence of a stationary point \cite{brandwood1983complex,remmert2012theory,sorber2012unconstrained}. Using Wirtinger calculus (using $z_j=\Fcurly \bfQ_j\bfz$), the partial gradients for~\eqref{eq:misfit1} can be computed as % , let $z_j=\Fcurly \bfQ_j\bfz$, 
\begin{equation}
  \begin{split}
    \nabla_\bfz r_j=J_j=\hf \overline{\diag{z_j/\left| z_j \right|}}\Fcurly \bfQ_j, \\
        \nabla_\bfz f =\sum_{j=1}^N J_j^*r_j=\sum_{j=1}^N \hf \bfQ_j^*\Fcurly^*\diag{z_j/\left| z_j \right|}\left(\left| z_j \right|-d_j\right),
      \end{split}
\end{equation}
\paragraph{Hessian.}
To compute the Hessian matrix, let \[T_{1,j}=\bfQ_j^*\Fcurly^*\diag{d_j/\left| z_j \right|}\Fcurly\bfQ_j,\] and
\[T_{2,j}=\bfQ_j^*\Fcurly^*\diag{d_j\odot z_j^2/\left| z_j \right|^3}\overline{\Fcurly}\overline{\bfQ_j},\] then
\begin{equation}
\renewcommand{\arraystretch}{1.6}
\ds
\b{H}=\sum_{j=1}^N
\begin{bmatrix}
\ds
\bfQ_j^*\bfQ_j-\hf \Refr{T_{1,j}}+\hf \Refr{T_{2,j}} & \hf \Im{T_{1,j}}+\hf\Im{T_{2,j}}\\
\hf \Im{T_{1,j}^*}+\hf\Im{T_{2,j}^*} & \bfQ_j^*\bfQ_j-\hf \Refr{T_{1,j}}-\hf \Refr{T_{2,j}}
\end{bmatrix},
\end{equation}
where the known 2$^{\text{nd}}$ derivatives are 
$
\b{K} = \sum_{j=1}^N
\bigg[
			\begin{smallmatrix} 
				\b{Q}^*_j \b{Q}_j 	& \\
							& \b{Q}^*_j \b{Q}_j
			\end{smallmatrix} 
\bigg]
$ and the remaining block elements of $ \b{H} $ are estimated.}
\jjbc{Defining the vectorization from complex to % $ \b{x} = \text{vec}(\text{Real}(\b{z}), \text{Img}(\b{z})) \equiv \text{vecR}(\b{z}) $
real variables by $ \b{x} =\text{vecR}(\b{z}) \equiv \text{vec}(\Refr{\b{z}}, \Im{\b{z}}) $, where
$ \text{vec}(\biidx{1}{x},\biidx{2}{x}) = [ \: \biidx{1}{x}^T \: \biidx{2}{x}^T \: ]^T $, we define the vectors for the structured BFGS methods by $ \bk{y} = \text{vecR}(\nabla f(\bko{z})) - 
\text{vecR}(\nabla f(\bk{z}))$, 
$ \bk{s}  = \text{vecR}(\bko{z}) - \text{vecR}(\bk{z}) $, $ \bhk{u} 	= \bk{y} - \b{K}\bk{s} $ and
\begin{equation*}
%	\bar{\b{u}}_k 	= \bk{y} - \b{K}\bk{s}, \quad \quad
	\bk{u} 		= \bhk{u} + \b{K} \bk{s} = \bk{y}.
\end{equation*}
Using these vectors, we can form the compact structured BFGS matrices. % using the known derivatives in $ \b{K} $.
In this experiment, we compare a limited memory structured BFGS method (L-S-BFGS) and limited memory BFGS % the the full memory structured BFGS (S-BFGS),
%a limited memory version (L-S-BFGS), the full memory BFGS
(L-BFGS) method in Figure \ref{fig:EX_V}. The image dimensions are $ \widehat{n} = 50 $ so
that the total number of real variables is $ n = 2\cdot \widehat{n}^2 = 5,000  $. Moreover, $ \widehat{m} = 16 $ so
that $ \widehat{m}^2 \times \widehat{n}^2 = 256 \times 5,000 = 1,280,000 $, and $ N = 16 $. Because of the structure of 
the known Hessian $l=1$ (cf. Table 1), and solves with this matrix are done on the order of $ \mathcal{O}(n) $. \jjbc{Because of the
size of this problem, and the corresponding computational/memory requirements the recursive update formulas from eqs. \eqref{eq:intro_sbfgsm} and \eqref{eq:intro_sbfgsp} are not applicable, yet the limited memory techniques threshold the required computational resources.}} 
 
%\begin{align*}
%	\bk{y} 		&= \nabla f(\bko{z}) - \nabla f(\bk{z}), \\
%	\bk{s} 		&= \bko{z} - \bk{z}, \\
%	\bar{\bk{u}} 	&= \bk{y} - \b{K}\bk{s}, \\
%	\bk{u} 		&= \bk{y} + \b{K} \bk{s} = \bk{y}.
%\end{align*}}
%However, when additional $2^{\text{nd}}$ derivatives
%of the objective function can efficiently be computed, it may be desirable to exploit this information
%in the optimization method. Even if the problem}

\section{Conclusions}
\label{sec:concl}
In this article we develop the compact representations of the structured
BFGS formulas proposed in Petra et al. \cite{PetraChiangAnitescu19}.
Limited memory versions of the compact representations with four 
non-constant initialization strategies are implemented in two line search algorithms.
\jjb{The proposed limited memory compact representations enable efficient
search direction computations by the Sherman-Morrison-Woodbury formula and
the use of efficient initialization strategies.} The proposed methods are compared in a collection of experiments, which 
include the original full-memory methods. The structured methods typically
require fewer total iterations than do the unstructured approaches. Among the four
proposed initialization strategies, initializations 1 and 2 appear best for the
structured minus methods (L-S-BFGS-M), whereas initializations
4 and 2 appear robust for the structured plus (L-S-BFGS-P) methods. \jjbc{In
an array of applications, including a large-scale real world imaging problem, the proposed
structured limited memory methods obtain better numerical results than conventional
unstructured methods.}

%\section*{Acknowledgments}
%\section*{Acknowledgments}
\begin{acknowledgements}
This work was supported by the U.S. Department of Energy, Office of Science,
Advanced Scientific Computing Research, under Contract DE-AC02-06CH11357
at Argonne National Laboratory.
through the Project "Multifaceted Mathematics for Complex Energy Systems."
This work was also performed under the auspices of the U.S. Department of Energy by 
Lawrence Livermore National Laboratory under Contract DE-AC52-07NA27344. 
\end{acknowledgements}

\section*{Appendix A: Initialization Comparison with $ \phi = 1000 $}
\jjbc{In Section \ref{subsec:EX_I}, the four L-S-BFGS-M initializations were compared on structured quadratic objective functions with eigenvalues
clustered around 1, whereas in this section the eigenvalues are clustered around 1000. In particular, $ \phi = 1000 $, and the elements of $ \b{D}_i $ are uniformly distributed  in the interval $[-999,0]$.
The results are displayed in Figure \ref{fig:EX_IA_LRG}.
\begin{figure*}[t!]
	\begin{minipage}{0.48\textwidth}
		\includegraphics[trim=0 0 20 15,clip,width=\textwidth]{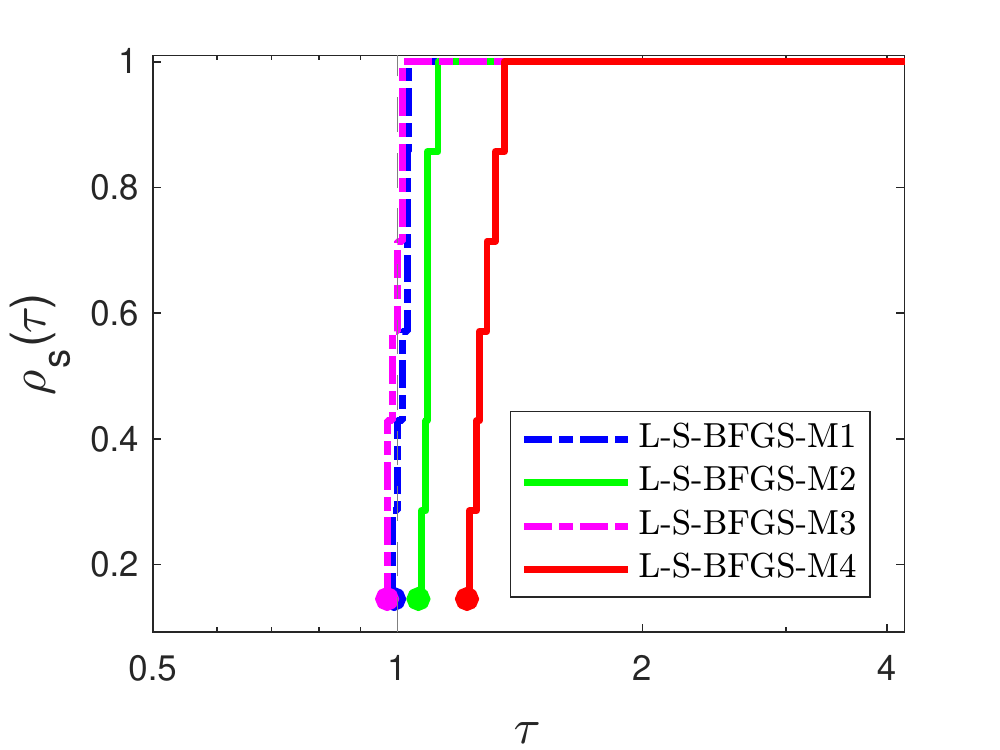}
	\end{minipage}
		\hfill
	\begin{minipage}{0.48\textwidth}
		\includegraphics[trim=0 0 20 15,clip,width=\textwidth]{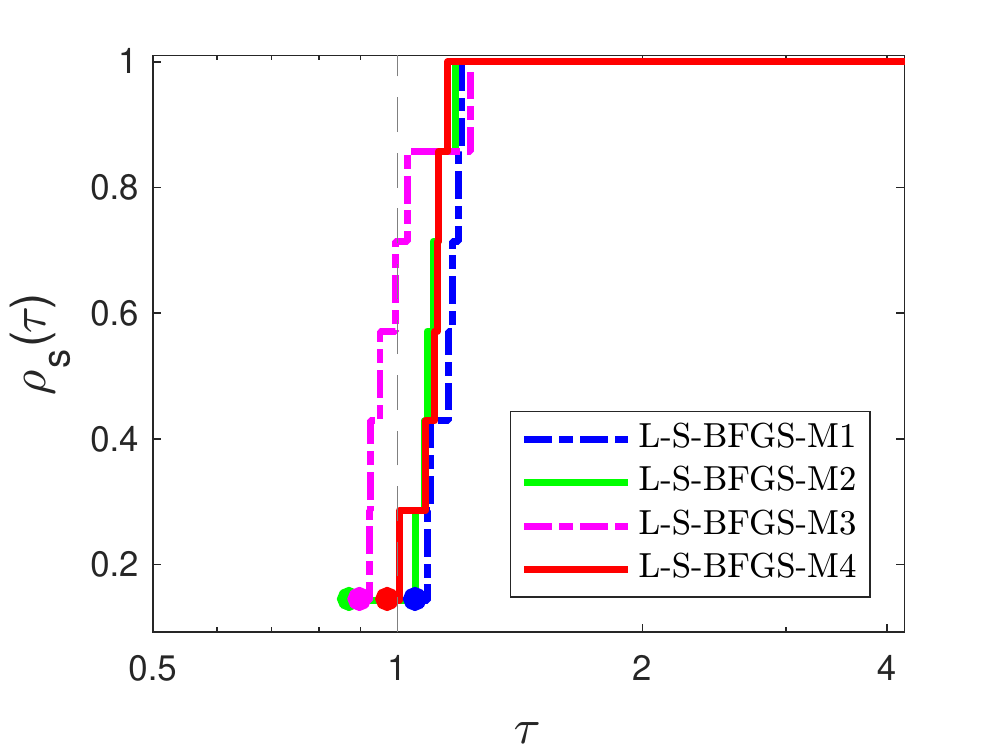}
	\end{minipage}
	\caption{Comparison of initialization strategies for L-S-BFGS-M on problems with eigenvalues clustered around 1,000 with 
	$1 \le \lambda_r \le 1000 $ and $ \lambda_{r+1} = \cdots = \lambda_{n} = 1000 $. Left: number of iterations; right: time.}
		\label{fig:EX_IA_LRG}       
\end{figure*}
For the large clustered eigenvalues Init. 1 and 3 require the fewest iterations, while
Init. 3 appears fastest overall. For L-S-BFGS-P the computations with $ \phi = 1000 $
are in Figure \ref{fig:EX_IA_LRG_P}
\begin{figure*}[t!]
	\begin{minipage}{0.48\textwidth}
		\includegraphics[trim=0 0 20 15,clip,width=\textwidth]{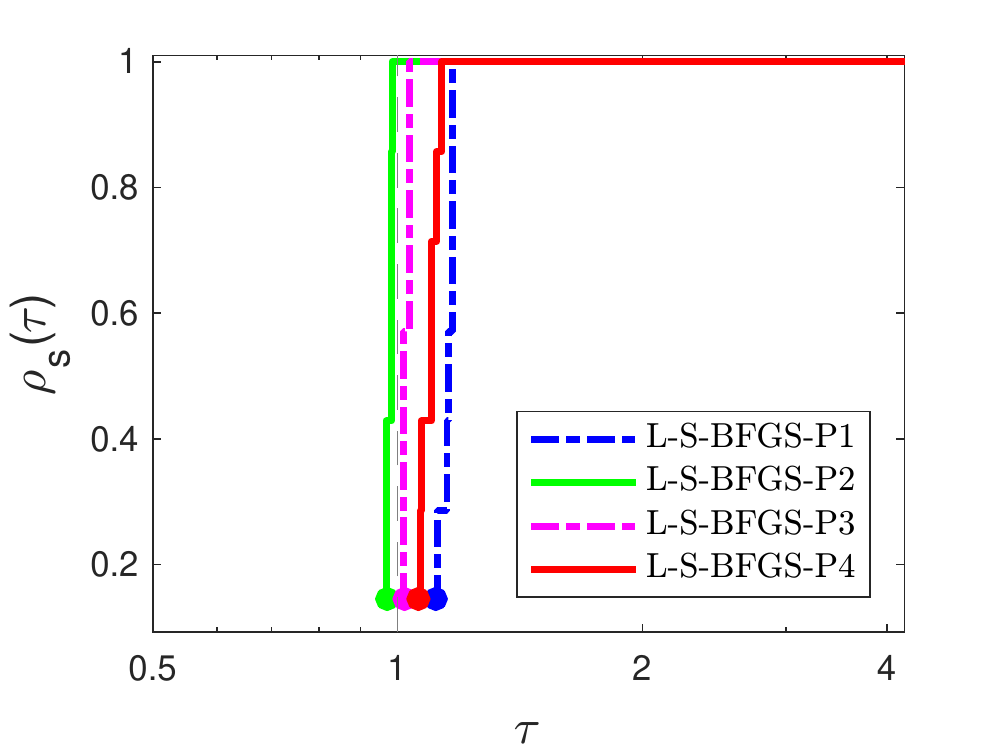}
	\end{minipage}
		\hfill
	\begin{minipage}{0.48\textwidth}
		\includegraphics[trim=0 0 20 15,clip,width=\textwidth]{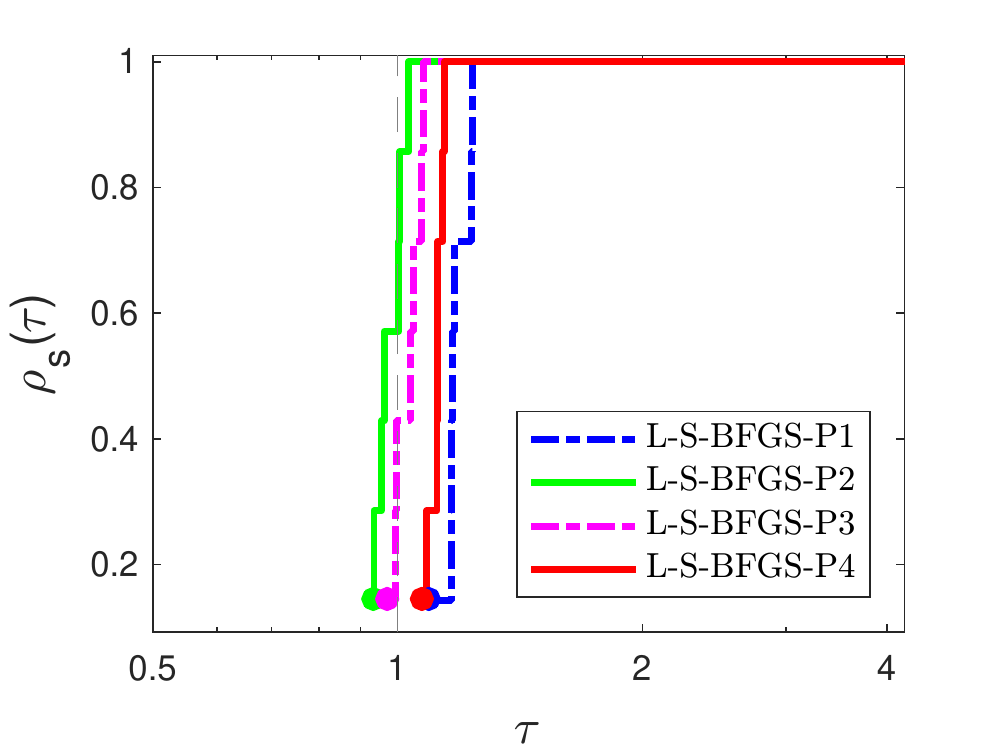}
	\end{minipage}
	\caption{Comparison of initialization strategies for L-S-BFGS-P on problems with eigenvalues clustered around 1,000 with 
	$1 \le \lambda_r \le 1000 $ and $ \lambda_{r+1} = \cdots = \lambda_{n} = 1000 $. Left: number of iterations; right: time.}
		\label{fig:EX_IA_LRG_P}       
\end{figure*}
In the comparison of L-S-BFGS-P, Init. 2 and Init. 3 do best in iterations.}

\bibliographystyle{spmpsci}      % mathematics and physical sciences
\bibliography{myrefs}   % name your BibTeX data base

%\begin{thebibliography}{}
%%
%% and use \bibitem to create references. Consult the Instructions
%% for authors for reference list style.
%\end{thebibliography}

\end{document}

%% file: cover.tex
\pagestyle{empty}
  
\vspace{1.75in}

\begin{centering}

ARGONNE NATIONAL LABORATORY

9700 South Cass Avenue

Argonne, Illinois  60439

\vspace{1.5in}

{\large \textbf{Compact Representations of Structured BFGS Matrices}}

\vspace{.5in}

\textbf{J. J. Brust, W. Di, S. Leyffer, and C. G. Petra}

\vspace{.5in}

%% Argonne Leadership Facility and
Mathematics and Computer Science Division

\vspace{.25in}

Preprint ANL/MCS-P9279-0120

\vspace{.5in}

August 2020

\end{centering}

\vspace{2.0in}

\bigskip

\par\noindent
\footnotetext [1]
{
This work was supported by the U.S. Department of Energy, Office of Science,
Advanced Scientific Computing Research, under Contract DE-AC02-06CH11357
at Argonne National Laboratory.
through the Project "Multifaceted Mathematics for Complex Energy Systems."
This work was also performed under the auspices of the U.S. Department of Energy by 
Lawrence Livermore National Laboratory under Contract DE-AC52-07NA27344. 
}

\newpage

\vspace*{\fill}
\begin{center}
\fbox{
\parbox{4in}{
The submitted manuscript has been created by UChicago Argonne, LLC, Operator of Argonne 
National Laboratory (``Argonne''). Argonne, a U.S. Department of Energy Office of Science 
laboratory, is operated under Contract No. DE-AC02-06CH11357. The U.S. Government retains 
for itself, and others acting on its behalf, a paid-up nonexclusive, irrevocable worldwide 
license in said article to reproduce, prepare derivative works, distribute copies to the 
public, and perform publicly and display publicly, by or on behalf of the Government. 
The Department of Energy will provide public access to these results of federally 
sponsored research in accordance with the DOE Public Access 
Plan. \texttt{http://energy.gov/downloads/doe-public-accessplan}
}}
\end{center}
\vfill

\newpage
\pagestyle{plain}
\setcounter{page}{1}